\documentclass[11pt,reqno]{amsart}
\usepackage[T1]{fontenc}
\usepackage[utf8]{inputenc}
\usepackage{lmodern}
\usepackage{amsmath,amssymb,amsthm,mathtools,mathrsfs}
\usepackage[expansion=false]{microtype}
\usepackage[a4paper,margin=28mm]{geometry}
\usepackage{xcolor}
\usepackage{enumitem}
\usepackage{booktabs}
\usepackage[hidelinks]{hyperref}
\usepackage{zref-clever}
\zcsetup{cap=true,nameinlink=true}
\zcRefTypeSetup{item}{Name-sg={},name-sg={},Name-pl={},name-pl={},namesep={},
  refbounds={(,,,)},rangesep={\textendash}}

\setlist[enumerate]{itemsep=2pt,topsep=4pt}

\providecommand{\germ}{\mathfrak}

\numberwithin{equation}{section}
\theoremstyle{plain}
\AddToHook{env/proposition/begin}{\zcsetup{countertype={theorem=proposition}}}
\AddToHook{env/lemma/begin}{\zcsetup{countertype={theorem=lemma}}}
\AddToHook{env/corollary/begin}{\zcsetup{countertype={theorem=corollary}}}
\AddToHook{env/definition/begin}{\zcsetup{countertype={theorem=definition}}}
\AddToHook{env/remark/begin}{\zcsetup{countertype={theorem=remark}}}
\AddToHook{env/example/begin}{\zcsetup{countertype={theorem=example}}}
\AddToHook{env/convention/begin}{\zcsetup{countertype={theorem=convention}}}
\AddToHook{env/conjecture/begin}{\zcsetup{countertype={theorem=conjecture}}}
\zcRefTypeSetup{conjecture}{Name-sg={Conjecture},name-sg={conjecture},Name-pl={Conjectures},name-pl={conjectures}}
\zcRefTypeSetup{convention}{Name-sg={Convention},name-sg={convention},Name-pl={Conventions},name-pl={conventions}}
\newtheorem{theorem}{Theorem}[section]
\newtheorem{proposition}[theorem]{Proposition}
\newtheorem{lemma}[theorem]{Lemma}
\newtheorem{corollary}[theorem]{Corollary}
\newtheorem{conjecture}[theorem]{Conjecture}
\theoremstyle{definition}
\newtheorem{definition}[theorem]{Definition}
\newtheorem{convention}[theorem]{Convention}
\theoremstyle{remark}
\newtheorem{remark}[theorem]{Remark}
\newtheorem{example}[theorem]{Example}
\theoremstyle{plain}
\newtheorem{introtheorem}{Theorem}

\zcsetup{countertype={introtheorem=theorem}}

\newcommand{\Z}{\mathbb{Z}}      % the integers
\newcommand{\Q}{\mathbb{Q}}      % the rationals
\newcommand{\C}{\mathbb{C}}      % the complex numbers
\newcommand{\Gm}{\mathbb{G}_m}   % the multiplicative group
\newcommand{\glf}{\Bbbk}         % a global field (\zcref{sec:global}); k is the residue field of F

\newcommand{\OO}{\mathcal{O}}    % ring of integers of F
\newcommand{\pp}{\mathfrak{p}}   % its maximal ideal
\newcommand{\Lat}{\mathcal{L}}   % lattice in \g
\newcommand{\dLat}{\Lambda}      % lattice in \gd (the dual of \Lat)
\newcommand{\Aff}{\mathbb{A}}    % affine space
\newcommand{\bG}{\mathbf{G}}     % the connected reductive group over F
\newcommand{\bA}{\mathbf{A}}     % a split torus
\newcommand{\bP}{\mathbf{P}}     % a parabolic subgroup of \bG
\newcommand{\bQ}{\mathbf{Q}}     % another parabolic subgroup of \bG
\newcommand{\bM}{\mathbf{M}}     % a Levi subgroup
\newcommand{\bL}{\mathbf{L}}     % another Levi subgroup
\newcommand{\bU}{\mathbf{U}}     % a unipotent radical
\newcommand{\bT}{\mathbf{T}}     % a maximal torus
\newcommand{\bZ}{\mathbf{Z}}     % the center of \bG; Z=\bZ(F)
\newcommand{\bH}{\mathbf{H}}     % a group of type B in characteristic 2 (\zcref{lem:oddorthogonal})
\newcommand{\Godd}{\bG_{\mathrm{odd}}}
\newcommand{\Gord}{\bG_{\mathrm{ord}}}
\newcommand{\godd}{\g_{\mathrm{odd}}}
\newcommand{\gord}{\g_{\mathrm{ord}}}
\newcommand{\g}{\mathfrak{g}}    % the Lie algebra of G
\newcommand{\gd}{\mathfrak{g}^*} % its dual
\newcommand{\uu}{\mathfrak{u}}   % Lie algebra of a unipotent radical
\newcommand{\m}{\mathfrak{m}}    % Lie algebra of a Levi subgroup
\newcommand{\liet}{\mathfrak{t}}  % Lie algebra of a torus
\newcommand{\lieh}{\mathfrak{h}}  % Lie algebra of \bH
\newcommand{\Ell}[2][]{#2_{#1\mathrm{ell}}}   % use with the dual Lie algebra: \Ell{\gd}
\newcommand{\Par}[2][]{#2_{#1\mathrm{par}}}
\newcommand{\sclab}{\mathrm{sc}}
\newcommand{\dslab}{\mathrm{ds}}
\newcommand{\dens}[2]{\delta_{#1}(#2)}  % \dens{\ast}{K}
\newcommand{\scd}[1]{\dens{\sclab}{#1}}  % supercuspidal density
\newcommand{\dsd}[1]{\dens{\dslab}{#1}}  % discrete-series density
\newcommand{\fmass}[1]{f^{#1}}
\newcommand{\fsc}{\fmass{\sclab}}
\newcommand{\fds}{\fmass{\dslab}}
\newcommand{\fdeg}[1]{d_{#1}}  % formal degree of a discrete series #1, a Haar measure on G
\newcommand{\elld}[2][]{\delta_{\mathrm{ell}}^{#1}(#2)}  % optional argument: the group
\newcommand{\Bad}{\mathcal{B}}   % the boundary nets \Bad_h (\zcref{thm:sard})
\newcommand{\bnets}{boundary nets}

\newcommand{\Hk}{\mathcal{H}}    % Hecke algebra
\newcommand{\Pstd}{\mathscr{P}^{\mathrm{st}}}    % all proper standard parabolic subgroups (\zcref{sec:parabolics})
\newcommand{\Pssall}{\mathscr{P}^{\mathrm{ss}}}  % all proper semistandard parabolic subgroups
\newcommand{\Pc}{\mathscr{P}^{\mathrm{st}}_{\mathrm{max}}}    % the maximal elements of \Pstd
\newcommand{\Pcss}{\mathscr{P}^{\mathrm{ss}}_{\mathrm{max}}}  % the maximal elements of \Pssall
\newcommand{\Nss}{N_{\bG}}  % the number |\Pcss| of maximal semistandard parabolics (\zcref{sec:parabolics})
\newcommand{\Lcl}{\mathscr{L}}   % the G-conjugacy classes of Levi subgroups (\zcref{sec:parabolics})
\newcommand{\Lsup}[1]{\mathfrak{l}_{#1}}
\newcommand{\Lform}[1]{\mathfrak{l}(#1)}
\newcommand{\Lnc}[1]{\mathfrak{l}_{\mathrm{nc}}(#1)}
\newcommand{\GhatL}[2]{\widehat{G}_{#1}#2}
\newcommand{\densL}[3]{\delta_{#1}^{#2}(#3)}
\newcommand{\wlab}{\mathrm{w}}   % label of the weak Jacquet module r_P^{\wlab}
\newcommand{\Pim}[1]{\gd_{#1}}
\newcommand{\Estr}{\mathcal{E}}
\newcommand{\levid}[2]{\delta^{#1}(#2)}
\newcommand{\Crvmax}{\Crv^{\max}}
\newcommand{\Badmax}{\Bad^{\max}}
\newcommand{\eexpall}{\tau_0}
\newcommand{\dGall}{d^{\,0}_{\bG}}
\newcommand{\weakjac}{weak Jacquet module}

\newcommand{\levistratum}{Levi stratum}
\newcommand{\levistrata}{Levi strata}

\newcommand{\cchain}{contracting chain}
\newcommand{\Cchain}{Contracting chain}
\newcommand{\Fl}{\mathcal{F}}    % flag manifold \Fl_P=\bFl_P(F)
\newcommand{\bFl}{\boldsymbol{\Fl}}  % flag variety \bFl_P=\bG/\bP
\newcommand{\sS}{\mathfrak{s}^+}   % enlarged positive lattice
\newcommand{\sT}{\mathfrak{s}^-}   % enlarged negative lattice
\newcommand{\Ghat}[1]{\widehat{G}_{#1}}    % \Ghat{\ast}: the class of type \ast in the unitary dual
\newcommand{\Gsc}{\Ghat{\sclab}}    % supercuspidal unitary dual
\newcommand{\Gds}{\Ghat{\dslab}}    % discrete series modulo the center
\newcommand{\Gtemp}{\widehat{G}_{\mathrm{temp}}} % tempered unitary dual
\newcommand{\Pps}{\delta_{\mathrm{ps}}(K_n)}       % principal series density (\zcref{sec:gl2exact})
\newcommand{\Pst}{\delta_{\mathrm{St}}(K_n)}
\newcommand{\Cps}{c_{\mathrm{ps}}}
\newcommand{\Cst}{c_{\mathrm{St}}}
\newcommand{\St}{\mathrm{St}}                   % Steinberg representation
\newcommand{\Iw}{\mathcal{I}}                   % Iwahori subgroup of GL_2
\newcommand{\Xc}{\mathcal{X}}                   % unitary characters of F^\times
\newcommand{\cG}{\mathcal{G}}                   % smooth (in the unramified case, \zcref{sec:integraldata}: reductive) group scheme over O
\newcommand{\cP}{\mathcal{P}}                   % parabolic subgroup scheme
\newcommand{\cM}{\mathcal{M}}                   % Levi subgroup scheme
\newcommand{\cU}{\mathcal{U}}                   % unipotent radical scheme
\newcommand{\cT}{\mathcal{T}}                   % maximal torus scheme
\newcommand{\cS}{\mathcal{S}}                   % maximal split torus scheme
\newcommand{\one}{\mathbf{1}}
\newcommand{\der}{\operatorname{der}}

\newcommand{\abs}[1]{\lvert#1\rvert}
\newcommand{\norm}[1]{\lVert#1\rVert}
\newcommand{\etacompat}{$\eta$-compatible}
\newcommand{\adapted}{adapted}
\newcommand{\Adapted}{Adapted}

\DeclareMathOperator{\vol}{vol}
\DeclareMathOperator{\tr}{tr}
\DeclareMathOperator{\diag}{diag}
\DeclareMathOperator{\Ad}{Ad}
\DeclareMathOperator{\ad}{ad}
\DeclareMathOperator{\Hom}{Hom}
\DeclareMathOperator{\val}{val}
\DeclareMathOperator{\Lie}{Lie}
\DeclareMathOperator{\GL}{GL}
\DeclareMathOperator{\PGL}{PGL}
\DeclareMathOperator{\SL}{SL}
\DeclareMathOperator{\Gal}{Gal}
\DeclareMathOperator{\SO}{SO}
\DeclareMathOperator{\Mat}{M}
\DeclareMathOperator{\Crit}{Crit}

\DeclareMathOperator{\rad}{rad}
\DeclareMathOperator{\id}{id}
\DeclareMathOperator{\chr}{char}   % characteristic of a field
\DeclareMathOperator{\Irr}{Irr}
\DeclareMathOperator{\End}{End}
\DeclareMathOperator{\disc}{disc}

\DeclareMathOperator{\cInd}{c-Ind}  % compact induction
\DeclareMathOperator{\cond}{cond}   % conductor exponent of a character of F^\times

\newcommand{\peculiar}{peculiar}
\newcommand{\Peculiar}{Peculiar}
\newcommand{\ordinary}{ordinary}
\newcommand{\Ordinary}{Ordinary}
\newcommand{\oddorth}{odd orthogonal}
\newcommand{\Oddorth}{Odd orthogonal}
\newcommand{\eexp}{\tau}
\newcommand{\dG}{d_{\bG}}
\newcommand{\erad}{R_{\mathrm{ell}}}
\newcommand{\lf}{\lambda}
\newcommand{\cell}{\omega}
\newcommand{\Bf}{\mathsf{B}}
\newcommand{\Avg}{\mathsf{P}}
\newcommand{\nchart}{normalized chart}
\newcommand{\npar}{normalized parametrization}  % the inverse of a normalized chart
\newcommand{\nf}{\zeta}
\newcommand{\Inc}{\mathcal{Y}}
\newcommand{\bInc}{\boldsymbol{\Inc}}
\newcommand{\incm}{\mathrm{pr}}
\newcommand{\Crv}{\Sigma}
\newcommand{\Vf}{Y}
\newcommand{\XL}{X_{\dLat}}
\newcommand{\Incc}[1]{\Inc^{\circ}_{#1}}  % big cell of \Inc_P (\zcref{sec:levi})
\newcommand{\iatlas}{integral atlas}

\newcommand{\Plprob}{\mathbb P}

\begin{document}

\title[Supercuspidal, discrete, and elliptic densities]
{Supercuspidal, discrete, and elliptic densities for reductive groups over non-archimedean local fields}

\author{Erez Lapid}
\address{Department of Mathematics, Weizmann Institute of Science, Rehovot 7610001, Israel}
\email{erez.lapid@weizmann.ac.il}

\author{Mark Shusterman}
\address{Department of Mathematics, Weizmann Institute of Science, Rehovot 7610001, Israel}
\email{mark.shusterman@weizmann.ac.il}

\date{November 13, 2026}

\begin{abstract}
Let $G$ be the group of points of a connected reductive group over a
non-archimedean local field $F$ of arbitrary characteristic, with residue
field of size $q$. For a principal congruence tower $K_n$ attached to a Lie
algebra lattice $\Lat$, the supercuspidal and discrete series proportions
of the Plancherel mass of the $K_n$-fixed vectors both converge, with error
$O(q^{-n/4})$, to the probability that a Haar-random form in the dual lattice $\dLat=\Lat^\vee$ is
elliptic. For $\GL_d$, the supercuspidal limit confirms a prediction of Howe.
This probability is positive, and we characterize when it is one.
More generally, outside a union of congruence cells of total mass $O(q^{-n/4})$,
each cell determines a common Levi class for the supercuspidal and
discrete series supports of every tempered representation containing its
character. The limiting mass of each class is the probability
that a Haar-random form in $\dLat$ lies in the corresponding stratum. The proof combines finite Fourier
analysis, quantitative estimates for parabolic incidence maps, and a
Jacquet-module argument that detects unitary exponents.
\end{abstract}

\maketitle

% Table of contents at the section level only; raise tocdepth to 2 to
% list subsections as well.
\setcounter{tocdepth}{1}
\tableofcontents

\section{Introduction}\label{sec:intro}

Let $\bG$ be a connected reductive group over a non-archimedean local
field $F$ of arbitrary characteristic, with ring of integers $\OO$,
maximal ideal $\pp$ and residue field of size $q$. Write $G=\bG(F)$.
For a compact open subgroup $K$ of $G$, the Plancherel formula defines a
probability measure on the tempered dual:
\[
d\Plprob_K(\pi)=\vol(K)\dim\pi^K\,d\mu_G(\pi).
\]
Here $\mu_G$ is the Plancherel measure for the Haar measure on $G$
used to define $\vol(K)$.
We study this measure for a principal congruence tower $K_n$.
For $\GL_d$, Howe anticipated the limit of its supercuspidal mass
\cite[\S4]{MR0492088}. The limit is geometric, and the geometry also
determines the distribution of Levi classes for the two notions of support.

\subsection{The main result}

Let $\g=\Lie(\bG)(F)$ and let $\Lat\subset\g$ be a lattice, with dual lattice
$\dLat\subset\gd$. A normalized analytic chart $\kappa$ at the identity
defines $K_n=\kappa^{-1}(\pp^n\Lat)$ for $n$ sufficiently large.
For a smooth affine model $\cG$ over $\OO$, this includes the kernels of
reduction modulo $\pp^n$ with $\Lat=\Lie(\cG)(\OO)$.
A form $\lf\in\gd$ is \emph{elliptic} if $\lf|_{\uu_P}\ne0$ for every
proper parabolic $F$-subgroup $\bP$, where $\uu_P$ is the Lie algebra of
its unipotent radical. Write $\elld{\dLat}$ for the probability of
ellipticity under the probability Haar measure $\nu_{\dLat}$.

Let $\Lcl$ be the set of $G$-conjugacy classes of Levi subgroups, ordered
by inclusion up to conjugacy. Write $\Lsup{\sclab}(\pi)$ and
$\Lsup{\dslab}(\pi)$ for the Levi classes of the supercuspidal and discrete
series supports of a tempered representation $\pi$. On the geometric
side, consider the Levi classes of all parabolic subgroups $\bP$ such
that $\lf|_{\uu_P}=0$, including $\bG$. This set has a least element
$\Lform{\lf}$ by \zcref{lem:meet}. Write $\Estr_{[M]}$ for the set of
forms with $\Lform{\lf}=[M]$. In particular, $\Estr_{[G]}=\Ell{\gd}$.

\begin{introtheorem}[Plancherel mass and Levi classes]\label{thm:A}
For every $[M]\in\Lcl$,
\[
\Plprob_{K_n}\{\pi:\Lsup{\sclab}(\pi)=\Lsup{\dslab}(\pi)=[M]\}
=\nu_{\dLat}(\dLat\cap\Estr_{[M]})+O(q^{-n/4}).
\]
The total mass where the two support classes differ is $O(q^{-n/4})$.
In particular, the supercuspidal and discrete series masses both
converge to $\elld{\dLat}>0$ with error $O(q^{-n/4})$.
\end{introtheorem}

The precise statement in \zcref{thm:main} gives improved exponents when
the relevant incidence maps are generically noncritical. For the
supercuspidal and discrete series masses, only the maximal proper
parabolic subgroups are needed.
When the center is compact, the supercuspidal and discrete series
masses can also be written as the sums of
$\fdeg{\pi}(K_n)\dim\pi^{K_n}$ over the corresponding representations,
where $\fdeg{\pi}$ is the formal degree of $\pi$ viewed as a Haar
measure on $G$.
There are localized estimates for unions of congruence cells
(\zcref{cor:localized}) and explicit uniform bounds for reductive models
(\zcref{thm:uniform,thm:global}).

\subsection{The elliptic density}

For $\GL_2$, the lattice $\Mat_2(\OO)$ gives elliptic density $1/2$.
The congruence tower attached to the Iwahori Lie lattice has limiting
density $q/(q+1)$; the tower attached to its trace dual has limiting
density $1/(q+1)$ (\zcref{prop:gl2,prop:iwahori}).

In characteristic two, let $\Godd$ be the direct factor of $\bG$
comprising the geometric adjoint factors of type $\mathrm B_d$,
including $\mathrm B_1=\mathrm A_1$, and write
$\bG=\Godd\times\Gord$. In other characteristics put $\Gord=\bG$.

\begin{introtheorem}[Full elliptic density]\label{thm:B}
The elliptic density is one for some lattice, equivalently for every
lattice, if and only if $\Gord$ is anisotropic modulo its center.
\end{introtheorem}

The simplest isotropic example is split $\PGL_2$ in characteristic
two. The bracket of opposite root vectors is scalar and hence zero in
$\mathfrak{pgl}_2$. The incidence map of the Borel subgroup is everywhere
critical, and its image is null. The exact Plancherel calculation in
\zcref{prop:pgl2} exhibits the approach to full density.

\subsection{Idea of the proof}\label{sec:idea}

For $2r\le n$, the characters of $K_{n-r}/K_n$ are indexed by the
cosets of $\pp^r\dLat$ in $\dLat$, which we call \emph{cells}.
Finite Fourier analysis decomposes $\Plprob_{K_n}$ into contributions
from these cells, with each cell carrying its Haar measure
(\zcref{sec:cells}).

For a parabolic subgroup $\bP$, the incidence map
$G\times^P\uu_P^\perp\to\gd$ sends $[g,\nf]$ to $\Ad^*(g)\nf$.
Its image records which forms vanish on the Lie algebra of the
unipotent radical of a conjugate of $\bP$. The two manifolds have equal dimension. A quantitative
Sard estimate gives exceptional sets of measure $O(q^{-h})$ made up of
cosets of $\pp^{2h}\dLat$. Outside them, each such coset either misses
the image or admits lifts with the inverse bracket pairing
$Q_{\nf,P}^{-1}$ bounded by $O(q^h)$. In the big cell the Jacobian
determinant is a polynomial of degree $\dim\uu_P$, giving the improved
rate when it is nonzero (\zcref{sec:incidence}).

If a cell misses the incidence image, every representation containing
its character has zero Jacquet module along $\bP$. If it has a lift
with the stated inverse bound, an averaging argument gives nonzero
Jacquet modules. Once the averaging construction stabilizes, its
normalized operator is unitary. On the image of the isotypic
subspace in the normalized Jacquet module, it induces the action of
a central element of the Levi subgroup. For tempered representations,
this image therefore lies in the weak Jacquet module (\zcref{prop:elliptic,thm:regular}).
Applying both alternatives to all proper parabolic subgroups recovers
both support classes from the same cell. The argument uses
$r\le n/2+O(1)$ and $r\ge2h+O(1)$, so we can take $h=n/4+O(1)$.

\subsection{Relation to earlier work}\label{sec:related}

For $\GL_d$, the cells are Howe's dual blobs of shallow characters.
Howe describes the representations containing sufficiently regular
shallow characters and relates their Plancherel mass to the geometry of
the blobs \cite[\S\S2--4]{MR0492088}. The geometric and spectral
comparison here is quantitative and works in arbitrary characteristic,
including for wildly ramified groups. The argument uses congruence
charts in place of the exponential map and requires no local character
expansions.

Bushnell, Henniart and Kutzko compare the Plancherel measure of
representations containing a fixed type with that of its spherical
Hecke algebra \cite[Theorem~B]{MR2767510}. Kim obtains Plancherel
comparisons for families specified by good minimal $K$-types and by
Yu's data under hypotheses on the field and group
\cite[Theorem~3.3.1]{MR2075916}\cite[Theorem~10.1]{MR2276772}.
Here a single character of an abelian congruence quotient suffices to
recover the support classes on the cells outside the exceptional sets.
The Jacquet argument is related to the theory of covers
\cite[\S\S6--8]{MR1643417}. Its stationary step detects the unitary
exponents needed for discrete series support.

\subsection*{AI disclosure}
In preparing this paper the authors used the AI systems Claude (Anthropic)
and ChatGPT (OpenAI), for drafting and revising text, checking
computations and searching the literature. The authors have reviewed the
resulting material and take full responsibility for the content of the
paper.

\subsection*{Acknowledgments}
The work was supported in part by a grant of access to OpenAI models
through the ChatGPT for Academic Researchers Program.
The second named author is supported by the European Union (ERC, Function
Fields, 101161909) and is the Dr.~A.~Edward Friedmann Career Development
Chair in Mathematics.

\section{Preliminaries}\label{sec:notation}

Throughout, $F$ is a non-archimedean local field of arbitrary
characteristic, with ring of integers $\OO$, maximal ideal $\pp$, uniformizer $\varpi$,
residue field $k=\OO/\pp$ of cardinality $q$ and normalized valuation $\val$.

We freely use standard terminology and results on algebraic groups over
$F$ (which is imperfect when $\chr F>0$) from
\cite{MR0207712,MR1102012,MR3729270,MR3362817}, on reductive group schemes over
$\OO$ from \cite{MR3362641,MR4520154}, on smooth representations from \cite{MR0425030,MR0579172,MR2567785},
on the Plancherel formula from \cite{MR1989693}, and on analytic manifolds and groups over $F$ from \cite{MR2179691}; for standard groups and
their principal congruence subgroups see \cite[Chapter~III, \S7, nos.~3--4]{MR1728312} and
\cite[Chapter~13]{MR1720368}.

\subsection{Notation}\label{sec:notationbasic}

Let $\bG$ be a connected reductive group over $F$,
$G=\bG(F)$, $\g=\Lie(\bG)(F)$, $\gd=\Hom_F(\g,F)$ its dual and $D=\dim\bG=\dim_F\g$.
Let $\bZ$ be the center of $\bG$ and $Z=\bZ(F)$, the center of $G$.

Denote by $\Ad$ the adjoint action of $G$ on $\g$ and by $\Ad^*$ the coadjoint action on $\gd$, so that $\Ad^*(g)\lf=\lf\circ\Ad(g^{-1})$.

For any $\lf\in\gd$ let $\tilde\Bf_\lf$ be the alternating form on $\g$ given by
\[
\tilde\Bf_\lf(X,Y)=\lf([X,Y]).
\]
The corresponding linear map $\g\to\gd$, $Y\mapsto\tilde\Bf_\lf(\cdot,Y)$, is
$Y\mapsto\ad^*(Y)\lf$, where $\ad^*(Y)\lf=-\lf\circ\ad(Y)$ is the coadjoint
action of $\g$; it is the differential at $1$ of $g\mapsto\Ad^*(g)\lf$.

Let $\Hk(G)$ be the Hecke algebra of $G$, the algebra of smooth compactly
supported distributions on $G$ under convolution. 
We fix a Haar measure $dg$ on $G$. Then $f\mapsto f\,dg$ is a linear isomorphism
$C_c^\infty(G)\to\Hk(G)$; we write $\varphi/dg\in C_c^\infty(G)$ for the preimage
of $\varphi\in\Hk(G)$.
For a compact open subgroup $J$ of $G$, let $e_J\in\Hk(G)$ be the Haar
probability measure on $J$, extended by zero to $G$, and for a character
$\chi$ of $J$ let $e_J^\chi=\chi^{-1}\cdot e_J$. These are idempotents of $\Hk(G)$.

We fix once and for all a character $\psi$ of the additive group of $F$ which is trivial on $\OO$ and nontrivial on $\pp^{-1}$.

For a lattice $\Lat_0$ in a finite-dimensional $F$-vector space and a
vector $X$ let $\val_{\Lat_0}(X)=\sup\{a\in\Z:X\in\pp^a\Lat_0\}$ be
the depth of $X$, and $\norm{X}_{\Lat_0}=q^{-\val_{\Lat_0}(X)}$ its norm (the
maximum norm of its coordinates in any basis of $\Lat_0$); we drop
the subscript $\Lat_0$ when the lattice is clear from the context.
Linear maps between spaces with lattices are given the operator norms.
The level and depth parameters $n,r,h,b,a,s$ used below are integers.

For a subset $X$ of a topological space, $\partial X$ denotes its boundary.

By an \emph{analytic manifold} we mean a Hausdorff analytic manifold
over $F$, and a \emph{ball} in such a manifold of dimension $m$ is an open subset analytically
isomorphic to the unit ball $\{y\in F^m:\norm{y}\le1\}$ of $F^m$, for the maximum norm.
Only this elementary notion is used (charts given by convergent power series on
compact open balls); no rigid analytic geometry is needed.
For a morphism $f:X\to X'$ of analytic manifolds of the same dimension, the \emph{critical set} $\Crit f$ is the closed set
of points at which $f$ is not \'etale, i.e.\ at which the differential
$df_x:T_xX\to T_{f(x)}X'$ is not bijective, and
$\Crv_f=f(\Crit f)$ is the set of \emph{critical values}.

\subsection{Parabolic subgroups}\label{sec:parabolics}

Parabolic subgroups of $\bG$ are always assumed defined over $F$, and for such a $\bP$ we write $\uu_P$ for the Lie algebra of its unipotent radical.

For a cocharacter $\eta:\Gm\to\bG$ defined over $F$ we write $P_\bG(\eta)$ for the parabolic subgroup
\[
P_\bG(\eta)=\{g:\lim_{t\to0}\eta(t)g\eta(t)^{-1}\text{ exists}\}
\]
with unipotent radical
\[
U_\bG(\eta)=\{g:\lim_{t\to0}\eta(t)g\eta(t)^{-1}=1\}
\]
and Levi decomposition $P_\bG(\eta)=Z_\bG(\eta)\ltimes U_\bG(\eta)$.
Every parabolic subgroup is of this form.
The multiplication map $U_\bG(-\eta)\times P_\bG(\eta)\to\bG$ is an open immersion whose image is the big cell.

We fix a maximal $F$-split torus $\bA_0$ of $\bG$ and a minimal parabolic subgroup $\bP_0\supset\bA_0$.
Put $A_0=\bA_0(F)$.
A standard (resp., semistandard) parabolic subgroup is one containing $\bP_0$ (resp., $\bA_0$).
There are finitely many semistandard parabolic subgroups.
Denote by $\Pstd$ (resp., $\Pssall$) the set of all proper standard (resp.,
semistandard) parabolic subgroups, and by $\Pc\subset\Pstd$ (resp.,
$\Pcss\subset\Pssall$) the subset of the maximal ones.
Every proper parabolic subgroup is $G$-conjugate to a unique element of $\Pstd$.
Every semistandard parabolic subgroup $\bP$ has a unique Levi subgroup $\bM_P$ containing $\bA_0$:
we can write $\bP=P_\bG(\eta)$ for a cocharacter $\eta\in X_*(\bA_0)$ which is positive on its unipotent roots and zero on its Levi roots.
Then $\bM_P=Z_\bG(\eta)$ and we write $U_P^-=U_\bG(-\eta)$ for the opposite unipotent radical; it depends only on $\bP$.
Any parabolic subgroup in $\Pssall$ (resp., $\Pcss$) is of the form $\dot w\bP\dot w^{-1}$ with $\bP\in\Pstd$ (resp., $\bP\in\Pc$) and $\dot w$ in the normalizer $N_G(A_0)$ of $\bA_0$ in $G$.

By a Levi subgroup we mean a Levi subgroup, defined over $F$, of a
parabolic subgroup of $\bG$; it is semistandard if it contains $\bA_0$,
i.e., if it is $\bM_P$ for a semistandard $\bP$. Every Levi subgroup is
$G$-conjugate to a semistandard one. Let $\Lcl$ be the finite set of
$G$-conjugacy classes of Levi subgroups, and write $[M]$ for the class of
$\bM$. These are classes of subgroups of $\bG$ up to conjugation by $G$, not
isomorphism classes. For $[L],[M]\in\Lcl$ write $[L]\preceq[M]$ if a
$G$-conjugate of $\bL$ is contained in $\bM$ (it is then a Levi subgroup of
$\bM$), and $[L]\prec[M]$ if moreover $[L]\ne[M]$. Comparing dimensions shows
that $\preceq$ is a partial order on $\Lcl$; its largest element is $[G]$.

We write $\Nss=\abs{\Pcss}$ for the number of maximal semistandard parabolic
subgroups; thus $\Nss=0$ if and only if $\Pc=\emptyset$, i.e., if and only if
$\bG^{\der}$ is anisotropic over $F$. We have
\[
\Nss=\sum_{\bP\in\Pc}[W_F:W_{F,M_P}],
\]
where $W_F=N_G(A_0)/Z_G(A_0)$ (resp., $W_{F,M_P}$) is the relative Weyl group of
$G$ (resp., of $M_P$). For example, for $\bG=\GL_d$, $\Nss=2^d-2$.

\begin{remark}[Counting maximal parabolic subgroups]\label{rem:countpar}
We have $N_{\bG}\le N_{\bG_{\bar F}}$.
Indeed, let $\bT$ be a maximal $F$-torus containing $\bA_0$. Every semistandard
parabolic subgroup of $\bG$ contains $\bT$.
The absolute Galois group $\Gal(\bar F/F)$ acts on the
Coxeter complex $\mathscr{C}_W$ of the Weyl group of $(\bG_{\bar F},\bT_{\bar F})$.
The maximal proper parabolic subgroups of $\bG_{\bar F}$ containing $\bT$ correspond
to the vertices of $\mathscr{C}_W$, and the maximal proper semistandard parabolic
$F$-subgroups to the minimal nonempty Galois-stable faces of $\mathscr{C}_W$.
The latter are mutually disjoint, so their number is at most the number of vertices, which is $N_{\bG_{\bar F}}$.
\end{remark}

\subsection{Representations}\label{sec:representations}
By a representation of $G$ we tacitly mean a smooth representation on a complex vector space.

For a representation $(\pi,V)$ of $G$ and $\varphi\in\Hk(G)$ let $\pi(\varphi)=\int_G\pi(g)\,d\varphi(g)$.
For a parabolic subgroup $\bP=\bM\bU$ let $V_U=V/\langle\pi(u)v-v\rangle$ be the (unnormalized) Jacquet module.
For a compact open subgroup $J$ and a character $\chi$ of $J$, the operator $\pi(e_J^\chi)$ is the projection onto the
$\chi$-isotypic subspace $V^{J,\chi}$.

Let $\Irr(G)$ be the set of isomorphism classes of irreducible representations of $G$, let
$\Gsc\subset\Gds\subset\Gtemp\subset\Irr(G)$
be the subclasses of unitarizable supercuspidal, discrete series and
tempered representations. By discrete series we mean unitarizable
irreducible representations that are square-integrable modulo $Z$.
Recall that every irreducible representation is admissible.

For a parabolic subgroup $\bP=\bM\bU$ let
$\delta_P(m)=\abs{\det(\Ad(m)|_{\uu_P})}$, $m\in M$, be the modulus
character of $P$, let $i_P^G$ be the normalized parabolic induction, and for a
representation $(\pi,V)$ of $G$ let
\[
r_P(V)=\delta_P^{-1/2}\otimes V_U
\]
be the normalized Jacquet module, a representation of $M$; we also write
$r_P(\pi)$. Thus $r_P(\pi)\ne0$ if and only if $V_U\ne0$, and $r_G(\pi)=\pi$. The functor $r_P$ is exact and left
adjoint to $i_P^G$, and both functors preserve finite length.
Let $A_M$ be the group of $F$-points of the maximal split torus in the center of
$\bM$. If $V$ is of finite length, then $r_P(V)$ is the direct sum of its
generalized eigenspaces with respect to $A_M$; the characters of $A_M$ which occur are the
\emph{exponents} of $V$ along $P$. If moreover $V$ is tempered, its
\emph{\weakjac{}} $r_P^{\wlab}(V)$ is the sum of the generalized eigenspaces of
the unitary exponents. It is a direct summand
of $r_P(V)$, and $r_G^{\wlab}(V)=V$. The conditions $r_P(\pi)\ne0$ and $r_P^{\wlab}(\pi)\ne0$ depend only on the
$G$-conjugacy class of $\bP$.

We recall Casselman's criteria for temperedness and square-integrability. For a parabolic subgroup
$\bP=\bM\bU$ let $\Delta(P)$ be the set of simple roots of $A_M$ in $\uu_P$.
Then $\pi\in\Irr(G)$ is tempered if and only if for every parabolic subgroup
$\bP=\bM\bU$ of $\bG$, including $\bG$, and every exponent $\chi$ of $\pi$ along
$P$, the character $\abs\chi$ of $A_M$ is trivial on $A_G$ and
\begin{equation}\label{eq:casselman}
\abs{\chi(a)}=\prod_{\alpha\in\Delta(P)}\abs{\alpha(a)}^{c_\alpha}\quad(a\in A_M)
\qquad\text{with real numbers $c_\alpha\ge0$.}
\end{equation}
Moreover, $\pi$ is a discrete series representation if and only if in
addition $c_\alpha>0$ for all $\bP$, $\chi$ and $\alpha\in\Delta(P)$. In particular,
$r_P^{\wlab}(\pi)=0$ for $\pi\in\Gds$ and $\bP$ proper.

For every $\pi\in\Irr(G)$ there are a parabolic subgroup $\bQ=\bL\bU_Q$ and
an irreducible supercuspidal representation $\sigma$ of $L$ such that $\pi$
embeds in $i_Q^G\sigma$, and the pair $(\bL,\sigma)$ is unique up to $G$-conjugacy. We write
$\Lsup{\sclab}(\pi)=[L]\in\Lcl$. Similarly, for every $\pi\in\Gtemp$ there are
$\bQ=\bL\bU_Q$ and a discrete series representation $\delta$ of $L$ such
that $\pi$ is a direct summand of the unitary representation $i_Q^G\delta$, and $(\bL,\delta)$ is unique up to
$G$-conjugacy. We write $\Lsup{\dslab}(\pi)=[L]$.
The supercuspidal support of $\pi$ is that of $\delta$, so that
\begin{equation}\label{eq:labels}
\Lsup{\sclab}(\pi)\preceq\Lsup{\dslab}(\pi),\quad\text{with equality if and only if $\delta$ is supercuspidal.}
\end{equation}
The support classes can be recovered as the least classes in the
corresponding Jacquet profiles:
\begin{equation}\label{eq:least}
\begin{aligned}
\Lsup{\sclab}(\pi)&=\min_{\preceq}\{[M_P]:r_P(\pi)\ne0\},\\
\Lsup{\dslab}(\pi)&=\min_{\preceq}\{[M_P]:r_P^{\wlab}(\pi)\ne0\}
\qquad(\pi\in\Gtemp).
\end{aligned}
\end{equation}
The weak Jacquet modules along the maximal proper
parabolic subgroups already detect whether a tempered representation is discrete series.

For $[M]\in\Lcl$ and $\ast\in\{\sclab,\dslab\}$ put
\begin{equation}\label{eq:spectralsets}
\GhatL{\ast}{[M]}=\{\pi\in\Gtemp:\Lsup{\ast}(\pi)=[M]\}.
\end{equation}
For each $\ast$ these sets form a partition of $\Gtemp$ into finitely many
measurable subsets (unions of components of $\Gtemp$, where the components are indexed by the $G$-conjugacy classes
of pairs $(\bL,\mathcal O)$ consisting of a Levi subgroup $\bL$ and an orbit $\mathcal O$ of discrete series representations
of $L$ under twisting by unitary unramified characters), and $\GhatL{\ast}{[G]}=\Ghat{\ast}$.

We use Harish-Chandra's Plancherel formula
\begin{equation}\label{eq:plancherel}
(\varphi/dg)(1)=\int_{\Gtemp}\tr\pi(\varphi)\,d\mu_G(\pi),\quad\varphi\in\Hk(G),
\end{equation}
where $\mu_G$ is the Plancherel measure.
Rescaling the Haar measure $dg$ by $c$ rescales $\mu_G$ by $c^{-1}$.
If $e\in\Hk(G)$ is an idempotent, then $\tr\pi(e)=\dim\pi(e)V$, and
\[
0\le \int_A\tr\pi(e)\,d\mu_G(\pi)\le(e/dg)(1)
\]
for any measurable subset $A$ of $\Gtemp$.
In particular, for a compact open subgroup $K$ of $G$ and $\ast\in\{\sclab,\dslab\}$ let
\[
\dens{\ast}{K}=\vol(K)\int_{\Ghat{\ast}}\dim\pi^K\,d\mu_G(\pi),
\]
which does not depend on $dg$. We call $\scd{K}$ and $\dsd{K}$ the
\emph{supercuspidal density} and the \emph{discrete series density} of $K$.
They are the proportions of the Plancherel mass of $\vol(K)e_K$ carried by
supercuspidal and discrete series representations, respectively.
Thus,
\[
0\le\scd{K}\le\dsd{K}\le1.
\]
More generally, for $[M]\in\Lcl$ let
\begin{equation}\label{eq:densityM}
\densL{\ast}{[M]}{K}=\vol(K)\int_{\GhatL{\ast}{[M]}}\dim\pi^K\,d\mu_G(\pi),
\end{equation}
so that $\dens{\ast}{K}=\densL{\ast}{[G]}{K}$ and, by \eqref{eq:plancherel} for
$\varphi=e_K$, $\sum_{[M]\in\Lcl}\densL{\ast}{[M]}{K}=1$ for each $\ast$.

\subsection{Elliptic linear forms}\label{sec:ellipticforms}

Let
\[
\Par{\gd}=\bigcup_{\bP}\uu_P^\perp,
\]
the union over all proper parabolic subgroups $\bP$ of $\bG$. Thus,
\begin{equation}\label{eq:N}
\Par{\gd}=\bigcup_{\bP\in\Pc}\Ad^*(G)\uu_P^\perp.
\end{equation}
The complement
\begin{equation}\label{eq:E}
\Ell{\gd}=\gd\setminus\Par{\gd}
\end{equation}
is the \emph{elliptic} part of $\gd$: it consists of the linear forms whose restriction to $\uu_P$ is nonzero
for every proper $\bP$. This is the notion of elliptic element of $\gd$ of
Beuzart-Plessis \cite[Appendix~D]{MR4325284}. For a lattice $\dLat\subset\gd$, the
\emph{elliptic density} of $\dLat$ is
\[
\elld{\dLat}=\elld[\bG]{\dLat}=\nu_{\dLat}(\dLat\cap\Ell{\gd}),
\]
where $\nu_{\dLat}$ is the probability Haar measure on $\dLat$.
Null sets in $\gd$ are taken with respect to a Haar measure.

In characteristic zero, a non-degenerate invariant bilinear form on $\g$ identifies $\gd$ with $\g$.
Under this identification, $\uu_P^\perp$ is $\Lie(\bP)(F)$, and up to a null set, $\Ell{\gd}$
is the set of regular semisimple elements whose centralizer is an elliptic maximal torus.
Indeed, the regular semisimple elements form a conull set, a regular
semisimple element of $\Lie(\bP)$ is conjugate under $\bU_P(F)$ to an
element of $\m$, and a maximal torus lies in a proper Levi subgroup if and
only if it is not elliptic.
However, in characteristic two, $\Par{\gd}$ is a null set for $\PGL_2$ even though $\PGL_2$ is split (\zcref{ex:pgl2}).

\subsection{The incidence maps}\label{sec:incidencemaps}
For $\bP\in\Pstd$ let $P=\bP(F)$ and let $\bFl_P=\bG/\bP$ be the (projective)
variety of parabolic subgroups conjugate to $\bP$.
Recall that the map $G\to\bFl_P(F)$ is surjective and that
parabolic $F$-subgroups which are conjugate over $\bar F$ are conjugate under $G$.
Thus $\Fl_P=\bFl_P(F)=G/P$ is the (compact, $F$-analytic) manifold of parabolic
$F$-subgroups of $\bG$ conjugate to $\bP$.
Let $\bInc_P$ be the \emph{incidence variety}, the closed subvariety of
$\bFl_P\times\Lie(\bG)^*$ consisting of the pairs $(\bQ,\lf)$ such that $\lf$
vanishes on the Lie algebra of the unipotent radical of $\bQ$, and let
$\incm_P:\bInc_P\to\Lie(\bG)^*$ be the projection. Via the first projection, the
incidence variety $\bInc_P$ is the vector bundle
$\bG\times^{\bP}\Lie(\bU_P)^\perp$ over $\bFl_P$, where $\bU_P$ is the unipotent
radical of $\bP$; in particular it is smooth. Its $F$-points are
\begin{equation}\label{eq:Y}
\Inc_P=\bInc_P(F)=\{(\bQ,\lf)\in\Fl_P\times\gd:\lf|_{\uu_Q}=0\},\qquad \incm_P(\bQ,\lf)=\lf.
\end{equation}
This is the associated bundle $G\times^P\uu_P^\perp$, i.e., the
coadjoint version of the partial Grothendieck--Springer resolution
(cf.~\cite{MR737927} for the adjoint version $G\times^P\Lie(\bP)(F)\to\g$ over $\C$,
and \cite[\S2.4]{Chi26} for the coadjoint version for Borel subgroups);
its subbundle $G\times^P\Lie(\bP)(F)^\perp$, with fibre
$\Lie(\bP)(F)^\perp\subset\uu_P^\perp$, is the cotangent bundle of $\Fl_P$
\cite[Lemma~1.4.9]{MR1433132}.
In particular, $G$ acts on $\Inc_P$ on the left by
$g\cdot(\bQ,\lf)=(g\bQ g^{-1},\Ad^*(g)\lf)$, and $\incm_P$ is $G$-equivariant.

$\Inc_P$ is a vector subbundle of $\Fl_P\times\gd$ of rank $D-\dim\uu_P$
over the manifold $\Fl_P$ of dimension $\dim\uu_P$.
Hence, $\Inc_P$ is an analytic manifold of dimension $D$ and $\incm_P$ is proper.
Thus, the image
\begin{equation}\label{eq:Pim}
\Pim{P}=\incm_P(\Inc_P)=\Ad^*(G)\uu_P^\perp
\end{equation}
is closed. It is the set of linear forms which vanish on $\uu_Q$ for some
conjugate $\bQ$ of $\bP$; it is stable under $\Ad^*(G)$ and under
multiplication by $F^\times$. We also put $\Pim{G}=\gd$, and use the notation
$\Pim{P}$ for an arbitrary parabolic subgroup $\bP$; it depends only on the
conjugacy class of $\bP$.

Taking the disjoint union over the maximal parabolic subgroups in $\Pc$ we obtain an analytic manifold
\[
\Inc=\bInc(F),\qquad \bInc=\bigsqcup_{\bP\in\Pc}\bInc_P,
\]
of dimension $D$ and a proper map $\incm:\Inc\to\gd$ with image $\Par{\gd}$ by \eqref{eq:N}.
In particular, $\Ell{\gd}$ is open. Let
\begin{equation}\label{eq:Sigma}
\Crvmax=\Crv_{\incm}=\bigcup_{\bP\in\Pc}\incm_P(\Crit \incm_P)
\qquad\text{and}\qquad
\Crv=\bigcup_{\bP\in\Pstd}\incm_P(\Crit \incm_P)\supset\Crvmax
\end{equation}
be the set of critical values of $\incm$ and the set of critical values of the
incidence maps of all proper parabolic subgroups, respectively. Since the
maps $\incm_P$ are proper, $\Crvmax$ and $\Crv$ are closed, and
they are null sets by Sard's theorem (see \zcref{thm:sard}). They are
stable under $\Ad^*(G)$ and under multiplication by $F^\times$.

If $\lf\in\Pim{P}\setminus\incm_P(\Crit\incm_P)$, then $\lf=\incm_P(y)$ with $d\incm_P$ invertible at $y$,
so $\lf$ is an interior point of $\Pim{P}$ by the inverse function theorem. Hence
\begin{equation}\label{eq:boundary}
\partial\Pim{P}\subset\incm_P(\Crit\incm_P)\subset\Crv\quad(\bP\in\Pstd),\qquad
\partial\Ell{\gd}\subset\Crvmax;
\end{equation}
in particular, $\partial\Ell{\gd}$ is a null set.

For $[M]\in\Lcl$ let
\begin{equation}\label{eq:PimM}
\Pim{[M]}=\bigcup_{\bP:[M_P]=[M]}\Pim{P},
\end{equation}
the union over the parabolic subgroups $\bP$, including $\bG$, which have a
Levi subgroup in the class $[M]$; it is a finite union, over the corresponding elements
of $\Pstd\cup\{\bG\}$. Thus $\Pim{[G]}=\gd$ and
$\Par{\gd}=\bigcup_{[M]\ne[G]}\Pim{[M]}$.

\begin{lemma}[The least Levi class of a linear form]\label{lem:meet}
Let $\bP$ and $\bP'$ be parabolic subgroups and let $\lf\in\gd$ vanish on
$\uu_P$ and on $\uu_{P'}$. Then there is a parabolic subgroup $\bQ\subset\bP$
such that $\lf$ vanishes on $\uu_Q$, $[M_Q]\preceq[M_P]$ and $[M_Q]\preceq[M_{P'}]$.
Consequently, for every $\lf\in\gd$ the set $\{[M]\in\Lcl:\lf\in\Pim{[M]}\}$ has
a least element, which we denote by $\Lform{\lf}\in\Lcl$.
\end{lemma}

\begin{proof}
Recall that $\bP\cap\bP'$ contains the centralizer of a maximal split torus
$\bA$ of $\bG$, and that $\bQ=(\bP\cap\bP')\bU_P$ is a parabolic subgroup of
$\bG$, contained in $\bP$. Let $\bM_P$ and $\bM_{P'}$ be the Levi subgroups of
$\bP$ and $\bP'$ containing $\bA$. Then $\bQ=(\bM_P\cap\bP')\ltimes\bU_P$, and
$\bM_P\cap\bP'$ is a parabolic subgroup of $\bM_P$ with Levi subgroup
$\bM_P\cap\bM_{P'}$ and unipotent radical $\bM_P\cap\bU_{P'}$. Hence
$\bM_Q=\bM_P\cap\bM_{P'}$ and
\[
\uu_Q=\uu_P\oplus(\m_P\cap\uu_{P'})\subset\uu_P+\uu_{P'},
\]
so that $\lf$ vanishes on $\uu_Q$. The last part follows from finiteness.
\end{proof}

The notation is parallel to $\Lsup{\sclab}(\pi)$ and $\Lsup{\dslab}(\pi)$ of
\eqref{eq:least}. The two are related in \zcref{prop:partition}.

The \emph{\levistratum{}} of $[M]$ is
\begin{equation}\label{eq:strata}
\Estr_{[M]}=\{\lf\in\gd:\Lform{\lf}=[M]\}
=\Pim{[M]}\setminus\bigcup_{[L]\prec[M]}\Pim{[L]}.
\end{equation}
Thus the \levistrata{} form a partition of $\gd$ into subsets which are stable
under $\Ad^*(G)$ and under multiplication by scalars in $F^\times$, each of them is the difference of two closed
sets, and $\Estr_{[G]}=\Ell{\gd}$. By \eqref{eq:boundary}, the boundary of each
\levistratum{} is contained in the null set $\Crv$, and the sets
$\Estr_{[M]}\setminus\Crv$ are open and closed in $\gd\setminus\Crv$. Moreover, for
$\lf\notin\Crv$ the classes $[M]$ with $\lf\in\Pim{[M]}$ are exactly those with
$\Lform{\lf}\preceq[M]$ (\zcref{prop:upset}). For a lattice $\dLat\subset\gd$ put
\begin{equation}\label{eq:levid}
\levid{[M]}{\dLat}=\nu_{\dLat}(\dLat\cap\Estr_{[M]}),
\end{equation}
so that $\levid{[G]}{\dLat}=\elld{\dLat}$ and $\sum_{[M]\in\Lcl}\levid{[M]}{\dLat}=1$.

\begin{remark}[Levi strata in characteristic zero]\label{rem:toral}
Suppose that $\chr F=0$ and identify $\gd$ with $\g$ by a non-degenerate
invariant bilinear form, as in \zcref{sec:ellipticforms}. For a regular
semisimple $\lf$ let $\bT_\lf$ be its centralizer, a maximal torus, and let
$\bA_\lf$ be the maximal split subtorus of $\bT_\lf$. Then $\Lform{\lf}$ is the
class of the Levi subgroup $Z_{\bG}(\bA_\lf)$. Since the regular semisimple
elements form a conull set, this describes the \levistrata{} up to null sets.
Indeed, $\lf\in\Pim{P}$ if and only if a conjugate of $\lf$ lies in
$\Lie(\bP)(F)$, and a regular semisimple element of $\Lie(\bP)(F)$ is conjugate
under $U_P$ to an element of $\m_P$. Hence $\lf\in\Pim{[M']}$ if and only if a
conjugate of $\bT_\lf$ is contained in $\bM'$. A Levi subgroup containing
$\bT_\lf$ is the centralizer of a split torus contained in $\bA_\lf$, hence contains
$Z_{\bG}(\bA_\lf)$; conversely $\bT_\lf\subset Z_{\bG}(\bA_\lf)$. This description will not be used in the sequel.
\end{remark}

\begin{definition}[Peculiar groups]\label{def:peculiar}
We say that $\bG$ is \emph{\peculiar{}} if for some $\bP\in\Pc$ the incidence
map $\incm_P$ is everywhere critical, i.e., $\Crit \incm_P=\Inc_P$.
\end{definition}

\Peculiar{} groups exist only when $\chr F=2$.
We give a precise characterization in \zcref{prop:peculiar}.

Let $\dG{}$ be the maximal dimension of the unipotent radical of a
maximal proper parabolic subgroup of $\bG$, with the convention
$\dG{}=1$ if there is none, and put
\begin{equation}\label{eq:tau}
\eexp=\eexp_{\bG}=
\begin{cases}
1+1/\dG{}&\text{if $\bG$ is not \peculiar{},}\\
1&\text{if $\bG$ is \peculiar{}.}
\end{cases}
\end{equation}
This will appear in the error exponent for the supercuspidal and discrete series
densities. For the other Levi classes all proper parabolic subgroups
intervene: let $\dGall{}=\dim\uu_{P_0}$ be the maximal dimension of the
unipotent radical of a proper parabolic subgroup of $\bG$, with the convention
$\dGall{}=1$ if there is none, and put
\begin{equation}\label{eq:tauall}
\eexpall=
\begin{cases}
1+1/\dGall{}&\text{if $\bG$ is not \peculiar{},}\\
1&\text{if $\bG$ is \peculiar{}.}
\end{cases}
\end{equation}
Thus $1\le\eexpall\le\eexp\le2$.

\begin{example}[$\GL_d$]\label{ex:gld}
Let $\bG=\GL_d$ and identify $\gd$ with $\Mat_d(F)$ by the trace
pairing. For the maximal parabolic subgroup stabilizing a nonzero
proper subspace $W$, the Lie algebra $\uu_P=\Hom(F^d/W,W)$ has
annihilator equal to the set of endomorphisms preserving $W$.
Hence $\Ell{\gd}$ consists of the matrices with no nonzero proper
invariant subspace, or equivalently with characteristic polynomial
irreducible over $F$ (in any characteristic, including inseparable
polynomials). For $\Lat=\Mat_d(\OO)$ we have $\dLat=\Mat_d(\OO)$, and
reduction modulo $\pp$ gives the lower bound
\[
\elld[\GL_d]{\Mat_d(\OO)}\ge\frac{I_d(q)\abs{\GL_d(k)}}{q^{d^2}(q^d-1)},
\]
where $I_d(q)=\frac1d\sum_{e\mid d}\mu(e)q^{d/e}$ is the number of
monic irreducible polynomials of degree $d$ over $k$: each of them is
the characteristic polynomial of $\abs{\GL_d(k)}/(q^d-1)$ matrices over
$k$, and every lift of such a matrix has irreducible characteristic
polynomial over $F$.
\end{example}

\begin{example}[Levi classes for $\GL_d$]\label{ex:gldlevi}
Let $\bG=\GL_d$ and identify $\gd$ with $\Mat_d(F)$ by the trace pairing, as
in \zcref{ex:gld}. The set $\Lcl$ is in bijection with the partitions of $d$:
the partition $d=d_1+\dots+d_k$ corresponds to the class of
$\GL_{d_1}\times\dots\times\GL_{d_k}$. For the parabolic subgroup $\bP$
stabilizing a flag of subspaces of $F^d$, the annihilator $\uu_P^\perp$ is the set of
endomorphisms preserving this flag. Hence, for $\lf\in\Mat_d(F)$ and $[M]$
the class of $\GL_{d_1}\times\dots\times\GL_{d_k}$, we have $\lf\in\Pim{[M]}$
if and only if the $F[\lf]$-module $F^d$ admits a filtration by submodules whose
successive quotients have dimensions $d_1,\dots,d_k$, in some order, or equivalently,
the characteristic polynomial $\chi_\lf$
of $\lf$ admits a factorization into polynomials over $F$ of degrees $d_1,\dots,d_k$.

Consequently, $\Lform{\lf}$ is the class of the partition formed by the degrees
of the irreducible factors of $\chi_\lf$, counted with multiplicity (these are
the dimensions of the composition factors of $F^d$), and for every $\lf$ the
classes $[M]$ with $\lf\in\Pim{[M]}$ are exactly those with
$\Lform{\lf}\preceq[M]$ (compare \zcref{prop:upset}). Thus the \levistratum{}
$\Estr_{[M]}$ is the set of matrices whose characteristic polynomial is a
product of irreducible polynomials, not necessarily distinct, of degrees
$d_1,\dots,d_k$, and for every lattice
$\dLat\subset\Mat_d(F)$ the density $\levid{[M]}{\dLat}$ is the probability
that the characteristic polynomial of a random element of $\dLat$ is such a product.
\end{example}

\subsection{Standard groups}\label{sec:standardgroups}
By a \emph{\nchart{}} of $G$ we mean an $F$-analytic isomorphism $\kappa$ from a
neighborhood of $1$ in $G$ onto a neighborhood of $0$ in $\g$ with
$\kappa(1)=0$ and $d\kappa_1=\id$; \nchart{}s exist since $G$ is an
$F$-analytic manifold ($\bG$ being smooth over $F$).
We can take $\kappa$ to be the restriction, to a suitable neighborhood of $1$
in $G$, of an algebraic map $\bG\to\g$ which is \'etale at $1$; both the model
charts and the \nchart{} attached to a linear embedding in \zcref{sec:models}
are of this form. Such a map is in general not an isomorphism of any Zariski
open neighborhood of $1$ onto its image, but only an analytic isomorphism near
$1$; the inverse and the group law in these coordinates are given by power
series.
We shall only use that $\kappa$ is analytic and normalized.
We call the inverse of a \nchart{} a \emph{\npar{}}.
We use similar terminology for smooth $F$-subgroups of $\bG$.
If $\bU$ is the unipotent radical of a parabolic $F$-subgroup of $\bG$, then $\bU$ is $F$-isomorphic, as a variety, to an affine space.
Composing such an isomorphism with a translation and a linear change of variables we
obtain a \npar{} $\xi:\Lie(\bU)(F)\to\bU(F)$, defined on all of $\Lie(\bU)(F)$, which is
an isomorphism of $F$-varieties.

Let $(\kappa,\Lat)$ be a pair consisting of a \nchart{} $\kappa$ and a lattice $\Lat$ of $\g$.
We say that $(\kappa,\Lat)$ is \emph{standard} if the image of $\kappa$ contains $\pp\Lat$,
$K_1=\kappa^{-1}(\pp\Lat)$ is a (compact open) subgroup of $G$ and the group law in $K_1$ is given by a formal group law with coefficients in $\OO$:
\[
\kappa(xy)=\Phi(\kappa(x),\kappa(y)),\ x,y\in K_1
\]
where, in the coordinates $X=(X_1,\dots,X_D)$ of a basis of $\Lat$,
the components $\Phi_1,\dots,\Phi_D\in\OO[[X,Y]]$ are formal power series
in the $2D$ variables $X,Y$. This also implies $[\Lat,\Lat]\subset\Lat$:
if the quadratic term of $\Phi(X,Y)$ is $\beta(X,Y)$, then
$[X,Y]=\beta(X,Y)-\beta(Y,X)$, and $\beta$ has integral coefficients.

For every \nchart{} $\kappa$ and lattice $\Lat$ there is $n_0\ge0$ such that $(\kappa,\pp^{n_0}\Lat)$ is standard.

Suppose that $(\kappa,\Lat)$ is standard. Then the sets
\begin{equation}\label{eq:tower}
K_n=\kappa^{-1}(\pp^n\Lat),\qquad n\ge1,
\end{equation}
are compact open subgroups of $G$ which are normal in $K_1$, and $\Ad(K_1)\Lat=\Lat$.
The groups $K_n$ are called the principal congruence subgroups with respect to $\Lat$.
For any $a,b\ge 1$, $x\in K_a$, $y\in K_b$, we have
\begin{subequations}
\begin{gather}
\kappa(xy)-\kappa(x)-\kappa(y)\in\pp^{a+b}\Lat,\qquad
\kappa(x^{-1})+\kappa(x)\in\pp^{2a}\Lat,\label{eq:law}\\
[K_a,K_b]\subset K_{a+b}.\label{eq:commutator}
\end{gather}
\end{subequations}
For $0\le s\le a$, $\kappa$ induces an isomorphism of Abelian groups
\[
K_a/K_{a+s}\rightarrow\pp^a\Lat/\pp^{a+s}\Lat.
\]
Moreover, the pushforward under $\kappa$ of a Haar measure on $K_1$ is a Haar measure on $\pp\Lat$
so that
\begin{equation}\label{eq:volume}
\vol(K_{a+s})=q^{-sD}\vol(K_a)\qquad(a\ge1,\ s\ge0).
\end{equation}

For a fixed lattice $\Lat$, the principal congruence subgroups defined
by two standard pairs $(\kappa,\Lat)$ and $(\kappa',\Lat)$ coincide for all sufficiently large $n$.
More generally, if $(\kappa',\Lat')$ is any standard pair, with principal congruence
subgroups $K'_n$, then there exists $c\ge0$ such that
$K'_{n+c}\subset K_n\subset K'_{n-c}$ for all sufficiently large $n$.
These inclusions compare levels; the limiting densities in our results depend on the lattice.

For the rest of the paper we fix a standard pair $(\kappa,\Lat)$.
This determines the principal congruence subgroups $K_n$, the dual lattice
$\dLat=\Lat^\vee=\{\lf\in\gd:\lf(\Lat)\subset\OO\}$ and the cell characters below.

Auxiliary lattices introduced later do not change these data.

As a general rule, constants denoted $C$, $c$, \dots\ may depend on $F$, $\bG$, $(\kappa,\Lat)$ and
auxiliary choices but not on running variables and parameters $n,r,h$.

\subsection{Models}\label{sec:models}

Let $\cG$ be a smooth affine group scheme of finite type over
$\OO$ with generic fiber $\bG$, let $\OO[\cG]$ be its coordinate ring, and let
$I=\ker(\varepsilon\colon\OO[\cG]\to\OO)$ be the augmentation ideal.
Since $\cG$ is smooth, $I/I^2$ is a free $\OO$-module of rank $D$, canonically
dual to $\Lat=\Lie(\cG)(\OO)$; thus $I/I^2=\Lat^\vee$. Choose an $\OO$-linear
section $s\colon\Lat^\vee\to I$ of the projection $I\to I/I^2$, and define
$\kappa_{\cG}\colon\ker(\cG(\OO)\to\cG(k))\to\pp\Lat$ by $\lf(\kappa_{\cG}(g))=s(\lf)(g)$ for
$\lf\in\Lat^\vee$; note that $f(g)\in\pp$ for $f\in I$ and $g\in\ker(\cG(\OO)\to\cG(k))$.
We call $\kappa_{\cG}$ a \emph{model chart} of $\cG$.
It depends on $s$, but nothing below will depend on it: two choices of $s$ differ by a map into $I^2$, so the corresponding
model charts agree modulo $\pp^{2a}\Lat$ on $\ker(\cG(\OO)\to\cG(\OO/\pp^a))$.

Since $\cG$ is smooth, $s$ extends to an isomorphism of the completed symmetric algebra
of $\Lat^\vee$ onto the $I$-adic completion of $\OO[\cG]$,
so the group law becomes a formal group law over $\OO$, and $\kappa_{\cG}$ is a bijection of
$\ker(\cG(\OO)\to\cG(k))$ onto $\pp\Lat$ which
identifies $\ker(\cG(\OO)\to\cG(\OO/\pp^n))$ with $\kappa_{\cG}^{-1}(\pp^n\Lat)$
for $n\ge1$. Hence $(\kappa_{\cG},\Lat)$ is standard. Hensel's lemma makes $\cG(\OO)\to\cG(k)$ surjective.
For future reference we record these facts in the following lemma, where we also use \eqref{eq:volume}.

\begin{lemma}[Smooth models]\label{lem:model}
Let $\kappa_{\cG}$ be a model chart of $\cG$. Then $\kappa_{\cG}$ is a \nchart{},
$(\kappa_{\cG},\Lat)$ is standard, and $K_n=\ker(\cG(\OO)\to\cG(\OO/\pp^n))$ for all
$n\ge1$. Moreover $[\cG(\OO):K_1]=\abs{\cG(k)}$ and $\vol(K_n)=\abs{\cG(k)}^{-1}q^{-(n-1)D}\vol(\cG(\OO))$.
\end{lemma}

\begin{example}[Parahoric towers]\label{ex:parahoric}
Let $x$ be a point of the Bruhat--Tits building of $\bG$ and let $\cG_x$ be the parahoric group scheme of
Bruhat--Tits (with connected special fibre). Its congruence kernels
$\ker(\cG_x(\OO)\to\cG_x(\OO/\pp^n))$, $n\ge1$,
form the tower of the smooth model $\cG_x$ (\zcref{lem:model}) with
lattice $\Lie(\cG_x)(\OO)$. If $\bG$ splits over a tamely ramified extension
of $F$, they are the Moy--Prasad groups
$G_{x,n}$ of integral depth $n$ (not $G_{x,n+}$)
\cite[Proposition~8.5.16]{MR4520154}. Thus our results apply to the Moy--Prasad
filtration at integral depths, with a limit depending on this lattice.
\end{example}

Another source of towers is a realization of $\bG$ as an algebraic subgroup
of $\GL_d$ over $F$, that is, a homomorphism $\bG\to\GL_d$ over $F$ with trivial
scheme-theoretic kernel; such a homomorphism is a closed immersion, and we
write $\bG\subset\GL_d$.
This gives rise to a lattice and principal congruence subgroups
\begin{equation}\label{eq:embeddingtower}
\Lat^\natural=\g\cap\Mat_d(\OO),\qquad
K^\natural_n=\begin{cases}
G\cap\bigl(1+\pp^n\Mat_d(\OO)\bigr)&n\ge1,\\
G\cap\GL_d(\OO)&n=0.\end{cases}
\end{equation}
The groups $K^\natural_n$ are the kernels of reduction modulo $\pp^n$ on the
$\OO$-points of the schematic closure $\cG^\natural$ of $\bG$ in
$\GL_{d,\OO}$, since $\cG^\natural(\OO)=G\cap\GL_d(\OO)$, and
$\Lat^\natural=\Lie(\cG^\natural)(\OO)$.

\begin{lemma}[Smooth congruence models of a linear embedding]\label{lem:embedding}
For $n\ge1$, let $\mathcal H_n$ be the affine matrix group with coordinate ring
\[
\OO[x_{ij},\det(1+\varpi^nX)^{-1}],\qquad X=(x_{ij})_{1\le i,j\le d}.
\]
The map $X\mapsto 1+\varpi^nX$ identifies $(\mathcal H_n)_F$ with
$\GL_{d,F}$. Via this identification, let $\cG^\natural_n$ be the
schematic closure of $\bG$ in $\mathcal H_n$. For all sufficiently large $n$,
$\cG^\natural_n$ is a smooth affine group scheme over $\OO$. Under
the generic-fiber identification with $\bG$, its Lie lattice is
$\pp^n\Lat^\natural$ and its congruence kernels are
\[
\ker\bigl(\cG^\natural_n(\OO)\to
          \cG^\natural_n(\OO/\pp^m)\bigr)=K^\natural_{n+m}
\qquad(m\ge1).
\]
\end{lemma}

\begin{proof}
The ambient matrix group has product
$(X,Y)\mapsto X+Y+\varpi^nXY$. The schematic closure is flat over
$\OO$ and inherits the group law, since its generic fiber is
schematically dense in it and in its fiber square over $\OO$.

We first choose equations near the identity in the original matrix
coordinates. Since $\bG$ and $\GL_{d,F}$ are smooth, the ideal of
$\bG$ is generated near $1$ by polynomials $f_1,\dots,f_c$ with
independent differentials at $1$, where $c=d^2-\dim\bG$.
Choose their linear terms $\ell_i=df_i|_1$ to form an $\OO$-basis
of the annihilator of $\Lat^\natural$ in $\Mat_d(\OO)^\vee$.
This is possible by an invertible $F$-linear change of generators,
since $\Lat^\natural$ is saturated. Let $D(u)$ be a principal
open neighborhood of $1$ in $\GL_{d,F}$ on which these equations
generate the ideal, with $u$ a polynomial and $u(1)=1$.

For large $n$, the polynomials
\[
F_{i,n}(X)=\varpi^{-n}f_i(1+\varpi^nX)
\]
have integral coefficients and satisfy $F_{i,n}\equiv\ell_i$ modulo
$\pp$. Also $u(1+\varpi^nX)$ is integral and reduces to $1$.
The common zero scheme of the $F_{i,n}$ on
$D(u(1+\varpi^nX))$ has generic fiber $\bG$ on this open set.
Since the reductions of the $\ell_i$ are independent, some
$c\times c$ Jacobian minor is a unit at every point of the special
fiber. On the open set where this minor is invertible, the equations
give a smooth, hence flat, scheme by the Jacobian criterion. It is
the schematic closure of its generic fiber there. The two open
conditions contain the entire special fiber of the ambient matrix
group, so $\cG^\natural_n$ is smooth over $\OO$.

The linear terms of the $F_{i,n}$ are exactly the $\ell_i$, so the
tangent lattice at $X=0$ is $\Lat^\natural$ in $X$ coordinates.
The differential of $X\mapsto1+\varpi^nX$ sends it to
$\pp^n\Lat^\natural$. The $\OO$-points of the closure are the
integral matrices $X$ with $1+\varpi^nX\in G$, since an integral
point factors through the closure whenever its generic point does.
The kernel of reduction modulo $\pp^m$ consists of these points
with $X\in\pp^m\Mat_d(\OO)$, which gives the asserted kernel.
\end{proof}

Apply \zcref{lem:model} at a fixed smooth level $n$ to get the shifted
tower $(K^\natural_{n+m})_{m\ge1}$ from the lattice $\pp^n\Lat^\natural$.
The fixed shift changes only the implied constants in the asymptotic
estimates. The limiting Levi densities are those for
$(\Lat^\natural)^\vee$, since the strata are invariant under scaling
by $F^\times$.

\subsection{Unramified groups}\label{sec:integraldata}
As usual, we say that $\bG$ is \emph{unramified} if
it is quasi-split and splits over an unramified extension of $F$.
Equivalently, $\bG$ has a \emph{reductive model}, that is, a reductive group
scheme $\cG$ over $\OO$ (smooth and affine, with connected reductive geometric
fibres) with generic fibre $\bG$; the groups $\cG(\OO)$ are then exactly the
hyperspecial parahoric subgroups $\cG_x(\OO)$ of $G$ (\zcref{ex:parahoric}), and the
corresponding $\cG_x$ are exactly the reductive models of $\bG$.
For these facts, for the structure theory of reductive group schemes over $\OO$ recalled below, and for the
facts on groups over global fields at the end of this subsection, we
refer to \cite{MR3362641} and \cite{MR4520154}.

By the \emph{unramified case} we mean the following situation: $\bG$ is
unramified, $\cG$ is a reductive model of $\bG$, $\Lat=\Lie(\cG)(\OO)$, and
$\kappa=\kappa_{\cG}$ is the model chart of \zcref{lem:model}, so that
$(\kappa,\Lat)$ is standard and
\begin{equation}\label{eq:integralKn}
K_n=\kappa^{-1}(\pp^n\Lat)=\ker\bigl(\cG(\OO)\to\cG(\OO/\pp^n)\bigr),\qquad n\ge1.
\end{equation}
The model contains closed $\OO$-tori $\cS\subset\cT$ with $\cS$
split and maximal among split tori, and $\cT$ a maximal torus which splits
over the maximal unramified extension.

A Borel subgroup $B$ of $\bG$ containing $T=\cT_F$ is a minimal
parabolic subgroup, and it extends uniquely to a Borel subgroup scheme
of $\cG$ containing $\cT$. Let $\Delta$ be the set of
simple roots of $\bG_{\bar F}$ with respect to $(T,B)$ and $W$ the
absolute Weyl group; the Galois group of $F$ acts on $\Delta$ through
the Galois group of $k$, and the group $X_*(\cT)$ of cocharacters of the
$\OO$-torus $\cT$, i.e., of $T$ over $F$, is $X_*(T_{\bar F})^{\Gal(\bar F/F)}$.
Every parabolic subgroup $\bP$ of $\bG$ is $G$-conjugate to a standard
one, $\bP\supset B$; we take the representatives $\bP\in\Pc$ standard, i.e., $\bA_0=\cS_F$ and $\bP_0=B$ in \zcref{sec:parabolics}.
The type of a standard maximal parabolic subgroup is a Galois orbit
$I\subset\Delta$, the simple roots not in its Levi subgroup. Let $\ell_I$
be the smallest positive integer for which there is a
Galois-invariant cocharacter $\eta_I\in X_*(T_{\bar F})$, determined up to a central cocharacter, satisfying
\begin{equation}\label{eq:cocharchoice}
\langle\alpha,\eta_I\rangle=
\begin{cases}\ell_I&\alpha\in I,\\0&\alpha\in\Delta\setminus I.\end{cases}
\end{equation}
Fix one such $\eta_I$.
It extends to a cocharacter of $\cT$ and satisfies $\bP=P_{\bG}(\eta_I)$.
The closed subgroup schemes
$\cP=P_{\cG}(\eta_I)$, $\cM=Z_{\cG}(\eta_I)$ and $\cU^\pm=U_{\cG}(\pm\eta_I)$
are smooth with generic fibres $\bP$, $\bM$, $\bU^\pm$, $\cP=\cM\ltimes\cU$,
and the multiplication map
\begin{equation}\label{eq:integralcell}
\cU\times\cM\times\cU^-\to\cG
\end{equation}
is an open immersion. Moreover the action of $\Gm$ on $\cU$ by
conjugation through $\eta_I$ extends to an action of the multiplicative
monoid $\Aff^1$, so that conjugation by $t=\eta_I(\varpi)$ is an endomorphism of the
$\OO$-scheme $\cU$; similarly for $t^{-1}$ on $\cU^-$.

The quotient
$\cG/\cP$ is a smooth projective $\OO$-scheme, $\Fl_P=\bFl_P(F)=(\cG/\cP)(\OO)$,
and the big cell $\cU^-\cP/\cP\simeq\cU^-$ is an open subscheme, and the Iwasawa decomposition $G=P\,\cG(\OO)$ holds.

Let $W_k=N_{\cG}(\cS)(k)/Z_{\cG}(\cS)(k)$ be the relative Weyl group of $\cG_k$ and $W_{k,P}$ that of $\cM_k$.
We may lift the Bruhat decomposition of $\cG_k$ to $\OO$.
More precisely, first, we may lift any $w\in W_k/W_{k,P}$ to $\dot w\in N_{\cG}(\cS)(\OO)$.
Put $\bP_w=\dot w\bP\dot w^{-1}$, and let $\cP_w$ and $\cU^\pm_w$ be the conjugates of $\cP$ and $\cU^\pm$ by $\dot w$. Then
\begin{equation}\label{eq:integralcover}
(\cG/\cP)(\OO)=\bigcup_{w\in W_k/W_{k,P}}\dot w\,\cU^-(\OO)\cP(\OO)/\cP(\OO);
\end{equation}
that is, every parabolic subgroup conjugate to $\bP$ is $u\bP_wu^{-1}$ for
some $w$ and some $u\in\cU^-_w(\OO)$. Moreover, the $\bP_w$,
$w\in W_k/W_{k,P}$, are pairwise distinct elements of $\Pcss$.

Finally, $\cU^\pm$ is isomorphic as an $\OO$-scheme to an affine space.
Write $d_P=\dim_F\uu_P=\dim_F\uu_P^-$.
We fix such isomorphisms $\xi_\pm:\Aff^{d_P}\to\cU^\pm$ such that $\xi_\pm(0)=1$ and $d\xi_\pm$ is the identity with
respect to $\Lat_\pm=\Lie(\cU^\pm)(\OO)$; the lattices $\Lat_\pm$
are direct summands of $\Lat=\Lie(\cG)(\OO)$ and are the direct sums
of their intersections with the weight spaces of $\eta_I$.

In the coordinates $\kappa$, the group law, the inverse, the conjugation by
any element of $\cG(\OO)$, the maps \eqref{eq:integralcell} and its inverse,
and the two contracting conjugations, by $t$ on $\cU$ and by $t^{-1}$ on
$\cU^-$, are given by power series with coefficients in $\OO$: they are
morphisms of $\OO$-schemes between neighborhoods of the identity, computed in
the formal completions along the identity sections.

Define
\begin{equation}\label{eq:mG}
m_{\cG}=\max_{\bP\in\Pc}\ \max_{\langle\alpha,\eta_I\rangle>0}\langle\alpha,\eta_I\rangle,
\end{equation}
the largest weight of $\eta_I$ on $\uu_P$, the inner maximum being over
the roots $\alpha$ of $\bG_{\bar F}$ and $I$ the type of $\bP$. If $\Pc=\emptyset$, then $\bG$, being quasi-split, is a torus, and we set $m_{\cG}=0$.

For instance, for $\bG=\GL_d$, $d\ge2$, one has $m_{\cG}=1$, using $\eta_j(t)=\diag(t\one_j,\one_{d-j})$,
while for $\bG=\SL_d$ one has $m_{\cG}=d$: the cocharacter attached to the Levi subgroup $\mathrm S(\GL_j\times\GL_{d-j})$
is $t\mapsto(t^{(d-j)/g}\one_j,t^{-j/g}\one_{d-j})$ with $g=\gcd(j,d)$, of weight $d/g$ on $\uu_P$.

The invariants $\Nss$ and $m_{\cG}$ depend only on the based absolute root datum of $(\bG,B,T)$
with the action of $\Gal(\bar F/F)$ on it.
They control the constants in the estimates of
\zcref{lem:integralatlas}, \zcref{prop:regularunram}, and \zcref{thm:uniform}.

Finally, we record what is needed in the global setting of \zcref{thm:global}.
Let $\glf$ be a global field, $\bG$ a connected reductive group over
$\glf$ and $\bG\subset\GL_d$ a closed embedding over $\glf$. For a finite place
$v$ of $\glf$, with completion $\glf_v$, ring of integers $\OO_v$, maximal ideal $\pp_v$ and residue
field of cardinality $q_v$, let $\Lat_v$ and $K_{v,n}$ $(n\ge0)$ be the lattice $\Lat^\natural$ and the groups
$K^\natural_n$ of \eqref{eq:embeddingtower} for $F=\glf_v$ and the embedding
$\bG_{\glf_v}\subset\GL_{d,\glf_v}$, and let $\cG_v$ be the schematic closure of
$\bG_{\glf_v}$ in $\GL_{d,\OO_v}$. Let $S$ be the set of finite places $v$ of
$\glf$ for which $\cG_v$ is not a reductive group scheme. It is a finite set.
Indeed, the schematic closure of $\bG$ in $\GL_d$ over the ring of integers of
$\glf$ (or, if $\glf$ is a function field, over the ring of elements of $\glf$
which are integral outside an auxiliary place) is a flat affine group scheme of
finite type whose fibres are connected and reductive over a dense open subset,
and its formation commutes with the flat base change to $\OO_v$.
For $v\notin S$, the group $\bG_{\glf_v}$ is unramified with reductive model
$\cG_v$, we have $\Lat_v=\Lie(\cG_v)(\OO_v)$ and
$K_{v,n}=\ker(\cG_v(\OO_v)\to\cG_v(\OO_v/\pp_v^n))$ for $n\ge1$, as noted after
\eqref{eq:embeddingtower}, and $K_{v,0}=\bG(\glf_v)\cap\GL_d(\OO_v)=\cG_v(\OO_v)$
is a hyperspecial parahoric subgroup of $\bG(\glf_v)$.

The absolute Galois group of $\glf$ acts on the based absolute root datum of $\bG$ through a finite quotient $\Gamma$ (for quasi-split
$\bG$, this is the action considered above), and its restriction to the decomposition group at $v$ is
the corresponding action for $\bG_{\glf_v}$. Hence, if $\bG_{\glf_v}$ is unramified, its invariant
\eqref{eq:mG} depends only on the image of the decomposition group in $\Gamma$.
Let $m_{\bG}$ be the maximum of these invariants over all such $v$. It
is finite since $\Gamma$ has finitely many subgroups, and it depends only on $\bG$. Thus
$m_{\cG_v}\le m_{\bG}$ for $v\notin S$. Moreover, the number of maximal
semistandard parabolic subgroups satisfies $N_{\bG_{\glf_v}}\le N_{\bG_{\bar\glf}}$ for every $v$, by
\zcref{rem:countpar}.

\subsection{Cells}\label{sec:cells}
Let $n\ge1$. For $\lf\in \dLat$, let
\begin{equation}\label{eq:chi}
\chi_{n,\lf}(g)=\psi\bigl(\varpi^{-n}\lf(\kappa(g))\bigr),\qquad g\in K_{\lceil n/2\rceil}.\footnote{Replacing $\varpi$
by $u\varpi$ with $u\in\OO^\times$ replaces $\chi_{n,\lf}$ by $\chi_{n,u^{-n}\lf}$.}
\end{equation}
This is a character of $K_{\lceil n/2\rceil}$, trivial on $K_n$. For any $\gamma\in K_1$ we have
\begin{equation}\label{eq:K1action}
\chi_{n,\Ad^*(\gamma)\lf}(x)=\chi_{n,\lf}(\gamma^{-1}x\gamma),\qquad
x\in K_{\lceil n/2\rceil}.
\end{equation}
Let $r\ge0$ with $n\ge2r$. The restriction of $\chi_{n,\lf}$ to
$K_{n-r}\subset K_{\lceil n/2\rceil}$, which we still denote by
$\chi_{n,\lf}$, depends only on $\cell=\lf+\pp^r\dLat$, and $\cell\mapsto\chi_{n,\lf}$ is a bijection between
$\dLat/\pp^r\dLat$ and the character group of $K_{n-r}/K_n$.

For $\GL_d$, these are Howe's dual blobs \cite{MR0492088}: his
shallow characters of $K_v=1+\pp^v\Mat_d(\OO)$ are those trivial on
$K_{2v}$, and their dual blobs are cosets of $\pp^{-v}\Mat_d(\OO)$ in
$\pp^{-2v}\Mat_d(\OO)$ \cite[\S2]{MR0492088}; up to the scaling in
\eqref{eq:chi}, these are our cells. Moy and Prasad \cite{MR1253198}
developed the corresponding filtration and dual-coset framework for general
reductive groups. At positive depth their unrefined minimal $K$-types use
nondegenerate characters of the quotients $G_{x,r}/G_{x,r+}$. We use
characters of the larger abelian congruence quotients $K_{n-r}/K_n$ and
impose no minimality condition on individual cells
(cf.~\zcref{ex:parahoric}).

Set
\[
e_{n,\cell}=e_{K_{n-r}}^{\chi_{n,\lf}}.
\]
We have
\begin{equation}\label{eq:cellformula}
\frac{e_{n,\cell}}{dg}(g)=
\frac{\one_{K_1}(g)}{\vol(K_n)}\int_\cell
\psi\bigl(-\varpi^{-n}\lf(\kappa(g))\bigr)\,d\nu_{\dLat}(\lf),\qquad g\in G.
\end{equation}
More generally, for a union $\cell\subset\dLat$ of cosets of $\pp^r\dLat$, let
$e_{n,\cell}=\sum_{\cell'}e_{n,\cell'}$, the sum over the cosets $\cell'$ of $\pp^r\dLat$
contained in $\cell$. Thus \eqref{eq:cellformula} continues to hold and in particular
$e_{n,\cell}$ is well defined (it does not depend on $r$, as long as $\cell$ is $\pp^r\dLat$ invariant).
In particular, $e_{n,\dLat}=e_{K_n}$.

By Fourier inversion on the finite abelian group $K_{n-r}/K_n$,
transported to $\Hk(G)$ by \cite[\S4.4]{MR2234120},
\begin{equation}\label{eq:partition}
e_{K_n}=\sum_{\cell\in \dLat/\pp^r\dLat}e_{n,\cell},
\end{equation}
and $e_{n,\cell}$, $\cell\in\dLat/\pp^r\dLat$, are mutually orthogonal idempotents.
By the Plancherel formula, for every union $\cell$ of cosets of $\pp^r\dLat$
\begin{equation}\label{eq:cellmass}
\vol(K_n)\int_{\Gtemp}\tr\pi(e_{n,\cell})\,d\mu_G(\pi)=\nu_{\dLat}(\cell).
\end{equation}
Compare \cite[Theorem~B]{MR2767510}, which describes the Plancherel measure
of the representations containing a given irreducible representation of a
compact open subgroup.

For $n,r$ as above, $\ast\in\{\sclab,\dslab\}$ and a cell $\cell\in \dLat/\pp^r\dLat$ let
\begin{equation}\label{eq:f}
\fmass{\ast}_{n,r}(\cell)=q^{rD}\vol(K_n)\int_{\Ghat{\ast}}\tr\pi(e_{n,\cell})\,d\mu_G(\pi)
\end{equation}
be the proportion of the Plancherel mass of $\cell$ carried by $\Ghat{\ast}$, whose
total mass is $q^{-rD}=\nu_{\dLat}(\cell)$ by \eqref{eq:cellmass}; we regard $\fmass{\ast}_{n,r}$
as a function on $\dLat$ constant on cells. Since $\tr\pi(e_{n,\cell})\ge0$ and
$\Gsc\subset\Gds$, \eqref{eq:cellmass} shows that
$0\le\fsc_{n,r}\le\fds_{n,r}\le1$, and \eqref{eq:partition} gives
\begin{equation}\label{eq:Pscdecomp}
\dens{\ast}{K_n}=\sum_{\cell\in \dLat/\pp^r\dLat}\nu_{\dLat}(\cell)\fmass{\ast}_{n,r}(\cell)=\int_\dLat \fmass{\ast}_{n,r}\,d\nu_{\dLat}.
\end{equation}
More generally, for $[M]\in\Lcl$ let
\begin{equation}\label{eq:fM}
\fmass{\ast}_{[M],n,r}(\cell)=q^{rD}\vol(K_n)\int_{\GhatL{\ast}{[M]}}\tr\pi(e_{n,\cell})\,d\mu_G(\pi)
\end{equation}
be the proportion of the Plancherel mass of $\cell$ carried by the class
$\GhatL{\ast}{[M]}$ of \eqref{eq:spectralsets}, so that
$\fmass{\ast}_{n,r}=\fmass{\ast}_{[G],n,r}$. For each $\ast$, these are
nonnegative functions on $\dLat$, constant on cells, with
$\sum_{[M]\in\Lcl}\fmass{\ast}_{[M],n,r}=1$ by \eqref{eq:cellmass}, and as before
\begin{equation}\label{eq:PscdecompM}
\densL{\ast}{[M]}{K_n}=\int_\dLat \fmass{\ast}_{[M],n,r}\,d\nu_{\dLat}.
\end{equation}

\subsection{Elliptic cells}\label{sec:elliptic}

For $\lf\in \dLat\cap\Ell{\gd}$ let $\erad(\lf)$ be the minimal integer $R\ge0$ with
$\lf+\pp^R\dLat\subset\Ell{\gd}$, the \emph{ellipticity radius} of $\lf$; it is
finite since $\Ell{\gd}$ is open (\zcref{sec:incidencemaps}), and
$\erad(\Ad^*(\gamma)\lf)=\erad(\lf)$ for $\gamma\in K_1$.

The following is a quantitative analogue of Moy--Prasad's criterion
for anisotropic unrefined minimal $K$-types \cite[Corollary~4.9]{MR1371680}.
(See \cite[Proposition~1]{MR522009}, \cite[Th\'eor\`eme~4.2]{MR760676}, \cite[Th\'eor\`eme~5]{MR396859},
\cite[Proposition~2.1]{MR0579176} for earlier special cases.)
In characteristic zero the same mechanism was used by Beuzart-Plessis,
together with the exponential map: in \cite{MR3463423}, in the form that the
Fourier transform of a function supported in the elliptic set is a cusp form,
to prove the existence of supercuspidal representations, and in
\cite[Appendix~D]{MR4325284}, where he shows that the characters attached to
suitable cosets of elliptic elements of $\gd$ are contained only in
supercuspidal representations. For $\SL_2$ see also \cite[\S6.1]{MR3753921}.

We state it for a single parabolic subgroup at a time, since this is
needed for the Levi classes.

\begin{proposition}[Cells and vanishing of Jacquet modules]\label{prop:elliptic}
There is a constant $C$ with the following property. Let $\lf\in\dLat$, let
$\bP$ be a proper parabolic subgroup and let $r\ge0$ and $n\ge 2r+C$ be
integers. Suppose that the cell $\lf+\pp^r\dLat$ is disjoint from
$\Pim{P}=\Ad^*(G)\uu_P^\perp$. Then
\begin{equation}\label{eq:nontrivialonU}
\chi_{n,\lf}|_{K_{n-r}\cap U_Q}\ne1\qquad\text{for every conjugate $\bQ$ of $\bP$,}
\end{equation}
and every irreducible representation $\pi$ of $G$ which contains the
character $\chi_{n,\lf}$ of $K_{n-r}$ satisfies $r_P(\pi)=0$.

In particular, if $\lf\in\dLat\cap\Ell{\gd}$ and $r\ge\erad(\lf)$, then every
irreducible representation of $G$ which contains the character $\chi_{n,\lf}$ of $K_{n-r}$ is
supercuspidal, and $\fsc_{n,r}(\lf)=1$.
In the unramified case (\zcref{sec:integraldata}) one can take $C=0$.
\end{proposition}

\begin{proof}
There are only finitely many $K_1$-conjugacy classes of proper parabolic subgroups,
and $\Ad^*(K_1)$ preserves $\dLat$.
Therefore, by \eqref{eq:K1action} it is enough to show that for a fixed proper parabolic subgroup $\bP$
there is a constant $C$ such that $\chi_{n,\lf}|_{K_{n-r}\cap U_P}\ne1$
for all $\lf\in \dLat$ and integers $r\ge0$, $n\ge2r+C$ such that $\lf+\pp^r\dLat$ is disjoint from $\uu_P^\perp$.

Let $\uu=\uu_P$, $\Lat_{\uu}=\Lat\cap\uu$ and
let $\rho:\gd\to\uu^*$ be the restriction. Since $\dLat=\Lat^\vee$, we have
$\rho(\dLat)=\Lat_{\uu}^\vee$. Fix a \npar{} $\xi$ of $U_P$. Since $\kappa\circ\xi(Y)=Y+O(Y^2)$, there is $c$ such that
$\kappa(\xi(Y))-Y\in\pp^{2a-c}\Lat$ for $Y\in\pp^a\Lat_{\uu}$,
$a\ge c$. In particular $\xi(\pp^b\Lat_{\uu})\subset K_b\cap U_P$,
$b=n-r$, and if $2b-c\ge n$, then
\[
\chi_{n,\lf}(\xi(Y))=\psi\bigl(\varpi^{-n}\lf(Y)\bigr),\qquad Y\in\pp^b\Lat_{\uu}.
\]
This character of $\pp^b\Lat_{\uu}$ is trivial if and only if
$\varpi^{-r}\rho(\lf)\in\Lat_{\uu}^\vee=\rho(\dLat)$, i.e., if and only if
$\lf+\pp^r\dLat$ meets $\uu^\perp$, which is excluded by assumption.
This proves \eqref{eq:nontrivialonU}.

Let $\Pi=\cInd_{K_b}^G\chi_{n,\lf}$ and $\chi=\chi_{n,\lf}$. For the
parabolic subgroup $P=MU$, the Mackey formula gives
\[
\Pi|_P\simeq\bigoplus_{g\in K_b\backslash G/P}\cInd_{P\cap g^{-1}K_bg}^P\chi^g,
\qquad\chi^g(x)=\chi(gxg^{-1}).
\]
The $U$-coinvariants of $\cInd_{P\cap g^{-1}K_bg}^P\chi^g$ vanish as soon as
$\chi^g$ is nontrivial on $U\cap g^{-1}K_bg$, i.e., as soon as $\chi$ is
nontrivial on $K_b\cap gUg^{-1}$. By \eqref{eq:nontrivialonU}, applied to the
parabolic subgroup $gPg^{-1}$, this holds for every $g$, so $\Pi_U=0$
(cf.\ \cite[4.5,~Corollary]{MR950594}). On the other hand, by Frobenius
reciprocity, an irreducible representation $\pi$ containing $\chi_{n,\lf}$
is a quotient of $\Pi$. Therefore, $r_P(\pi)=0$ by the exactness of the Jacquet functor.

If $\lf\in\Ell{\gd}$ and $r\ge\erad(\lf)$, then the cell $\lf+\pp^r\dLat$ is
disjoint from $\Par{\gd}$, hence from $\Pim{P}$ for every proper $\bP$. Thus
every proper Jacquet module of $\pi$ vanishes, i.e., $\pi$ is supercuspidal.
The assertion about $\fsc_{n,r}$ follows from \eqref{eq:f} and \eqref{eq:cellmass}.

In the unramified case, put $K=\cG(\OO)$. By the Iwasawa decomposition, $K$
acts transitively on each conjugacy class of parabolic subgroups, and
conjugation by $K$ preserves $\Lat$, $\dLat$, every $K_b$ and, in the integral
coordinates, the cell characters when $2b\ge n$; so \eqref{eq:K1action} holds
for $\gamma\in K$, and it suffices to treat a standard $\bP$.
For any standard proper $\bP$, choose a cocharacter of $\cT$
defining $\bP$ and use the integral big-cell construction of
\zcref{sec:integraldata} for this cocharacter. The lattice
$\Lat_\uu=\Lie(\cU)(\OO)$ is a direct summand of $\Lat$, and the integral
parametrization $\xi_+$ satisfies $\kappa(\xi_+(Y))-Y\in\pp^{2a}\Lat$ for
$Y\in\pp^a\Lat_\uu$, $a\ge1$. Thus the argument above works with $c=0$, and
$n\ge2r$ suffices.
\end{proof}

\begin{remark}[Existence of an irreducible representation containing a type]\label{rem:Ktype}
For a compact open subgroup $K$ of $G$ and an irreducible representation
$\rho$ of $K$ there is an irreducible representation $\pi$ of $G$
containing $\rho$. Indeed, $\cInd_K^G\rho$ is finitely generated -- it is generated by
the functions supported on $K$ -- so it has an irreducible quotient $\pi$, and $\Hom_K(\rho,\pi)=\Hom_G(\cInd_K^G\rho,\pi)\ne0$.
In particular, the representations in \zcref{prop:elliptic} exist.
\end{remark}

\section{The Plancherel limit}\label{sec:statements}

For the standard pair $(\kappa,\Lat)$ fixed above, write $K_n$ for its
principal congruence tower. The exponents $\eexp$ and $\eexpall$ are
defined in \eqref{eq:tau} and \eqref{eq:tauall}. The proof is given in
\zcref{sec:proofsgeneral}.

\begin{theorem}[Plancherel mass of the supercuspidal representations and of the Levi classes]\label{thm:main}\label{thm:levi}
We have $0<\elld{\dLat}\le1$ and, for $\ast\in\{\sclab,\dslab\}$,
\begin{equation}\label{eq:main}
\dens{\ast}{K_n}=\elld{\dLat}+O(q^{-\eexp n/4})\qquad(n\to\infty),
\end{equation}
and hence $0\le\dsd{K_n}-\scd{K_n}=O(q^{-\eexp n/4})$.
If $Z$ is compact, let $\fdeg{\pi}$ denote the formal degree of
$\pi\in\Gds$, viewed as a Haar measure on $G$. Then, for
$\ast\in\{\sclab,\dslab\}$,
\begin{equation}\label{eq:formal}
\sum_{\pi\in\Ghat{\ast}}\fdeg{\pi}(K_n)\dim\pi^{K_n}=\elld{\dLat}+O(q^{-\eexp n/4}).
\end{equation}
More generally, for every $[M]\in\Lcl$ and $\ast\in\{\sclab,\dslab\}$,
\begin{equation}\label{eq:mainlevi}
\densL{\ast}{[M]}{K_n}=\levid{[M]}{\dLat}+O(q^{-\eexpall n/4})\qquad(n\to\infty).
\end{equation}
The same estimate \eqref{eq:mainlevi} holds for the joint event
$\Lsup{\sclab}(\pi)=\Lsup{\dslab}(\pi)=[M]$. Moreover,
\begin{equation}\label{eq:exceptional}
\vol(K_n)\int_{\{\pi\in\Gtemp:\Lsup{\sclab}(\pi)\ne\Lsup{\dslab}(\pi)\}}
\dim\pi^{K_n}\,d\mu_G(\pi)=O(q^{-\eexpall n/4}).
\end{equation}
\end{theorem}

The implied constants depend on $F$, $\bG$ and $\Lat$; the exponents
$\eexp/4\ge\eexpall/4\ge1/4$ depend only on $\bG$ (see \zcref{rem:sharpness}).

For $[M]=[G]$ we have $\GhatL{\ast}{[G]}=\Ghat{\ast}$ and
$\levid{[G]}{\dLat}=\elld{\dLat}$, so that \eqref{eq:mainlevi} for $[M]=[G]$ is
\eqref{eq:main} with $\eexpall$ in place of $\eexp$.
The exponent in \eqref{eq:main} can be larger because the
$[G]$ estimate uses only maximal proper parabolic subgroups.
The density $\levid{[M]}{\dLat}$ may vanish for $[M]\ne[G]$ (\zcref{ex:pgl2}), whereas
$\levid{[G]}{\dLat}>0$.
By \eqref{eq:labels}, the integral in \eqref{eq:exceptional} is over the
tempered representations whose discrete series support $(\bL,\delta)$ has
$\delta$ not supercuspidal: parabolic induction from non-supercuspidal
discrete series of all Levi subgroups, including $G$ itself, carries
asymptotically negligible Plancherel mass. A localized version of
\eqref{eq:mainlevi} is given in \zcref{cor:localized}.

For a smooth model $\cG$ as in \zcref{lem:model}, with
$\Lat=\Lie(\cG)(\OO)$ and $\vol(\cG(\OO))=1$, substituting
$\vol(K_n)^{-1}=\abs{\cG(k)}q^{(n-1)D}$ into \zcref{thm:main} gives
\[
\int_{\Gsc}\dim\pi^{K_n}\,d\mu_G(\pi)=\elld{\dLat}\abs{\cG(k)}q^{(n-1)D}+O\bigl(q^{(D-\eexp/4)n}\bigr).
\]

The elliptic density $\elld{\dLat}$ depends on the lattice and not
only on $\bG$; it is studied in \zcref{sec:positivity}, where we also
determine when it equals $1$ (\zcref{thm:endpoint}).

\section{Parabolic incidence maps}\label{sec:incidence}

Using the quantitative Sard theorem of \zcref{sec:sard} (\zcref{thm:sard}), we
construct two decreasing families of \bnets{} $\Badmax_h\subset\Bad_h\subset\dLat$,
each a union of cosets of $\pp^{2h}\dLat$, with
$\nu_{\dLat}(\Badmax_h)\ll q^{-\eexp h}$ and $\nu_{\dLat}(\Bad_h)\ll q^{-\eexpall h}$.
The first is attached to the maximal proper parabolic subgroups: every coset of
$\pp^{2h}\dLat$ outside $\Badmax_h$ consists entirely of elliptic forms or
entirely of nonelliptic forms. The second is attached to all proper parabolic
subgroups: every coset of $\pp^{2h}\dLat$ outside $\Bad_h$ is, for each $\bP$,
contained in the image $\Pim{P}$ of the incidence map of $\bP$ or disjoint from it. Whenever a coset
outside the nets is contained in such an image, the incidence map has a
lift with quantitatively controlled inverse differential. This is the
geometric input for the cellwise comparison of \zcref{prop:density}.

\subsection{The pairing and the geometric alternatives}\label{sec:levi}
Fix a semistandard parabolic subgroup $\bP=\bM\bU$ and its opposite
semistandard $\bP^-=\bM\bU^-$ with common semistandard Levi subgroup $\bM$. Write
$\g=\uu\oplus\m\oplus\uu^-$, with $\uu^+=\uu$ and $\bU^+=\bU$,
and put $d_P=\dim\uu=\dim\uu^-$.
For $\lf\in\gd$ define
\begin{equation}\label{eq:levipairing}
\Bf_\lf(X,Y)=\lf([X,Y]),\qquad X\in\uu,\ Y\in\uu^-.
\end{equation}
For $\lf\in\uu^\perp$, let $Q_{\lf,P}:\uu^-\to\uu^*$ and $Q_{\lf,P}^*:\uu\to(\uu^-)^*$ be the corresponding maps.
Suppressing $P$ when it is understood from the context, we have
\begin{equation}\label{eq:pairing}
(Q_\lf Y)(X)=(Q_\lf^*X)(Y)=\Bf_\lf(X,Y).
\end{equation}
Identify $\m^*$ with the forms vanishing on $\uu\oplus\uu^-$, and let
$\Omega_M\subset\m^*$ be the open, possibly empty, set on which $Q_\lf$ is invertible
(it is Zariski open, hence also open for the topology of $F$).

For each $\bP\in\Pssall$, endow $\operatorname{Hom}_F(\uu^*,\uu^-)$ with a norm once and for all.
The choice affects only the constant $c$ in the following proposition.

\begin{proposition}[Geometric alternatives outside the \bnets{}]\label{prop:geometric}
There are a compact subset $\mathcal C\subset G$, compact sets
$\mathcal A_P\subset\uu_P^\perp$ for $\bP\in\Pssall$, an integer $c$,
and decreasing sets $\Badmax_h\subset\Bad_h\subset\dLat$, $h\ge1$, satisfying the following conditions.
\begin{enumerate}[label=\textup{(\arabic*)},ref=\arabic*]
\item\label{item:geo-cosets} Each $\Bad_h$ and each $\Badmax_h$ is a union of cosets of $\pp^{2h}\dLat$.
\item\label{item:geo-measure} $\bigcap_{h\ge1}\Bad_h=\Crv\cap\dLat$ and
$\bigcap_{h\ge1}\Badmax_h=\Crvmax\cap\dLat$. Moreover
$\nu_{\dLat}(\Bad_h)\ll q^{-\eexpall h}$ and $\nu_{\dLat}(\Badmax_h)\ll q^{-\eexp h}$.
\item\label{item:geo-alt} Let either $\lf\in\dLat\setminus\Bad_h$ and $\bP\in\Pstd$, or
$\lf\in\dLat\setminus\Badmax_h$ and $\bP\in\Pc$. Then either $\lf+\pp^{2h}\dLat$ is disjoint from $\Pim{P}$, or
$\lf+\pp^{2h}\dLat\subset\Pim{P}$, in which case there are a semistandard conjugate
$\bP'$ of $\bP$, $g\in\mathcal C$ and $\nf\in\mathcal A_{P'}$ with
\[
\lf=\Ad^*(g)\nf,\qquad \norm{Q_{\nf,P'}^{-1}}\le q^{h+c}.
\]
\item\label{item:geo-ell} In particular, for every $\lf\in\dLat\setminus\Badmax_h$ either
$\lf+\pp^{2h}\dLat\subset\Ell{\gd}$ or $\lf+\pp^{2h}\dLat\subset\Par{\gd}$.
\item\label{item:geo-tube} For $r\ge2$,
\[
\nu_{\dLat}((\partial\Ell{\gd}+\pp^r\dLat)\cap\dLat)\le\nu_{\dLat}((\Crvmax+\pp^r\dLat)\cap\dLat)\ll q^{-\eexp r/2}
\text{ and }\nu_{\dLat}((\Crv+\pp^r\dLat)\cap\dLat)\ll q^{-\eexpall r/2}.
\]
\end{enumerate}
\end{proposition}

In the unramified case we will make the statement more precise in \zcref{prop:geometricunram}.

\subsection{The differential and the coordinate Jacobian}\label{sec:incidencejacobian}
For each $\bP\in\Pssall$, fix a normalized parametrization
$\xi_-:\uu^-\to U^-$ which is an isomorphism of $F$-varieties (\zcref{sec:standardgroups}).
We identify $\Inc_P$, defined as in \eqref{eq:Y}, with $\Inc_{P_1}$ for the
standard conjugate $\bP_1\in\Pstd$ of $\bP$, and write $\incm$ for the projection
$\incm_{P_1}$ when $\bP$ is understood. Let $\Incc{P}$ be the inverse image
in $\Inc_P$ of the big cell
$U^-P/P\subset\Fl_P$. The trivialization
\begin{equation}\label{eq:trivialization}
U^-\times\uu^\perp\longrightarrow\Incc{P},\qquad
(u,\nf)\longmapsto u\cdot(\bP,\nf)
\end{equation}
and $\xi_-$ give a polynomial parametrization
\begin{equation}\label{eq:bigcellchart}
\phi_P:\uu^-\times\uu^\perp\longrightarrow\Incc{P},\qquad
(y,\nf)\longmapsto\xi_-(y)\cdot(\bP,\nf).
\end{equation}
The incidence map in these coordinates is
\[
f_P=\incm\circ\phi_P:\uu^-\times\uu^\perp\longrightarrow\gd,\qquad
(y,\nf)\longmapsto\Ad^*(\xi_-(y))\nf.
\]
Fix bases of $\uu^-$, $\uu^\perp$ and $\uu^*$. Order the basis of
$\uu^-\oplus\uu^\perp$ with the basis of $\uu^\perp$ first, and use for $\gd$
the basis of $\uu^\perp$ followed by lifts of the basis of $\uu^*$ (via the
exact sequence $0\to\uu^\perp\to\gd\to\uu^*\to0$). A linear map
$\uu^-\oplus\uu^\perp\to\gd$ which is the identity on $\uu^\perp$ and induces
$Q:\uu^-\to\uu^*$ on the quotients then has determinant $\det Q$ in these
bases.\footnote{In characteristic zero, a
non-degenerate invariant symmetric bilinear form on $\g$ identifies $\gd$
with $\g$ and $\uu^*$ with $\uu^-$, so that $Q_\nf$ becomes an endomorphism
of $\uu^-$ and $\det Q_\nf$ is defined without further choices. In positive
characteristic such a form need not exist.} After fixing bases,
$\nf\mapsto\det Q_\nf$ is a polynomial, well defined up to a nonzero constant.

\begin{lemma}[Differential and Jacobian of the incidence map]\label{lem:levi}
In the trivialization \eqref{eq:trivialization},
\begin{equation}\label{eq:dp}
(d\incm_P)_{(\bP,\nf)}(Y,\mu)=\mu+\ad^*(Y)\nf,
\qquad Y\in\uu^-,\ \mu\in\uu^\perp.
\end{equation}
Moreover,
\begin{equation}\label{eq:jacobianfiber}
\det Q_\nf=\det Q_{\nf|_\m},
\end{equation}
and hence
\begin{equation}\label{eq:levi}
(d\incm_P)_{(\bP,\nf)}\text{ is invertible}
\iff Q_\nf\text{ is invertible}\iff\nf|_\m\in\Omega_M.
\end{equation}
Finally, in the bases above,
\begin{equation}\label{eq:detfield}
\det df_P(y,\nf)=\det Q_\nf\qquad\text{for all }(y,\nf)\in\uu^-\times\uu^\perp.
\end{equation}
\end{lemma}

\begin{proof}
Differentiating the coadjoint action gives \eqref{eq:dp}. The map is
the identity on $\uu^\perp$ and induces $Q_\nf$ on the quotients.
For \eqref{eq:jacobianfiber}, let $\g=\bigoplus_k\g_k$ be the weight
decomposition for a cocharacter $\eta\in X_*(\bA_0)$ defining $\bP$, so that
$\uu=\bigoplus_{k>0}\g_k$, $\m=\g_0$ and $\uu^-=\bigoplus_{k>0}\g_{-k}$.
If $X\in\g_i$, $Y\in\g_{-j}$ with $i>j>0$, then $[X,Y]\in\g_{i-j}\subset\uu$,
so $\Bf_\nf(X,Y)=0$. Thus $Q_\nf$ is block triangular
in bases adapted to the weight decompositions, with weights ordered
by increasing $i$ on both sides. Its diagonal block $(Q_\nf)_{ii}$
is the restriction to $\g_{-i}$ followed by projection to $\g_i^*$.
These blocks only involve $\nf$ on $[\g_i,\g_{-i}]\subset\m$. Hence
$\det Q_\nf=\prod_{i>0}\det((Q_\nf)_{ii})=\det Q_{\nf|_\m}$ in these bases,
proving \eqref{eq:jacobianfiber} and hence \eqref{eq:levi}. The same
argument shows that $Q_\nf$ is invertible if and only if each pairing
$\g_i\times\g_{-i}\to F$, $(X,Y)\mapsto\nf([X,Y])$, is non-degenerate.

For the last assertion fix $y$, put $u=\xi_-(y)$. Writing multiplication by $u$ in the $\uu^-$ coordinates, we obtain a polynomial automorphism
$\theta_y$ of $\uu^-$ such that $\xi_-(y')=u\xi_-(\theta_y(y'))$ and in particular $\theta_y(y)=0$. Thus
\[
f_P(y',\nf')=\Ad^*(u)f_P(\theta_y(y'),\nf').
\]
Keeping $y$ and $u$ fixed and differentiating at $(y,\nf)$ gives by the chain rule
\[
df_P(y,\nf)=\Ad^*(u)\circ df_P(0,\nf)\circ(T_y\oplus1)
\]
where $T_y=d(\theta_y)_y$.
Since $\theta_y$ is a polynomial automorphism for every geometric point $y$, $\det T_y$ is a polynomial in $y$ which has no zeros over $\bar F$
and it equals $1$ at $y=0$. Hence $\det T_y\equiv1$.
Since $\bU^-$ is unipotent, $\det\Ad^*(u)=1$. The quotient calculation
above shows that $\det df_P(0,\nf)=\det Q_\nf$ in the bases
preceding the lemma, since $d\xi_-$ is the identity at $0$. This proves \eqref{eq:detfield}.
\end{proof}

As a polynomial in the fiber coordinate, $\det Q_{\bullet,P}$ is
either identically zero or homogeneous of degree $d_P$. The latter
occurs exactly when $\Omega_M\ne\emptyset$, equivalently when
$\incm_P$ has a noncritical point.

\subsection{Noncritical lifts}\label{sec:noncritical}
A noncritical point of an incidence map can be moved into the dual of the
Levi subalgebra. In characteristic zero this is the fact, recalled in
\zcref{sec:ellipticforms}, that a regular semisimple element of
$\Lie(\bP)(F)$ is conjugate under $U_P$ to an element of $\m$. The following
argument works in every characteristic. Recall that $\m^*$ is identified with the
forms vanishing on $\uu\oplus\uu^-$ (\zcref{sec:levi}).

\begin{lemma}[Normal form of a noncritical lift]\label{lem:normalform}
Let $\bP=\bM\bU\in\Pssall$, let $\nf\in\uu^\perp$ with $Q_\nf$ invertible, and
let $\nf_\m\in\m^*$ be the form which agrees with $\nf$ on $\m$. Then
$\nf=\Ad^*(u)\nf_\m$ for some $u\in U$.
\end{lemma}

\begin{proof}
Put $\mu=\nf_\m$; by \eqref{eq:jacobianfiber}, $Q_\mu$ is invertible. For
$u\in U$ the form $\Ad^*(u)\mu=\mu\circ\Ad(u^{-1})$ vanishes on $\uu$ and
agrees with $\mu$ on $\m$, since $\Ad(u^{-1})$ preserves $\uu$ and
$\Ad(u^{-1})X-X\in\uu$ for $X\in\m$. Hence it suffices to show that the
analytic map
\[
\Psi:U\longrightarrow(\uu^-)^*,\qquad\Psi(u)=(\Ad^*(u)\mu)|_{\uu^-}
\]
is surjective. We have $\Psi(1)=0$, and by \zcref{sec:notationbasic} and
\eqref{eq:pairing} the differential of $\Psi$ at $1$ is
$X\mapsto(\ad^*(X)\mu)|_{\uu^-}=-Q_\mu^*X$, which is invertible. By the inverse
function theorem, $\Psi(U)$ contains a neighborhood $W$ of $0$.

Let $\eta\in X_*(\bA_0)$ be a cocharacter defining $\bP$ and $\bM$ as in
\zcref{sec:parabolics}, and put $t=\eta(\varpi)$. Then $t$ centralizes $M$ and
normalizes $U$ and $U^-$, and $\Ad(t^{-1})$ acts on the weight space $\g_{-k}$,
$k>0$, by $\varpi^{k}$. Hence $\Ad^*(t)\mu=\mu$, and
$\Ad^*(t^j)\nu=\nu\circ\Ad(t^{-j})\to0$ as $j\to\infty$ for every
$\nu\in(\uu^-)^*$. Moreover, for $u\in U$ and $j\in\Z$,
\[
\Psi(t^jut^{-j})=(\Ad^*(t^j)\Ad^*(u)\Ad^*(t^{-j})\mu)|_{\uu^-}=\Ad^*(t^j)\Psi(u).
\]
Given $\nu\in(\uu^-)^*$, choose $j$ with $\Ad^*(t^j)\nu\in W$, say
$\Ad^*(t^j)\nu=\Psi(v)$ with $v\in U$. Then $\nu=\Psi(t^{-j}vt^j)$.
\end{proof}

\begin{proposition}[Incidence off the critical values]\label{prop:upset}
Let $\lf\in\gd\setminus\Crv$. Then for every parabolic subgroup $\bP$ of $\bG$,
\begin{equation}\label{eq:upset}
\lf\in\Pim{P}\iff\Lform{\lf}\preceq[M_P],
\end{equation}
where $\Lform{\lf}$ is the least class of \zcref{lem:meet}. In particular, for
$\lf\notin\Crv$ the condition $\lf\in\Pim{P}$ depends only on the class $[M_P]$ of
the Levi subgroups of $\bP$, and $\Pim{P}\setminus\Crv=\Pim{[M_P]}\setminus\Crv$.
\end{proposition}

\begin{proof}
If $\lf\in\Pim{P}$ then $\Lform{\lf}\preceq[M_P]$ by the definition of $\Lform{\lf}$.
Conversely, let $\bQ$ be a parabolic subgroup with $\lf\in\Pim{Q}$ and
$[M_Q]=\Lform{\lf}$, and suppose that $[M_Q]\preceq[M_P]$. We may assume that $\bQ$
is proper, since otherwise $\bP=\bG$. Both sides of \eqref{eq:upset} and the set $\Crv$ are invariant under
conjugation of $\lf$ and of $\bP$, so we may assume that $\bQ$ is semistandard, with
Levi subgroup $\bL=\bM_Q$, that $\lf\in\uu_Q^\perp$, and that $\bL$ is contained in a Levi subgroup
$\bM$ of $\bP$. Since $\lf\notin\Crv$, the point
$(\bQ,\lf)$ of $\Inc_Q$ is not critical, i.e., $Q_{\lf,Q}$ is invertible by
\eqref{eq:levi}. By \zcref{lem:normalform}, after conjugating $\lf$ by an element of $U_Q$ we may assume that $\lf$ vanishes on
$\uu_Q\oplus\uu_Q^-$, which is the sum of the weight spaces of $A_L$ on $\g$ with nonzero weight.
Now $A_M$ commutes with $L\subset M$, so $A_M\subset Z_G(A_L)=L$, and
therefore $A_M\subset A_L$. Since $\uu_P$ is a sum of weight spaces of $A_M$
with nonzero weight, it is contained in $\uu_Q\oplus\uu_Q^-$. Thus $\lf$ vanishes on
$\uu_P$, i.e., $\lf\in\Pim{P}$.
\end{proof}

\begin{remark}[Incidence profiles at critical values]
The equivalence \eqref{eq:upset} in fact holds for every $\lf$.
If the incidence map of a parabolic $\bQ$ realizing $\Lform{\lf}$
has a noncritical point, the noncritical forms are dense in
$\uu_Q^\perp$. The argument of \zcref{prop:upset} applied to these
forms gives the assertion by closedness of $\Pim{P}$.
The incidence profile and the Levi order both factor over a direct
product. By \zcref{lem:split,prop:peculiar}, it remains to treat odd
orthogonal factors in characteristic two. For such a factor use
the model $\mathfrak{so}(V)=\bigwedge^2 V$ of
\zcref{lem:oddorthogonal}, and write $A$ for the alternating form
corresponding to $\lf$ and $\beta$ for the polar form of $V$.
For an isotropic flag ending in $E$, incidence amounts to
$A(E,E^\perp)=0$ and stability of the flag under the endomorphism
$C$ of $E$ defined by $A(e,v)=\beta(Ce,v)$. The first condition
comes from the unipotent radical of the stabilizer of $E$. The
second is the trace-pairing condition in its $\GL(E)$ Levi factor.
It is easy to see that these conditions are preserved under enlargement of the Levi class.

The \emph{noncritical Levi class} $\Lnc{\lf}$ is the least class $[M_Q]$
for which $\lf$ has a noncritical lift to $\Inc_Q$, including $\bG$
with its identity incidence map. It is the class of the centralizer
of the largest split torus fixing $\lf$ and acting trivially on its
infinitesimal stabilizer $\g_\lf$. Thus, for every parabolic $\bP$,
\[
\lf\text{ has a noncritical lift to }\Inc_P
\quad\Longleftrightarrow\quad
\Lnc{\lf}\preceq[M_P].
\]
For $\GL_d$, the noncritical class has one block per primary component
of the characteristic polynomial, whereas the ordinary class has one
block per irreducible factor counted with multiplicity. For a single
nilpotent Jordan block the two classes are $[G]$ and $[T]$, respectively.
\end{remark}

\subsection{Construction of the \bnets{}}\label{sec:bad}
\begin{proof}[Proof of \zcref{prop:geometric}]
If $\Pc=\emptyset$, then $\Pstd=\emptyset$; take $\Bad_h=\Badmax_h=\emptyset$; all forms are elliptic.
Otherwise fix $\bP\in\Pstd$. The big cells of the semistandard conjugates of $\bP$
cover $\Inc_P$. Indeed, by the Bruhat decomposition
$G=\bigcup_{\dot w}P_0\dot wP$, together with
$\dot w^{-1}P_0\dot w\subset U_P^-P$, we have $G=\bigcup_{\dot w}\dot wU_P^-P$,
$\dot w\in N_G(A_0)$. Hence the translates $\dot wU_P^-P/P$ cover
$\Fl_P$, and $\dot wU_P^-P/P$ is the big cell of the semistandard
conjugate $\dot w\bP\dot w^{-1}$.

Let $x\in\Inc_{P,\dLat}=\incm_P^{-1}(\dLat)$, a compact open subset of $\Inc_P$. Choose a
semistandard conjugate $\bP_x$ of $\bP$ whose big cell contains $x$, write
$x=\phi_{P_x}(y_x,\nf_x^0)$, and fix lattices $H^-\subset\uu_{P_x}^-$ and
$H^0\subset\uu_{P_x}^\perp$ with bases. The map
$(z,\mu)\mapsto f_{P_x}(y_x+z,\nf_x^0+\mu)$ is polynomial with constant term
$\incm_P(x)\in\dLat$. Replacing $(H^-,H^0)$ by $(\pp^sH^-,\pp^sH^0)$, with the
rescaled bases, multiplies the coefficient of each monomial of degree $k$
by $\varpi^{sk}$. As there are finitely many coefficients, for $s$ large
all coefficients of positive degree lie in $\dLat$. Then
$(z,\mu)\mapsto\phi_{P_x}(y_x+z,\nf_x^0+\mu)$ maps $\pp^sH^-\oplus\pp^sH^0$
into $\Inc_{P,\dLat}$, and its image is an open neighborhood of $x$, because
$\phi_{P_x}$ is an open immersion. By compactness, finitely many of these
charts cover $\Inc_{P,\dLat}$; we write them as
\[
\psi_i:\Lat_i\longrightarrow\Inc_{P,\dLat},\qquad
\psi_i(z,\mu)=\phi_{P_i}(y_i+z,\nf_i^0+\mu),\qquad i\in I_P,
\]
where $\bP_i$ is a semistandard conjugate of $\bP$, $\Lat_i=H_i^-\oplus H_i^0$ is a product of lattices in
$\uu_{P_i}^-\oplus\uu_{P_i}^\perp$, and the polynomial maps
$f_i=\incm_P\circ\psi_i:\Lat_i\to\dLat$ have coefficients in $\dLat$.
Thus the $\psi_i$, $i\in I_P$, form an \iatlas{} of $\incm_P$ over $\dLat$ in the sense of \zcref{sec:bnets}.
Computing $\det df_i$ in the bases of $\Lat_i$ and $\dLat$, and $\det Q$ in
fixed bases of $\uu_{P_i}^-$ and $\uu_{P_i}^*$, \eqref{eq:detfield} gives
constants $a_i\in F^\times$ such that
\begin{equation}\label{eq:atlasjacobian}
\det df_i(z,\mu)=a_i\det Q_{\nf_i^0+\mu,P_i}.
\end{equation}

Let $\Bad_h(P)$ be the \bnets{} \eqref{eq:badf} for this atlas of $\incm_P$, and put
\begin{equation}\label{eq:netsunion}
\Bad_h=\bigcup_{\bP\in\Pstd}\Bad_h(P),\qquad
\Badmax_h=\bigcup_{\bP\in\Pc}\Bad_h(P).
\end{equation}
By \zcref{thm:sard}, applied to each $\incm_P$, the sets $\Bad_h(P)$ are unions of
cosets of $\pp^{2h}\dLat$ and decrease to $\incm_P(\Crit\incm_P)\cap\dLat$. Since the
unions in \eqref{eq:netsunion} are finite, $\Bad_h$ and $\Badmax_h$ decrease to
$\Crv\cap\dLat$ and $\Crvmax\cap\dLat$, respectively. By the same theorem, a coset of
$\pp^{2h}\dLat$ which is not contained in $\Bad_h(P)$ is either contained in
$\Pim{P}=\incm_P(\Inc_P)$ or disjoint from it. This gives the dichotomy in
\zcref{item:geo-alt}. Since $\Par{\gd}=\bigcup_{\bP\in\Pc}\Pim{P}$ by \eqref{eq:N}, it also gives
\zcref{item:geo-ell}:
\begin{equation}\label{eq:dichotomy}
\lf+\pp^{2h}\dLat\not\subset\Badmax_h\implies
\text{either }\lf+\pp^{2h}\dLat\subset\Ell{\gd}
\text{ or }\lf+\pp^{2h}\dLat\subset\Par{\gd}.
\end{equation}
The construction includes neighborhoods of critical values that
also have noncritical lifts.

Suppose now that $\lf\in\Pim{P}\cap\dLat$ lies outside $\Bad_h(P)$. Choose any $x\in \incm_P^{-1}(\lf)$ and a chart
containing $x$, say $x=\psi_i(z,\mu)$ with $i\in I_P$. Put $\bP'=\bP_i$,
$\nf=\nf_i^0+\mu$ and $g=\xi_-(y_i+z)$, where $\xi_-$ is the parametrization of $U_{P_i}^-$. Then
$\lf=\Ad^*(g)\nf$, with $g$ and $\nf$ in the fixed compact sets
\[
\mathcal C=\bigcup_i\xi_-(y_i+H_i^-),\qquad
\mathcal A_{P'}=\bigcup_{i:P_i=P'}(\nf_i^0+H_i^0),
\]
the unions being over the charts of all $\bP\in\Pstd$.
By \eqref{eq:sardchart} and \eqref{eq:atlasjacobian},
$\abs{\det Q_\nf}>\abs{a_i}^{-1}q^{-h}$.
The entries of $Q_\nf$ are bounded on these compact sets, so Cramer's
rule gives $\norm{Q_\nf^{-1}}\le q^{h+c}$, with $c$ depending on
the fixed reference norms and compact families. This completes the proof of
\zcref{item:geo-alt}.

For each $\bP\in\Pstd$ write $d_P=\dim\uu_P$. The estimate for a
single parabolic is
\[
\nu_{\dLat}(\Bad_h(P))\ll
\begin{cases}
q^{-(1+1/d_P)h} & \text{if }\det Q_{\bullet,P}\not\equiv0,\\
q^{-h} & \text{if }\det Q_{\bullet,P}\equiv0.
\end{cases}
\]
Indeed, \eqref{eq:sardmeasure} bounds this measure by
$q^{-h}\sum_{i\in I_P}\vol(S_h(f_i))$. Each source volume is at most one.
If $\det Q_{\bullet,P}\not\equiv0$, then \eqref{eq:atlasjacobian} is a
nonzero polynomial of degree $d_P$ in each chart, and
\zcref{lem:sublevel} bounds the source volumes by $O(q^{-h/d_P})$.

When $\bG$ is not \peculiar{}, $\Omega_{M_P}$ is nonempty for every $\bP\in\Pc$ by
\eqref{eq:levi} and \zcref{def:peculiar}, and in fact for every
$\bP\in\Pstd$ by \zcref{prop:peculiar}. Since $d_P\le\dG{}$ for $\bP\in\Pc$ and
$d_P\le\dGall{}$ for $\bP\in\Pstd$, the improved bounds
$\nu_{\dLat}(\Badmax_h)\ll q^{-(1+1/\dG{})h}$ and
$\nu_{\dLat}(\Bad_h)\ll q^{-(1+1/\dGall{})h}$ apply.
Otherwise the bound $O(q^{-h})$ suffices.
These are precisely the exponents in \eqref{eq:tau} and \eqref{eq:tauall}.

Finally, \zcref{item:geo-tube} follows from \zcref{item:geo-measure} and
\eqref{eq:bntube}, applied to each $\incm_P$, since $\partial\Ell{\gd}\subset\Crvmax$ by
\eqref{eq:boundary}.
\end{proof}

\subsection{The unramified case}\label{sec:integralatlas}
In the unramified case (\zcref{sec:integraldata}) the charts can be taken integral, and their number is bounded in terms of the root datum. Let
$\bP\in\Pc$, let $\dot w$, $\bP_w$, $\cP_w$ and $\cU^\pm_w$,
$w\in W_k/W_{k,P}$, be as in \zcref{sec:integraldata}, let $\xi^w_\pm$ be the
conjugates of $\xi_\pm$ by $\dot w$, and put
$\Lat^w_-=\Lie(\cU^-_w)(\OO)$. Consider the maps
\begin{equation}\label{eq:integralchart}
f_{P,w}:\Lat^w_-\times(\uu_{P_w}^\perp\cap \dLat)\to \dLat,\qquad
f_{P,w}(y,\nf)=\Ad^*(\xi^w_-(y))\nf,
\end{equation}
and the parametrizations $\psi_{P,w}(y,\nf)=\xi^w_-(y)\cdot(\bP_w,\nf)$ of $\Inc_P$.

\begin{lemma}[An \iatlas{} with a uniform number of charts]\label{lem:integralatlas}
The maps $f_{P,w}$ are polynomial with coefficients
in $\dLat$, and the $\psi_{P,w}$, for $\bP\in\Pc$ and $w\in W_k/W_{k,P}$, form
an \iatlas{} of $\incm$ over $\dLat$ with at most $\Nss$ charts, with
$\Lat_i=\Lat^w_-\oplus(\uu_{P_w}^\perp\cap \dLat)$. Consequently, the
\bnets{} \eqref{eq:badf} of this atlas, which we denote by $\Badmax_h$, satisfy
\begin{equation}\label{eq:badexplicit}
\nu_{\dLat}(\Badmax_h)\le \Nss\,q^{-h}\qquad(h\ge1).
\end{equation}
\end{lemma}

\begin{proof}
The map \eqref{eq:integralchart} is polynomial with integral
coefficients because $\cU^-_w$ acts on $\Lie(\cG)$ by an action of
$\OO$-group schemes and $\xi^w_-$ is an isomorphism of $\OO$-schemes
mapping $\Lat^w_-$ onto $\cU^-_w(\OO)$, and
$\psi_{P,w}$ is the trivialization \eqref{eq:trivialization} of $\Inc_P$ over the big cell of $\bP_w$, hence a parametrization of an open subset of $\Inc_P$. If
$(\bQ,\lf)\in \incm_P^{-1}(\dLat)$ then $\bQ=u\bP_wu^{-1}$ with $u\in\cU_w^-(\OO)$ by
\eqref{eq:integralcover}, and $\nf=\Ad^*(u^{-1})\lf\in \dLat$ vanishes on
$\uu_{P_w}$. The number of charts is at most $\Nss$, since the $\bP_w$ are
pairwise distinct elements of $\Pcss$ (\zcref{sec:integraldata}), and those
attached to different $\bP\in\Pc$ are not conjugate. The last assertion is \zcref{lem:tube} applied to the
$f_{P,w}$, which have coefficients in $\OO$.
\end{proof}

In the following proposition and its proof, the norms on $\g$, $\gd$ and their subspaces
are those defined by $\Lat$, by its direct summands $\Lie(\cU_P)(\OO)$, $\Lie(\cM_P)(\OO)$
and $\Lie(\cU_P^-)(\OO)$ for standard $\bP$, and by their duals.

\begin{proposition}[Geometric alternatives in the unramified case]\label{prop:geometricunram}
In the unramified case (\zcref{sec:integraldata}), the \bnets{} $\Badmax_h$ of
\zcref{lem:integralatlas} satisfy the assertions of \zcref{prop:geometric}
concerning $\Badmax_h$ and the maximal parabolic subgroups, with
$\nu_{\dLat}(\Badmax_h)\le \Nss q^{-h}$ for all $h\ge1$, $\mathcal C\subset\cG(\OO)$,
$\mathcal A_P=\uu_P^\perp\cap\dLat$ for $\bP\in\Pcss$, and $c=-1$.
\end{proposition}

\begin{proof}
The proof of \zcref{prop:geometric} applies to the atlas of \zcref{lem:integralatlas}, with
$\mathcal C=\bigcup_w\xi_-^w(\Lat_-^w)\subset\cG(\OO)$ and
$\mathcal A_{P_w}=\uu_{P_w}^\perp\cap\dLat$, and \eqref{eq:badexplicit} gives
$\nu_{\dLat}(\Badmax_h)\le \Nss q^{-h}$. It remains to see that $c=-1$. If
$\lf=f_{P,w}(y,\nf)$ lies outside $\Badmax_h$, then
$\abs{\det df_{P,w}(y,\nf)}>q^{-h}$ by \eqref{eq:sardchart}. In the
factorization of $df_{P,w}$ in the proof of \zcref{lem:levi}, the coadjoint
action of $u=\xi_-^w(y)$ and the differential $T_y$ are lattice
automorphisms, so $\abs{\det Q_\nf}=\abs{\det df_{P,w}(y,\nf)}>q^{-h}$. As
$Q_\nf$ is integral, Cramer's rule gives $\norm{Q_\nf^{-1}}<q^h$, hence
$\norm{Q_\nf^{-1}}\le q^{h-1}$, since the operator norm takes values in $q^{\Z}$.
\end{proof}

\section{Regular cells and Jacquet modules}\label{sec:regular}

We prove the counterpart of \zcref{prop:elliptic}: a cell admitting a
quantitatively regular lift to a proper parabolic subgroup $\bP$ occurs only in
representations whose Jacquet module along $\bP$ is nonzero, and it occurs
only in those tempered representations whose \weakjac{} along $\bP$ is
nonzero. In particular, it cannot occur in a discrete series
representation, let alone in a supercuspidal one. Throughout this section, fix $\bP=\bM\bU\in\Pssall$,
not necessarily maximal,
and a cocharacter $\eta:\Gm\to Z(\bM)$ defining $\bP$ and $\bM$.
Write $\bP^-=\bM\bU^-$ for the opposite parabolic defined by $-\eta$.
Recall that a form $\nf\in\uu^\perp$ defines the map
\[
Q_\nf:\uu^-\longrightarrow\uu^*,\qquad
(Q_{\nf}Y)(X)=\nf([X,Y]).
\]
We fix a norm on $\operatorname{Hom}_F(\uu^*,\uu^-)$ as in \zcref{sec:levi}.

\begin{theorem}[Regular cells]\label{thm:regular}
Fix compact sets $\mathcal C\subset G$ and $\mathcal K\subset\uu^\perp$.
There is a constant $C$ with the following property. Let $h\ge1$ be an integer and let
$\lf\in\dLat$ be of the form $\lf=\Ad^*(g)\nf$ with $g\in\mathcal C$,
$\nf\in\mathcal K$ and $\norm{Q_\nf^{-1}}\le q^h$. If $n$ and $r\ge1$ are integers with
\begin{equation}\label{eq:margins}
n-2r\ge C,\qquad 3r-n-2h\ge C,
\end{equation}
then every irreducible representation $\pi$ of $G$ which contains the character
\[
\chi_{n,\lf}(x)=\psi\bigl(\varpi^{-n}\lf(\kappa(x))\bigr),
\qquad x\in K_{n-r},
\]
has the following properties.
\begin{enumerate}[label=\textup{(\arabic*)},ref=\arabic*]
\item\label{item:reg-jacquet} $r_P(\pi)\ne0$; in particular, $\pi$ is not supercuspidal.
\item\label{item:reg-weak} If $\pi$ is tempered, then $r_P^{\wlab}(\pi)\ne0$; in particular, $\pi\notin\Gds$.
\end{enumerate}
In particular, $\fsc_{n,r}(\lf)=\fds_{n,r}(\lf)=0$.
\end{theorem}

The constant $C$ may depend on the compact sets $\mathcal C,\mathcal K$
and the chosen norm. Once these data are fixed for each
$\bP\in\Pssall$, finiteness of $\Pssall$ gives a common $C$, independent
of $n,r,h,\lf$ and $g$. In the unramified case the margins can be made explicit
(\zcref{prop:regularunram}).

The proof starts with a one-step averaging identity. In product
coordinates we enlarge the cell subgroup and normalize its character
to be trivial on the positive unipotent factor. We then iterate the
identity along an exhaustion of the unipotent radical. The coordinate
estimates that give separation and eventual stationarity are proved in
\zcref{sec:jacquetdetails}.

The averaging argument is related to the injectivity criterion for
Iwahori fixed vectors \cite[Lemma~4.7]{MR0444849} and to the
theory of covers \cite[\S\S6--8]{MR1643417}. Howe used related methods
for sufficiently regular shallow characters of $\GL_d$
\cite[\S3]{MR0492088}.

\subsection{The one-step averaging identity}
For a compact subgroup $C$ let $\Avg_C$ denote averaging with Haar
measure of total mass $1$. If $\theta$ is a character of $C$, let
$\Avg_C^\theta$ be the projection onto $V^{C,\theta}$.
The following elementary observation is the core of the proof. 

\begin{lemma}[One-step averaging]\label{lem:criterion}
Let $(\pi,V)$ be a smooth representation of $G$, $H$ a compact open
subgroup of $G$, $\theta$ a character of
$H$, and $S\subset S'$ compact subgroups with $S\subset H$ and
$\theta|_S=1$, where $[S':S]<\infty$. Suppose that for every
$x\in S'\setminus S$ there is $u\in H$ such that
$xux^{-1}\in H$ and $\theta(xux^{-1})\ne\theta(u)$.
Then, for $E=\Avg_H^\theta$,
\begin{equation}\label{eq:leftinverse}
E\Avg_{S'}v=[S':S]^{-1}v\qquad(v\in V^{H,\theta}).
\end{equation}
In particular, $\Avg_{S'}$ is injective on $V^{H,\theta}$.
\end{lemma}

\begin{proof}
For $x\in S'\setminus S$, choose $u$ as in the statement. Then
\[
\theta(u)E\pi(x)v=E\pi(x)\pi(u)v=E\pi(xux^{-1})\pi(x)v=\theta(xux^{-1})E\pi(x)v,
\]
so $E\pi(x)v=0$. On the other hand, since $v$ is $S$-invariant,
\[
\Avg_{S'}v=[S':S]^{-1}\sum_{x\in S'/S}\pi(x)v.
\]
Applying $E$ kills every coset except $S$, which contributes $v$.
\end{proof}

\subsection{The statement in the product chart}\label{sec:localsetting}

Retain the normalized
parametrization $\xi_-$ from \zcref{sec:incidencejacobian}, and choose
normalized parametrizations $\xi_+$ of $U$ and $\xi_M$ of $M$, with
$\xi_+$ an isomorphism of $F$-varieties. The \emph{product chart} $\kappa'$
is the \nchart{} defined near $1$ on the big cell $UMU^-$ by
\begin{equation}\label{eq:productchart}
\kappa'\bigl(\xi_+(X)\xi_M(X_0)\xi_-(Y)\bigr)=X+X_0+Y,
\end{equation}
for $X\in\uu$, $Y\in\uu^-$ and $X_0$ in the domain of $\xi_M$.

Call a lattice in $\g$ \emph{\etacompat{}} if it is the direct sum
of its intersections with the $\eta$-weight spaces. For such a lattice $\Lat'$ write
$\Lat'=\Lat_+\oplus\Lat_M\oplus\Lat_-$, where
$\Lat_\pm=\Lat'\cap\uu^\pm$ and $\Lat_M=\Lat'\cap\m$. Put $t=\eta(\varpi)$.

\begin{definition}[\Adapted{} pair]\label{def:adapted}
Let $\Lat'$ be an \etacompat{} lattice. The pair $(\kappa',\Lat')$ is
\emph{\adapted{}} if it is standard (\zcref{sec:standardgroups}) and
the conjugations by $t$ on $U$ and by $t^{-1}$ on $U^-$, written in the
coordinates $\xi_+$ and $\xi_-$, have integral coefficients with respect to
$\Lat_+$ and $\Lat_-$, respectively.
\end{definition}

\Adapted{} pairs exist: for every \etacompat{} lattice $\Lat'$, the pair
$(\kappa',\pp^k\Lat')$ is \adapted{} for $k$ large (\zcref{lem:adapted}).

Fix an \adapted{} pair $(\kappa',\Lat')$.
For $d\ge1$, write $U_d=\xi_+(\pp^d\Lat_+)$, $M_d=\xi_M(\pp^d\Lat_M)$ and
$U^-_d=\xi_-(\pp^d\Lat_-)$.
In this subsection, in \zcref{sec:groups,sec:jacquet} and in
\zcref{sec:jacquetdetails}, all norms
are induced by $\Lat'$ and its summands and duals; we use the same
notation $\norm{\cdot}$. Write $\val(\cdot)=-\log_q\norm{\cdot}$, with
$\val(0)=+\infty$, and measure group depths in these coordinates.
Let $w_+$ be the largest weight of $\eta$ on $\uu$.

For each integer $h\ge1$, define
\begin{equation}\label{eq:Qbound}
\mathcal A_h=\{\lf\in(\Lat')^\vee\cap\uu^\perp:
Q_\lf\text{ is invertible and }\norm{Q_\lf^{-1}}\le q^h\}.
\end{equation}
Fix $h\ge1$ and $\lf\in\mathcal A_h$, and put $Q=Q_\lf$.
The same inverse bound holds for $Q^*$. In particular, for every
$X\in\uu$ there is $Y\in\Lat_-$ with
$\abs{\Bf_\lf(X,Y)}\ge q^{-h}\norm{X}$, and the analogous statement
holds with the two factors interchanged.

Let $b,r\ge1$ be integers: $b$ is the depth of the
subgroup and $r$ the precision of the cell. Put $n=b+r$ and $a=r-h$. We assume
\begin{equation}\label{eq:used}
2a\ge b\ge r\ge2(h+w_+).
\end{equation}
The inequality $b\ge r$ makes the cell character well defined; the others
control the commutators in the enlarged groups used below.
By \zcref{lem:groups},
\[
J=U_bM_bU^-_b=\kappa'^{-1}(\pp^b\Lat')
\]
is a compact open subgroup, and its cell character associated with $\lf$ is given by
\begin{equation}\label{eq:thetaJ}
\theta_J\bigl(\xi_+(X)\xi_M(X_0)\xi_-(Y)\bigr)=\psi\bigl(\varpi^{-n}\lf(X_0+Y)\bigr).
\end{equation}

\begin{proposition}[Jacquet and \weakjac{} nonvanishing]\label{prop:local}
Let $(\kappa',\Lat')$ be an \adapted{} pair, let $h\ge1$ and
$\lf\in\mathcal A_h$, and let $b,r\ge1$ be integers satisfying \eqref{eq:used}.
\begin{enumerate}[label=\textup{(\arabic*)},ref=\arabic*]
\item\label{item:local-jacquet} Every smooth representation $V$ of $G$ containing $\theta_J$ has $V_U\ne0$.
\item\label{item:local-weak} Every irreducible tempered representation $\pi$ of $G$ containing $\theta_J$
has $r_P^{\wlab}(\pi)\ne0$; in particular, no representation in $\Gds$ contains $\theta_J$.
\end{enumerate}
\end{proposition}

The proof occupies \zcref{sec:groups,sec:jacquet}.
We will then deduce \zcref{thm:regular}.

\subsection{Extending and normalizing the character}\label{sec:groups}
We enlarge $U_b$ to the subgroup determined by the annihilator lattice
\begin{equation}\label{eq:Slattice}
\sS=\{X\in\uu:\Bf_\lf(X,\pp^b\Lat_-)
\subset\pp^n\}
=\pp^r(Q^*)^{-1}\Lat_-^\vee,\qquad S=\xi_+(\sS).
\end{equation}
Its size is controlled by
\[
\pp^r\Lat_+\subset\sS\subset\pp^{r-h}\Lat_+.
\]
Set
\begin{equation}\label{eq:JH}
R=U^-_b,\qquad H=SM_bR\supset J.
\end{equation}
The groups $J$ and $R$ are independent of $\lf$, whereas $S$ and $H$
depend on $\lf$ through $Q_\lf$; we suppress this dependence in the notation.
By \zcref{lem:groups}, $J$ is normal in $H$, the quotient
$H/J\simeq\sS/\pp^b\Lat_+$ is finite abelian, and the formula
\eqref{eq:thetaJ} extends $\theta_J$ to a character $\theta_\lf$ of
$H$ that is trivial on $S$.

Consequently $V^{J,\theta_J}$ decomposes into the spaces
$V^{H,\theta}$ as $\theta$ runs through the extensions of $\theta_J$.
Every extension other than $\theta_\lf$ is nontrivial on $S$.
However, any extension can be conjugated to one trivial on $S$: by
\zcref{lem:groups,lem:normalizationbounds}, for every character $\theta$ of
$H$ extending $\theta_J$ there are $y\in U^-$ normalizing $H$ and
$\lf'\in\mathcal A_h$, defining the same groups $S$ and $H$, such that
$\theta^y=\theta_{\lf'}$, where $\theta^y(z)=\theta(yzy^{-1})$ for $z\in H$.

Suppose that a representation $(\pi,V)$ of $G$ contains $\theta_J$. By the
decomposition above there is an extension $\theta$ of $\theta_J$ with
$V^{H,\theta}\ne0$, and for $y$ and $\lf'$ as above,
\[
\pi(y^{-1})V^{H,\theta}=V^{H,\theta^y}=V^{H,\theta_{\lf'}}\ne0.
\]
Thus, to prove \zcref{prop:local} we may replace $\lf$ by $\lf'$ and assume
that $V^{H,\theta_\lf}\ne0$.

\subsection{Averaging along a \cchain{}}\label{sec:jacquet}
To iterate this observation we use the following notion. Recall that $t=\eta(\varpi)$.

\begin{definition}[\Cchain{}]\label{def:cchain}
A \emph{\cchain{}} for $(\bP,t)$ is a sequence $(H_j,\theta_j)_{j\ge0}$,
where $H_j=S_jLR$ is a compact open subgroup of $G$ with compact subgroups
$S_j\subset U$, $L\subset M$ and $R\subset U^-$ such that $L$ normalizes $R$,
and $\theta_j$ is a character of $H_j$, such that for all $j\ge0$
\begin{enumerate}[label=\textup{(C\arabic*)},ref=C\arabic*]
\item\label{item:cchain-contract}
$tS_jt^{-1}\subset S_{j+1}$ and $t^{-1}Rt\subset R$;
\item\label{item:cchain-character}
$\theta_j|_{S_j}=1$ and $\theta_{j+1}(u)=\theta_j(t^{-1}ut)$ for $u\in LR$;
\item\label{item:cchain-separate}
for every $x\in t^{-1}S_{j+1}t\setminus S_j$ there is $u\in H_j$ with
$xux^{-1}\in H_j$ and $\theta_j(xux^{-1})\ne\theta_j(u)$;
\item\label{item:cchain-exhaust}
$\bigcup_{j\ge0}t^{-j}S_jt^j=U$.
\end{enumerate}
The \cchain{} is \emph{stationary} if there is $j_0$ such that
$(H_j,\theta_j)=(H_{j_0},\theta_{j_0})$ for all $j\ge j_0$; we then call
$(H_{j_0},\theta_{j_0})$ its \emph{limit}.
\end{definition}

Note that $S_j$ is open in $U$, $L$ in $M$ and $R$ in $U^-$, and the product
decomposition of $H_j$ is unique, by the big cell. In particular, for a
constant sequence \zcref{item:cchain-exhaust} holds automatically, since
conjugation by $t$ contracts $U$. Every tail $(H_j,\theta_j)_{j\ge j_0}$ of a
\cchain{} is a \cchain{}.

\begin{lemma}[Averaging along a \cchain{}]\label{lem:cchain}
Let $(H_j,\theta_j)_{j\ge0}$ be a \cchain{} for $(\bP,t)$, let $(\pi,V)$ be a
smooth representation of $G$ and let $v\in V^{H_0,\theta_0}$. Put
\[
v_j=\Avg_{S_j}\pi(t^j)v,\qquad E_j=\Avg_{H_j}^{\theta_j},\qquad C_j=t^{-j}S_jt^j.
\]
Then for all $j\ge0$,
\begin{enumerate}[label=\textup{(\arabic*)},ref=\arabic*]
\item\label{item:cchain-step} $v_j\in V^{H_j,\theta_j}$ and $v_{j+1}=\Avg_{S_{j+1}}\pi(t)v_j$;
\item\label{item:cchain-back} $E_j\pi(t^{-1})v_{j+1}=[t^{-1}S_{j+1}t:S_j]^{-1}v_j$; in particular, $v_j\ne0$ for all $j$ if $v\ne0$;
\item\label{item:cchain-telescope} $\pi(t^{-j})v_j=\Avg_{C_j}v$.
\end{enumerate}
Consequently, the canonical map $V^{H_0,\theta_0}\to V_U$ is injective.
\end{lemma}

\begin{proof}
Since $\theta_0|_{S_0}=1$, we have $v_0=v$. The nesting
$tS_jt^{-1}\subset S_{j+1}$ gives
$v_{j+1}=\Avg_{S_{j+1}}\pi(t)v_j$.
Inductively, $v_j\in V^{H_j,\theta_j}$: since $t\in Z(M)$ and
$t^{-1}Rt\subset R$, for $u\in LR$ we have
$\pi(u)\pi(t)v_j=\pi(t)\pi(t^{-1}ut)v_j=\theta_{j+1}(u)\pi(t)v_j$ by
\zcref{item:cchain-character}, and
\[
E_{j+1}=\Avg_{S_{j+1}}\Avg_{LR}^{\theta_{j+1}}.
\]
This factorization follows from the unique product decomposition
$H_{j+1}=S_{j+1}(LR)$ and the triviality of $\theta_{j+1}$ on
$S_{j+1}$, with all Haar measures normalized to have total mass one.
Hence $v_{j+1}=E_{j+1}\pi(t)v_j\in V^{H_{j+1},\theta_{j+1}}$, which proves
\zcref{item:cchain-step}. By \zcref{item:cchain-separate}, the
hypothesis of \zcref{lem:criterion} holds for $H_j$, $\theta_j$ and
$S_j\subset S_j'=t^{-1}S_{j+1}t$. Since
$\pi(t^{-1})v_{j+1}=\Avg_{S_j'}v_j$, \eqref{eq:leftinverse} gives
\zcref{item:cchain-back}. Finally, \zcref{item:cchain-telescope} is immediate
from the definition of $v_j$.

By \zcref{item:cchain-back}, $v_{j+1}=0$ implies $v_j=0$, so that $v_j\ne0$
for every $j$ if $v\ne0$. On the
other hand, the compact open subgroups $C_j$ of $U$ increase by
\zcref{item:cchain-contract} and exhaust $U$ by \zcref{item:cchain-exhaust},
so the kernel of $V\to V_U$ is the union of
the kernels of $\Avg_{C_j}$. Hence, by \zcref{item:cchain-telescope}, $v$ has
nonzero image in $V_U$ whenever $v\ne0$.
\end{proof}

We now return to the setting of \zcref{sec:groups}. Put
$\lf_j=\Ad^*(t^j)\lf$. For every $j\ge0$, \zcref{lem:chain} gives
$\lf_j\in\mathcal A_h$.
Applying the construction of \zcref{sec:groups} to $\lf_j$, with $b,r,n,h$ fixed,
defines $S_j$, $H_j=S_jM_bR$ and $\theta_j$, with $S_0=S$, $H_0=H$ and
$\theta_0=\theta_\lf$.

\begin{lemma}[The chain is contracting and stationary]\label{lem:stationary}
For every $\lf\in\mathcal A_h$, the sequence $(H_j,\theta_j)_{j\ge0}$ is a
stationary \cchain{} for $(\bP,t)$, with $L=M_b$ and $R=U^-_b$.
In fact, $(H_j,\theta_j)=(H_r,\theta_r)$ for every $j\ge r$.
\end{lemma}

\begin{proof}
By \zcref{lem:chain}, $(H_j,\theta_j)_{j\ge0}$ is a \cchain{}.
It remains to show that $(H_j,\theta_j)$ does not depend on $j$ for $j\ge r$.
Let $\lf_\infty$ be the form which agrees with $\lf$ on $\m$ and vanishes on
$\uu\oplus\uu^-$, and let $Q_\infty=Q_{\lf_\infty}$. In the notation of the
proof of \zcref{lem:chain}, the block of $Q_j$ from $\g_{-s'}$ to $\g_s^*$
vanishes for $s>s'$ and is $\varpi^{j(s'-s)}$ times the corresponding block of
$Q_0$ for $s\le s'$, while $Q_\infty$ consists of the blocks of $Q_0$ with
$s=s'$. Since $Q_0\Lat_-\subset\Lat_+^\vee$ and the lattices split into weight
spaces, $(Q_j-Q_\infty)\Lat_-\subset\pp^j\Lat_+^\vee$. On the other hand,
$\norm{Q_j^{-1}}\le q^h$ by \zcref{lem:chain}. Hence
\[
D_j:=Q_\infty^*(Q_j^*)^{-1}=1+(Q_\infty-Q_j)^*(Q_j^*)^{-1}
\quad\text{satisfies}\quad(D_j-1)\Lat_-^\vee\subset\pp^{j-h}\Lat_-^\vee,
\]
so that $D_j\Lat_-^\vee=\Lat_-^\vee$ for $j\ge h+1$; in particular, $Q_\infty$
is invertible. By
\eqref{eq:Slattice}, $\sS_j=\pp^r(Q_j^*)^{-1}\Lat_-^\vee
=\pp^r(Q_\infty^*)^{-1}\Lat_-^\vee$ for $j\ge h+1$, so $S_j$ and $H_j$ do not
depend on $j$ for $j\ge h+1$. Finally, $\theta_j$ is trivial on $S_j$, and on
$M_bR$ it is given by \eqref{eq:thetaJ} with $\lf_j$ in place of $\lf$. It
agrees with $\theta_0$ on $M_b$, and it is trivial on $R=U^-_b$ as soon as
$\lf_j(\pp^b\Lat_-)\subset\pp^n$. Since $\Ad(t^{-j})$ is multiplication by
$\varpi^{js}$ on $\g_{-s}$, $s\ge1$, and $\lf$ is integral on $\Lat'$, we have
$\lf_j(\pp^b\Lat_-)\subset\pp^{b+j}$, which is contained in $\pp^n$ for
$j\ge r$. As $r\ge h+1$ by \eqref{eq:used}, the lemma follows.
\end{proof}

Recall that $t=\eta(\varpi)\in A_M$ and that $\abs{\alpha(t)}<1$ for every
root $\alpha$ of $A_M$ in $\uu$. Thus
$\delta_P(t)=\abs{\det(\Ad(t)|_\uu)}<1$, and conjugation by $t$ multiplies a
Haar measure of $U$ by $\delta_P(t)$. Let $(\pi,V)\in\Gtemp$. By the temperedness
criterion \eqref{eq:casselman}, every exponent $\chi$ of $V$ along $P$
satisfies $\abs{\chi(t)}\le1$, with equality if and only if $\chi$ is unitary.
Consequently,
\begin{equation}\label{eq:unitaryeig}
\begin{minipage}{0.85\textwidth}
every finite-dimensional subspace of $r_P(V)$ which is stable under
$r_P(t)$, and on which $r_P(t)$ is diagonalizable with eigenvalues of absolute
value one, is contained in $r_P^{\wlab}(V)$.
\end{minipage}
\end{equation}

\begin{proposition}[Stationary \cchain{}s and the \weakjac{}]\label{prop:weak}
Let $(H_j,\theta_j)_{j\ge0}$ be a stationary \cchain{} for $(\bP,t)$ with
limit $(H,\theta)$, where $H=SLR$,
let $(\pi,V)$ be an irreducible unitarizable representation of $G$, and put
$W=V^{H,\theta}$.
\begin{enumerate}[label=\textup{(\arabic*)},ref=\arabic*]
\item\label{item:weak-unitary} The operator
\begin{equation}\label{eq:unitaryoperator}
A=\delta_P(t)^{-1/2}\Avg_S\pi(t)
\end{equation}
preserves $W$ and is unitary on $W$, for any invariant inner product on $V$.
\item\label{item:weak-inj} The canonical map $j_P:V\to V_U$ is injective on $W$, and
$j_P\circ A=r_P(t)\circ j_P$ on $W$.
\item\label{item:weak-weak} If $\pi$ is tempered, then $j_P(W)\subset r_P^{\wlab}(V)$. In
particular, if $W\ne0$ then $r_P^{\wlab}(\pi)\ne0$, and hence $\pi\notin\Gds$.
\end{enumerate}
Moreover, $W\ne0$ if $V^{H_0,\theta_0}\ne0$.
\end{proposition}

\begin{proof}
The last assertion follows from \zcref{item:cchain-step,item:cchain-back} of \zcref{lem:cchain}.
For the rest, passing to a tail, we may assume that
$(H_j,\theta_j)=(H,\theta)$ for all $j$.
Fix an invariant inner product on $V$. The space $W$ is finite-dimensional
since $\pi$ is admissible.
Let $E=\Avg_H^\theta$ be the orthogonal projection onto $W$. By
\zcref{item:cchain-step} of \zcref{lem:cchain}, the operator $A$
preserves $W$. Since $E\Avg_S=E$, the
restriction of $A$ to $W$ is the restriction of the operator
$\delta_P(t)^{-1/2}E\pi(t)E$ on $V$, whose
adjoint is $\delta_P(t)^{-1/2}E\pi(t^{-1})E$, as $E$ is self-adjoint and $\pi(t)$ is unitary.
Hence $A^*=\delta_P(t)^{-1/2}E\pi(t^{-1})|_W$. By
\zcref{item:cchain-back} of \zcref{lem:cchain},
\[
A^*Aw=\delta_P(t)^{-1}E\pi(t^{-1})\Avg_S\pi(t)w=\delta_P(t)^{-1}[t^{-1}St:S]^{-1}w
=w\qquad(w\in W).
\]
Hence $A^*A=1$ on the finite-dimensional space $W$, so that $A$ is unitary.
This proves \zcref{item:weak-unitary}.

The map $j_P$ is injective on $W$ by \zcref{lem:cchain}.
Moreover $j_P\circ\Avg_S=j_P$, so that
\[
r_P(t)j_P(w)=\delta_P(t)^{-1/2}j_P(\pi(t)w)=j_P(Aw)\qquad(w\in W).
\]
This proves \zcref{item:weak-inj}. Thus $j_P(W)$ is a finite-dimensional
$r_P(t)$-stable subspace of $r_P(V)$ on which $r_P(t)$ is conjugate to the
unitary operator $A$, and
the first part of \zcref{item:weak-weak} follows from \eqref{eq:unitaryeig}.
The last assertion follows, since the \weakjac{} of a discrete series
representation along the proper parabolic subgroup $\bP$ vanishes
(\zcref{sec:representations}).
\end{proof}

\begin{proof}[Proof of \zcref{prop:local}]
Let $(\pi,V)$ be a representation of $G$ containing $\theta_J$. By the
reduction at the end of \zcref{sec:groups} we may assume that
$V^{H,\theta_\lf}\ne0$. By \zcref{lem:stationary}, $(H_j,\theta_j)_{j\ge0}$ is
a stationary \cchain{} with $(H_0,\theta_0)=(H,\theta_\lf)$. Hence
$V^{H,\theta_\lf}$ injects into $V_U$ by \zcref{lem:cchain}, which gives
\zcref{item:local-jacquet}, and if $\pi$ is irreducible and tempered, then
$r_P^{\wlab}(\pi)\ne0$ by \zcref{prop:weak}, which gives
\zcref{item:local-weak}.
\end{proof}

\subsection{Return to the original cell}\label{sec:returncell}
\begin{proof}[Proof of \zcref{thm:regular}]
By \zcref{lem:adapted} we may fix an \adapted{} pair $(\kappa',\Lat')$ with
$\mathcal K\subset(\Lat')^\vee$.
We may prove the theorem using the operator norm induced by $\Lat'$ and its
summands and duals. Indeed, equivalence of norms gives a fixed integer
$e\ge0$ such that a bound $q^h$ in the original norm implies a bound
$q^{h+e}$ in the new norm. Replacing $h$ by $h+e$ and increasing $C$ by $2e$
in \eqref{eq:margins} absorbs this change.
Write $\lf=\Ad^*(g)\nf$ as in the hypothesis. Since
$\nf\in\mathcal K\subset(\Lat')^\vee$ and $\nf|_\uu=0$, the bound
$\norm{Q_\nf^{-1}}\le q^h$ gives $\nf\in\mathcal A_h$.
Set
\[
\kappa_g(x)=\Ad(g)\kappa'(g^{-1}xg),\qquad
\Lat_g=\Ad(g)\Lat'.
\]
Compactness gives a fixed $s\ge1$ with
$\pp^s\Lat_g\subset\Lat$ for every $g\in\mathcal C$.
Apply \zcref{lem:compare} to this family of transported pairs, with
$\mathcal B=\dLat$. Its constants are uniform in $g$.
Put $b=n-r$, $b'=b+s$ and $r'=r-s$, so that $b'+r'=n$. Then
\[
b'-r'=n-2r+2s,\qquad
2(r'-h)-b'=3r-n-2h-3s,\qquad
r'-2h=r-2h-s.
\]
Choose $C$ large enough for the chart comparison, and such that $C\ge3s$ and
$C\ge w_++s$. Since \eqref{eq:margins} gives $r-2h\ge2C$, the integers $b',r'$
satisfy \eqref{eq:used}, with the same $h$. Thus \zcref{prop:local} applies
at depth $b'$ and precision $r'$ to the character $\theta_{J'}$ of
$J'=U_{b'}M_{b'}U^-_{b'}$ attached to $\nf$ by \eqref{eq:thetaJ}.

After conjugation by $g$, this is the character
$x\mapsto\theta_{J'}(g^{-1}xg)=\psi(\varpi^{-n}\lf(\kappa_g(x)))$ of
$J_g=gJ'g^{-1}$. By \zcref{lem:compare}, $J_g\subset K_b$ and this
character is the restriction of $\chi_{n,\lf}$. Let $(\pi,V)$ be an
irreducible representation containing $\chi_{n,\lf}$. Then $V$ contains the
restriction of $\chi_{n,\lf}$ to $J_g$, and $\pi(g^{-1})$ maps the
corresponding isotypic subspace into $V^{J',\theta_{J'}}$. Thus $V$ contains
$\theta_{J'}$. Now \zcref{prop:local} gives \zcref{item:reg-jacquet,item:reg-weak}.
\end{proof}

\subsection{The unramified case}\label{sec:regularunram}
The norms are as before \zcref{prop:geometricunram}.

\begin{proposition}[Regular cells in the unramified case]\label{prop:regularunram}
In the unramified case (\zcref{sec:integraldata}), for $\bP\in\Pcss$ the conclusions of
\zcref{thm:regular} hold with $\mathcal C=\cG(\OO)$ and
$\mathcal K=\uu^\perp\cap\dLat$ whenever
\begin{equation}\label{eq:explicitmargins}
n\ge2r,\qquad 3r-n-2h\ge0,\qquad r-2h\ge2m_{\cG},
\end{equation}
with $m_{\cG}$ as in \eqref{eq:mG}.
\end{proposition}

\begin{proof}
We may assume that $\bP=\dot w\bP_I\dot w^{-1}$ with $\bP_I\in\Pc$ standard and
$\dot w\in N_{\cG}(\cS)(\OO)$, and conjugate the data of
\zcref{sec:integraldata} by $\dot w$. Take $\Lat'=\Lat$, and for $\kappa'$ the
product chart \eqref{eq:productchart} built from $\xi_\pm$ and from $\xi_M$,
the inverse of a model chart (\zcref{lem:model}) of the Levi subgroup scheme
$\cM=Z_{\cG}(\eta)$ of \zcref{sec:integraldata}. The chart $\kappa'$ has integral coefficients and, by
\eqref{eq:integralcell}, identifies $K_b$ with $U_bM_bU^-_b$ for $b\ge1$, where
$U^\pm_b=\cU^\pm(\OO)\cap K_b$ and $M_b=\cM(\OO)\cap K_b$. Its transition map
to $\kappa$ is $X+O(X^2)$ with integral coefficients, so \zcref{lem:compare}
holds with $b_0=1$, $s=0$ and $C_1=0$ for $(\kappa',\Lat)$ and its conjugates by
elements of $\cG(\OO)$. Every $\nf\in\uu^\perp\cap\dLat$ is integral on
$\Lat$, and the conjugations by $t$ on $U$ and by $t^{-1}$ on $U^-$ have
integral coefficients (\zcref{sec:integraldata}). Thus the pair
$(\kappa',\Lat)$ is \adapted{}, and no rescaling or change of norm is needed.
Moreover, $w_+\le m_{\cG}$ by the definition \eqref{eq:mG} of $m_{\cG}$,
since $\eta$ is conjugate to $\eta_I$. Hence, with $b=n-r$, the inequalities
\eqref{eq:used} follow from \eqref{eq:explicitmargins}, and
\zcref{prop:local} applies. Finally, $g\in\cG(\OO)$
and $s=0$, so the return to the original cell (\zcref{sec:returncell})
introduces no shift, and the proof concludes as before.
\end{proof}

The restriction to maximal parabolic subgroups in \zcref{prop:regularunram} is
for convenience only. The same proof applies to any $\bP\in\Pssall$, with
the data of \zcref{sec:integraldata} formed with respect to a cocharacter
$\eta$ of $\cT$ defining $\bP$, provided that $m_{\cG}$ is replaced in
\eqref{eq:explicitmargins} by the largest weight of $\eta$ on $\uu_P$.

\section{Cellwise comparison and proof of the Plancherel limit}\label{sec:proofs}

We now combine the incidence estimates with the two Jacquet criteria
and integrate the resulting cellwise comparison.

When $\bG$ is unramified (\zcref{sec:integraldata}), the constants can be made explicit. Let $\Nss$ be the number of maximal semistandard parabolic subgroups of
$\bG$ (\zcref{sec:parabolics}; for $\GL_d$ we have $\Nss=2^d-2$), and let
$m_{\cG}$ be the invariant \eqref{eq:mG} of the root datum of $\bG$ with its
Galois action, the largest weight on the unipotent radicals of the
cocharacters defining the maximal parabolic subgroups.

\begin{theorem}[Explicit uniform bound]\label{thm:uniform}
Suppose that $\bG$ is unramified, let $\cG$ be a reductive model of $\bG$ over
$\OO$, $\Lat=\Lie(\cG)(\OO)$ and $K_n=\ker(\cG(\OO)\to\cG(\OO/\pp^n))$. Then
for every $n\ge4m_{\cG}+5$ and $\ast\in\{\sclab,\dslab\}$,
\begin{equation}\label{eq:explicitbound}
\bigl|\dens{\ast}{K_n}-\elld{\dLat}\bigr|\le \Nss\,q^{\,m_{\cG}-n/4}.
\end{equation}
\end{theorem}

\begin{corollary}[Chevalley groups]\label{cor:chevalley}
Let $\cG$ be a Chevalley group scheme over $\Z$. For
every non-archimedean local field $F$, of any characteristic, with $\bG=\cG_F$,
$K_n=\ker(\cG(\OO)\to\cG(\OO/\pp^n))$ and $\Lat=\Lie(\cG)(\OO)$,
\[
\bigl|\dens{\ast}{K_n}-\elld{\dLat}\bigr|\le \Nss\,q^{\,m_{\cG}-n/4}\qquad(n\ge4m_{\cG}+5,\ \ast\in\{\sclab,\dslab\}),
\]
where $\Nss$ and $m_{\cG}$ depend only on the root datum.
\end{corollary}

For instance, for $\bG=\GL_d$, $d\ge2$, the estimate becomes
\[
\abs{\dens{\ast}{K_n}-\elld[\GL_d]{\Mat_d(\OO)}}\le(2^d-2)q^{1-n/4},\ \ n\ge9.
\]

For a reductive group over a global field with a fixed linear embedding, the
towers defined by the embedding come from reductive models at almost all
places. This gives a bound uniform over all places (\zcref{thm:global} in \zcref{sec:global}).

\subsection{Jacquet profiles and the least Levi class}\label{sec:profiles}

\begin{convention}[Choice of parameters]\label{rem:parameters}
For each $\bP\in\Pssall$, apply \zcref{thm:regular} with
$\mathcal K=\mathcal A_P$ and the compact set $\mathcal C$ supplied
by \zcref{prop:geometric}.
Take the maximum of the constants in \zcref{thm:regular} over the
finite set $\Pssall$. Fix $C\ge1$ larger than this maximum and the
constant in \zcref{prop:elliptic}, with an additional allowance of $2c$,
where we may assume that the constant $c$ of \zcref{prop:geometric} is nonnegative.
In the rest of this section, the constant in \eqref{eq:margins} is this $C$.

If integers $n,r,h$ with $h\ge1$ satisfy \eqref{eq:margins}, then $n\ge2r+C$ and
$r\ge2h+2C$. In particular, each cell in $\dLat/\pp^r\dLat$ is
either contained in $\Bad_h$ or disjoint from it, and likewise for $\Badmax_h$.
For suitable fixed nonnegative integers $C_1,C_2$ (for instance any $C_1\ge C/2$ and $C_2\ge(C+3C_1+2)/2$) the parameters
\begin{equation}\label{eq:rnhn}
r_n=\lfloor n/2\rfloor-C_1,\qquad h_n=\lfloor n/4\rfloor-C_2
\end{equation}
satisfy \eqref{eq:margins} and $h_n\ge1$ for all sufficiently large $n$,
and by \zcref{prop:geometric} $\nu_{\dLat}(\Bad_{h_n})\ll q^{-\eexpall n/4}$ and
$\nu_{\dLat}(\Badmax_{h_n})\ll q^{-\eexp n/4}$.
\end{convention}

Outside the \bnets{}, a cell determines which Jacquet modules of the
representations containing it vanish, and hence their supports.

\begin{proposition}[Jacquet profile and Levi classes of a cell]\label{prop:profile}\label{prop:partition}
Suppose that integers $n,r,h$ with $h\ge1$ satisfy \eqref{eq:margins}. Let
$\lf\in\dLat$ and let $\pi$ be an irreducible
representation of $G$ containing the character $\chi_{n,\lf}$ of $K_{n-r}$.
Let $\bP$ be a parabolic subgroup of $\bG$, and assume that either
$\lf\notin\Bad_h$, or $\lf\notin\Badmax_h$ and $\bP$ is a maximal proper
parabolic subgroup. Then
\begin{equation}\label{eq:profile}
\lf\in\Pim{P}\iff r_P(\pi)\ne0,
\end{equation}
and if $\pi$ is tempered, these conditions are also equivalent to $r_P^{\wlab}(\pi)\ne0$.
If $\lf\notin\Bad_h$, so that these equivalences hold for all $\bP$, then they
imply the equalities
\begin{equation}\label{eq:labelsagree}
\Lsup{\sclab}(\pi)=\Lform{\lf},\qquad\text{and}\qquad
\Lsup{\dslab}(\pi)=\Lform{\lf}\quad\text{if $\pi$ is tempered.}
\end{equation}
\end{proposition}

\begin{proof}
For $\bP=\bG$ all three conditions hold. Let $\bP$ be proper; since the
conditions depend only on the conjugacy class of $\bP$, we may assume that
$\bP\in\Pstd$, and $\bP\in\Pc$ in the second case of the hypothesis. By
\zcref{item:geo-alt} of \zcref{prop:geometric}, the coset $\lf+\pp^{2h}\dLat$, which
contains the cell $\lf+\pp^r\dLat$, is either disjoint from $\Pim{P}$ or
contained in it. In the first case $r_P(\pi)=0$ by \zcref{prop:elliptic}, and
a fortiori $r_P^{\wlab}(\pi)=0$ if $\pi$ is tempered. In the second case
$\lf=\Ad^*(g)\nf$ with $g\in\mathcal C$, $\nf\in\mathcal A_{P'}$ and
$\norm{Q_{\nf,P'}^{-1}}\le q^{h+c}$ for a semistandard conjugate $\bP'$ of
$\bP$. Apply \zcref{thm:regular} to $\bP'$ with $h+c$ in place of $h$; the extra $2c$ in
the choice of $C$ ensures its depth assumptions. Thus
$r_{P'}(\pi)\ne0$, and $r_{P'}^{\wlab}(\pi)\ne0$ if $\pi$ is tempered, and the
same holds for $\bP$.

Finally, if $\lf\notin\Bad_h$ then $\lf\notin\Crv$, and the least Levi classes on the two sides of the all-parabolic
comparison \eqref{eq:profile} agree by \eqref{eq:least} and \zcref{lem:meet}.
\end{proof}

Note that no noncriticality of the
incidence maps is assumed in \zcref{prop:profile}: if $\incm_P$ is everywhere critical then
$\Pim{P}\subset\Crv$, and the \bnets{} still have measure $O(q^{-h})$ by
\zcref{thm:sard}. Note also that all proper parabolic subgroups are treated
at once, with the same parameters $n,r,h$, so that no loss accumulates along
chains of Levi subgroups.

\subsection{Proof of \texorpdfstring{\zcref{thm:main}}{the main theorem}}\label{sec:proofsgeneral}

Recall that $\fmass{\ast}_{[M],n,r}$, defined in \eqref{eq:fM}, is constant on each
cell and records the proportion of its Plancherel mass carried by $\GhatL{\ast}{[M]}$,
and that $\fmass{\ast}_{n,r}=\fmass{\ast}_{[G],n,r}$.
In particular, $0\le\fmass{\ast}_{[M],n,r}\le1$ and
\[
\densL{\ast}{[M]}{K_n}=\int_\dLat \fmass{\ast}_{[M],n,r}\,d\nu_{\dLat}.
\]

We first compare the cellwise masses; the main estimates follow by integration.

\begin{proposition}[Cellwise comparison]\label{prop:density}
Suppose that integers $n,r,h$ with $h\ge1$ satisfy \eqref{eq:margins}. Then on
$\dLat$, for $\ast\in\{\sclab,\dslab\}$,
\begin{equation}\label{eq:cellcomparisonM}
\abs{\fmass{\ast}_{[M],n,r}-\one_{\Estr_{[M]}}}\le\one_{\Bad_h}\qquad([M]\in\Lcl),
\end{equation}
\begin{equation}\label{eq:cellcomparison}
\abs{\fmass{\ast}_{n,r}-\one_{\Ell{\gd}}}\le\one_{\Badmax_h},
\end{equation}
and the tempered representations $\pi$ with
$\Lsup{\sclab}(\pi)\ne\Lsup{\dslab}(\pi)$ carry no Plancherel mass of a cell
which is disjoint from $\Bad_h$. Consequently,
\begin{equation}\label{eq:masscomparison}
-\nu_{\dLat}(\Badmax_h\cap\Ell{\gd})\le\dens{\ast}{K_n}-\elld{\dLat}\le\nu_{\dLat}(\Badmax_h\cap\Par{\gd}).
\end{equation}
\end{proposition}

\begin{proof}
By \zcref{rem:parameters}, each cell $\cell=\lf+\pp^r\dLat$ is either contained
in $\Bad_h$ or disjoint from it, and likewise for $\Badmax_h$. Since the
functions $\fmass{\ast}_{[M],n,r}$ are constant on cells and take values in
$[0,1]$, it suffices to consider the cells which are disjoint from the net in
question. Recall that an irreducible representation $\pi$ contains the
character $\chi_{n,\lf}$ of $K_{n-r}$ if and only if $\pi(e_{n,\cell})\ne0$.

Let $\cell$ be disjoint from $\Bad_h$. By \zcref{item:geo-alt} of
\zcref{prop:geometric}, each of the sets $\Pim{P}$, hence each $\Pim{[M']}$,
either contains $\cell$ or is disjoint from it. Therefore
$\Lform{\lf'}$ is the same for all $\lf'\in\cell$; call it $[M]$, so that
$\cell\subset\Estr_{[M]}$. By \eqref{eq:labelsagree}, every
$\pi\in\Gtemp$ with $\pi(e_{n,\cell})\ne0$ belongs to
$\GhatL{\sclab}{[M]}\cap\GhatL{\dslab}{[M]}$. By \eqref{eq:fM} and \eqref{eq:cellmass}, this means
that $\fmass{\ast}_{[M'],n,r}=\one_{\Estr_{[M']}}$ on $\cell$
for all $[M']\in\Lcl$, and that the tempered representations with
$\Lsup{\sclab}(\pi)\ne\Lsup{\dslab}(\pi)$ carry no mass of $\cell$.

Let now $\cell$ be disjoint from $\Badmax_h$. By \zcref{item:geo-ell} of
\zcref{prop:geometric}, $\cell$ is contained in $\Ell{\gd}$ or in $\Par{\gd}$.
Let $\pi\in\Gtemp$ with $\pi(e_{n,\cell})\ne0$. By
\zcref{prop:profile} for maximal proper parabolic subgroups, in the first case $r_P(\pi)=0$ for every maximal proper
$\bP$, so that $\pi$ is supercuspidal, and in the second case
$r_P^{\wlab}(\pi)\ne0$ for some maximal proper $\bP$, so that $\pi\notin\Gds$.
Hence $\fmass{\ast}_{n,r}=\one_{\Ell{\gd}}$ on $\cell$.

Finally, the functions $\fmass{\ast}_{n,r}-\one_{\Ell{\gd}}$ are nonpositive on
$\Ell{\gd}$ and nonnegative on $\Par{\gd}$, so integration using
\eqref{eq:Pscdecomp} gives \eqref{eq:masscomparison}.
\end{proof}

\begin{proof}[Proof of \zcref{thm:main}]
Take $r=r_n$ and $h=h_n$ as in \eqref{eq:rnhn}. Integrating
\eqref{eq:cellcomparison} and \eqref{eq:cellcomparisonM} gives
\eqref{eq:main} and \eqref{eq:mainlevi}, using the bounds on
$\nu_{\dLat}(\Badmax_h)$ and $\nu_{\dLat}(\Bad_h)$.
The difference of the discrete series and supercuspidal masses is
nonnegative and has the same error bound by \eqref{eq:main}.
Positivity follows from \zcref{thm:endpoint}. For compact center the
discrete series are atoms of $\mu_G$, with
$\fdeg{\pi}=\mu_G(\{\pi\})\,dg$, giving \eqref{eq:formal}.

By \zcref{prop:density}, only cells in $\Bad_h$ can contribute to the
mass where the support classes disagree. Their total mass is
$\nu_{\dLat}(\Bad_h)$ by \eqref{eq:cellmass}, which proves
\eqref{eq:exceptional}. The joint event for $[M]$ differs from either
of the two individual support events for $[M]$ by at most this mass.
\end{proof}

For a fixed compact subset of $\dLat\setminus\Crv$, the \bnets{} are
eventually disjoint from it, so the cellwise comparison is eventually
exact there. Globally, the examples of $\GL_2$, and of $\PGL_2$ in
characteristic $2$, show that the error need not vanish (\zcref{rem:sharpness}).
Integrating \zcref{prop:density} over a union of cells gives a localized form of the estimate.

\begin{corollary}[Localization]\label{cor:localized}
Let $\ast\in\{\sclab,\dslab\}$ and $[M]\in\Lcl$. For every union
$\cell\subset\dLat$ of cosets of $\pp^{r_n}\dLat$ (see \eqref{eq:rnhn}), with
$e_{n,\cell}$ as in \zcref{sec:cells}, we have
\begin{equation}\label{eq:localizedM}
\vol(K_n)\int_{\GhatL{\ast}{[M]}}\tr\pi(e_{n,\cell})\,d\mu_G(\pi)
=\nu_{\dLat}(\cell\cap\Estr_{[M]})+O(q^{-\eexpall n/4})\qquad(n\gg1)
\end{equation}
and
\begin{equation}\label{eq:localized}
\vol(K_n)\int_{\Ghat{\ast}}\tr\pi(e_{n,\cell})\,d\mu_G(\pi)
=\nu_{\dLat}(\cell\cap\Ell{\gd})+O(q^{-\eexp n/4})\qquad(n\gg1),
\end{equation}
where the implied constants do not depend on $n$ and $\cell$. Moreover, if
$\cell$ is a union of cosets of $\pp^r\dLat$ for some $r$ and
$\cell\cap\Crv=\emptyset$ (resp., $\cell\cap\Crvmax=\emptyset$), then the error
term in \eqref{eq:localizedM} (resp., \eqref{eq:localized}) vanishes for all
sufficiently large $n$.
\end{corollary}

\begin{proof}
We have $\vol(K_n)\int_{\GhatL{\ast}{[M]}}\tr\pi(e_{n,\cell})\,d\mu_G=\int_\cell \fmass{\ast}_{[M],n,r_n}\,d\nu_{\dLat}$ by
\eqref{eq:fM}, and \eqref{eq:localizedM} follows from \eqref{eq:cellcomparisonM}
with $r=r_n$ and $h=h_n$ and from \zcref{rem:parameters}. If
$\cell\cap\Crv=\emptyset$, then $\cell\cap\Bad_h=\emptyset$ for $h\gg1$ by
\eqref{eq:bncompact}, applied to each $\incm_P$, so that $\fmass{\ast}_{[M],n,r_n}=\one_{\Estr_{[M]}}$ on $\cell$
for $n\gg1$ by \zcref{prop:density}. The assertions on $\Ghat{\ast}=\GhatL{\ast}{[G]}$
follow in the same way from \eqref{eq:cellcomparison}, with $\Badmax_h$ and $\Crvmax$ in place of
$\Bad_h$ and $\Crv$.
\end{proof}

\begin{remark}[The lower bound needs only elliptic cells]\label{rem:lower}
The first inequality in \eqref{eq:masscomparison}, and hence
the lower half of \eqref{eq:main}, does not require \zcref{thm:regular}.
In particular, if $\elld{\dLat}=1$, then
$0\le1-\scd{K_n}\ll q^{-n/4}$.
\end{remark}

\subsection{The unramified case}\label{sec:uniform}
We now make the constants explicit in the unramified case
(\zcref{sec:integraldata}). In general, the constants in the proof of
\zcref{thm:main} come from three sources: a finite \iatlas{} of the incidence
map $\incm$, the shifts needed to compare the original tower with the auxiliary
product charts and norms, and depth thresholds governed by the weights of
the defining cocharacters. In the unramified case integral coordinates
eliminate the comparison shifts, and the number of charts and the depth
thresholds are bounded in terms of the root datum, by
\zcref{lem:integralatlas}, the unramified case of \zcref{prop:elliptic}, and
\zcref{prop:geometricunram,prop:regularunram}.

\begin{proposition}[Explicit cellwise comparison]\label{prop:explicitcell}
In the unramified case, let $(\Badmax_h)$ be the \bnets{} of
\zcref{lem:integralatlas}, so that $\nu_{\dLat}(\Badmax_h)\le \Nss q^{-h}$. If integers
$n,r,h$ with $h\ge2$ satisfy
\begin{equation}\label{eq:explicitcellmargins}
n\ge2r,\qquad 3r-n-2h+2\ge0,\qquad r-2h+2\ge2m_{\cG},
\end{equation}
then the inequalities \eqref{eq:cellcomparison} hold on $\dLat$ for these nets.
\end{proposition}

\begin{proof}
If $\Pc=\emptyset$, then $\Ell{\gd}=\gd$, every irreducible representation
is supercuspidal and $\fsc_{n,r}=\fds_{n,r}=1$. Otherwise $m_{\cG}\ge1$, so
\eqref{eq:explicitcellmargins} gives $r-2h\ge2m_{\cG}-2\ge0$. Since $r\ge2h$, each cell $\lf+\pp^r\dLat$ is either contained
in $\Badmax_h$ or disjoint from it. Let it be disjoint. By
\zcref{prop:geometricunram}, either $\lf+\pp^{2h}\dLat\subset\Ell{\gd}$, and then
$\erad(\lf)\le2h\le r$ and \zcref{prop:elliptic} with $C=0$ gives
$\fsc_{n,r}(\lf)=\fds_{n,r}(\lf)=1$; or $\lf=\Ad^*(g)\nf$ with $g\in\cG(\OO)$,
$\nf\in\uu_P^\perp\cap\dLat$ and $\norm{Q_\nf^{-1}}\le q^{h-1}$, and then
\zcref{prop:regularunram} with $h-1$ in place of $h$, whose margins
\eqref{eq:explicitmargins} for $h-1\ge1$ are exactly \eqref{eq:explicitcellmargins}, gives $\fsc_{n,r}(\lf)=\fds_{n,r}(\lf)=0$.
Since $0\le \fsc_{n,r}\le \fds_{n,r}\le1$, the claim follows.
\end{proof}

\begin{proof}[Proof of \zcref{thm:uniform}]
If $\Pc=\emptyset$ there is nothing to prove. Otherwise $m_{\cG}\ge1$; take
$r=\lfloor n/2\rfloor$ and $h=\lfloor(n+3)/4\rfloor-m_{\cG}$, which is $\ge2$
for $n\ge4m_{\cG}+5$. Since $r\ge(n-1)/2$ and $h\le(n+3)/4-m_{\cG}$, we have
$n\ge2r$, $3r-n-2h+2\ge2m_{\cG}-1\ge0$ and $r-2h+2\ge2m_{\cG}$. By \zcref{prop:explicitcell}
and \eqref{eq:badexplicit},
\[
\int_\dLat\abs{\fmass{\ast}_{n,r}-\one_{\Ell{\gd}}}\,d\nu_{\dLat}\le\nu_{\dLat}(\Badmax_h)\le \Nss q^{-h}
\le \Nss q^{\,m_{\cG}-n/4},
\]
since $\lfloor(n+3)/4\rfloor\ge n/4$. The bound for $\dens{\ast}{K_n}$ follows as
in the proof of \zcref{thm:main}.
\end{proof}

\begin{proof}[Proof of \zcref{cor:chevalley}]
By base change, a Chevalley group scheme over $\Z$ gives a reductive model of
$\bG_F$ over $\OO$ for every $F$, and $\Nss$ and $m_{\cG}$ depend only on its
root datum. For $\bG=\GL_d$ we have $\Nss=2^d-2$ and $m_{\cG}=1$
(\zcref{sec:integraldata}), so that $4m_{\cG}+5=9$.
\end{proof}

\begin{remark}[Uniform localized bound]\label{rem:uniformlocalized}
In the setting of \zcref{thm:uniform}, let $n\ge4m_{\cG}+5$ and let
$\cell\subset\dLat$ be a union of cosets of $\pp^{\lfloor n/2\rfloor}\dLat$. Then the error term in \eqref{eq:localized} is at
most $\Nss\,q^{\,m_{\cG}-n/4}$. Indeed, as in the proof of
\zcref{cor:localized}, it equals
$\int_\cell(\fmass{\ast}_{n,r}-\one_{\Ell{\gd}})\,d\nu_{\dLat}$ with
$r=\lfloor n/2\rfloor$, which is bounded by the integral in the proof of
\zcref{thm:uniform}.
\end{remark}

\begin{remark}[A uniform bound with exponent $\eexp/4$ for $\GL_d$]
For hyperspecial $\GL_d$, the error is $O_d(q^{\eexp(1-n/4)})$ for $n\ge9$,
uniformly in $F$. In each integral chart the determinant polynomial takes
a unit value on the direction $\diag(0_i,\one_{d-i})$ under the trace pairing,
so the constant in \zcref{lem:sublevel} depends only on $d$.
Apply \eqref{eq:sardmeasure} with the same $h$ as in the proof of \zcref{thm:uniform}.
\end{remark}

\begin{remark}[Sharpness of the exponent]\label{rem:sharpness}
The exponent $\eexp/4$ in \zcref{thm:main} need not be optimal:
for $\GL_2$ and $\Lat=\Mat_2(\OO)$, where $\eexp=2$ and
$\elld{\dLat}=\frac12$ (\zcref{prop:gl2}), $\frac12-\scd{K_n}$ is of exact order $q^{-2n}$ (\zcref{prop:gl2exact}), while for $\PGL_2$
in characteristic $2$, where $\eexp=1$ and $\elld{\dLat}=1$, $1-\scd{K_n}$ is of exact
order $q^{-n/2}$ (\zcref{prop:pgl2}). We do not know the optimal exponent in general.
For the discrete series, the error vanishes identically for $\GL_2$, while
for $\PGL_2$ in characteristic $2$ it is again of exact order $q^{-n/2}$.

The explicit bound \eqref{eq:explicitbound} is weaker still: for $\GL_2$ it
gives $2q^{1-n/4}$. The general explicit bound uses the exponent $1/4$ in place of
$\eexp/4$ because the implied constant in \zcref{prop:geometric} depends on the
polynomials $\det Q_\nf$ through \zcref{lem:sublevel} and can depend on the field. The preceding remark gives a uniform stronger exponent for $\GL_d$.
For $\GL_2$ with the trace pairing, the \bnets{} $\Badmax_h$ of \zcref{lem:integralatlas} are
contained in the set of $\lf\in\Mat_2(\OO)$ with $\val(\disc\lf)\ge2h$, where
$\disc\lf=(a-d)^2+4bc$; for odd residual characteristic this set has measure
$\ll q^{-2h}$. This is consistent with $\eexp=2$ in \zcref{prop:geometric};
the loss in \eqref{eq:explicitbound} comes from using only the bound
$\Nss q^{-h}$. A further loss comes from the margins
\eqref{eq:margins}, which force $h\le n/4+O(1)$; we do not know how far they
can be relaxed.
\end{remark}

\subsection{Uniformity over a global field}\label{sec:global}
We use the setting and notation of the end of \zcref{sec:integraldata}: $\glf$ is a global
field, $\bG\subset\GL_d$ a connected reductive group over $\glf$ with a closed
embedding, and for every finite place $v$ of $\glf$ we have the lattice
$\Lat_v$ and the groups $K_{v,n}$, as well as the finite set $S$ and the
integer $m_{\bG}$. Let $\dens{\ast}{K_{v,n}}$, $\ast\in\{\sclab,\dslab\}$, be the
supercuspidal and discrete series densities of $K_{v,n}$ in $\bG(\glf_v)$.
For $v\notin S$, \zcref{thm:uniform} applies to the lattice $\Lat_v$ and
the tower $(K_{v,n})_{n\ge1}$. For the finitely many other $v$, the constants of
\zcref{thm:main} for $(K_{v,n})$ depend on the embedding and on $v$.

\begin{theorem}[Uniformity over a global field]\label{thm:global}
For every finite place $v\notin S$ and $\ast\in\{\sclab,\dslab\}$,
\begin{equation}\label{eq:globalgood}
\bigl|\dens{\ast}{K_{v,n}}-\elld[\bG_{\glf_v}]{\Lat_v^\vee}\bigr|\le N_{\bG_{\glf_v}}\,q_v^{\,m_{\bG}-n/4}
\qquad(n\ge4m_{\bG}+5).
\end{equation}
Moreover, there are constants $n_0$ and $C$, depending on $\bG$ and the embedding, such that
\begin{equation}\label{eq:globalall}
\bigl|\dens{\ast}{K_{v,n}}-\elld[\bG_{\glf_v}]{\Lat_v^\vee}\bigr|\le C\,q_v^{\,m_{\bG}-n/4}
\end{equation}
for all finite places $v$ and all $n\ge n_0$.
\end{theorem}

\begin{proof}
For $v\notin S$, \zcref{thm:uniform} applies to $\cG_v$
(\zcref{sec:integraldata}), and since $m_{\cG_v}\le m_{\bG}$, this gives \eqref{eq:globalgood}.

For $v\in S$, \zcref{lem:embedding} and \zcref{thm:main} give constants $C_v$, $n_v$ with
\[
\abs{\dens{\ast}{K_{v,n}}-\elld[\bG_{\glf_v}]{\Lat_v^\vee}}\le C_vq_v^{-n/4}
\]
for $n\ge n_v$. Thus, \eqref{eq:globalall} holds with
\[
n_0=\max(\{4m_{\bG}+5\}\cup\{n_v:v\in S\})
\]
and
\[
C=\max(\{N_{\bG_{\bar\glf}}\}\cup\{C_vq_v^{-m_{\bG}}:v\in S\}),
\]
since $N_{\bG_{\glf_v}}\le N_{\bG_{\bar\glf}}$.
If $S$ is empty, these are $4m_{\bG}+5$ and $N_{\bG_{\bar\glf}}$.
\end{proof}

\begin{remark}\label{rem:globalsplit}
Suppose that $\bG$ is split and that the embedding is the base change of a closed
embedding $\cG\subset\GL_{d}$ over $\Z$, where $\cG$ is a Chevalley group scheme. Then
$\cG_v=\cG\otimes\OO_v$ for every $v$, so that $S=\emptyset$; moreover $\Gamma$ is trivial, $m_{\bG}=m_{\cG}$ and $N_{\bG_{\bar\glf}}=\Nss$, so that
\eqref{eq:globalgood} holds for all $v$ and $n\ge4m_{\cG}+5$; this is \zcref{cor:chevalley} for the completions of $\glf$.
\end{remark}

\subsection{The supercuspidal part of the Hecke algebra}\label{sec:squares}
When $Z$ is compact, the cellwise analysis of \zcref{sec:proofsgeneral} also gives a lower bound
for the sum of the squares of the dimensions of the spaces of $K_n$-fixed
vectors in supercuspidal representations, i.e., the dimension of the supercuspidal part
of the Hecke algebra of $K_n$-biinvariant functions (which is spanned by the $K_n$-biinvariant matrix coefficients of $\pi\in\Gsc$).

\begin{theorem}[Supercuspidal part of the Hecke algebra]\label{thm:squares}
Suppose that $Z$ is compact, and let $K$ be a compact open
subgroup of $G$ such that $K_n$ is a normal subgroup of $K$ for all
$n\gg1$, e.g.\ $K=K_1$. Then
\begin{equation}\label{eq:squares}
\sum_{\pi\in\Gsc}(\dim\pi^{K_n})^2\ge\bigl(\elld{\dLat}-O(q^{-\eexp n/4})\bigr)[K:K_n].
\end{equation}
\end{theorem}

For a smooth model as in \zcref{lem:model}, one may take $K=\cG(\OO)$; then $[K:K_n]=\abs{\cG(k)}q^{(n-1)D}$.

\begin{proof}
Put $r=r_n$, $h=h_n$, and take $n$ large enough that
$\Gamma=K/K_n$ and $H=K_{n-r}/K_n$ are defined. The group $H$ is
abelian of order $q^{rD}$, with characters indexed by the cells in
$\dLat/\pp^r\dLat$. Let $\mathcal W$ be the cells contained in
$\Ell{\gd}\setminus\Badmax_h$. By \zcref{prop:profile} for maximal proper parabolic subgroups, every irreducible
representation of $G$ containing one of their characters is
supercuspidal, and
\[
\frac{|\mathcal W|}{|H|}
\ge\elld{\dLat}-\nu_{\dLat}(\Badmax_h).
\]

In $\C[H]$ let $e$ be the sum of the character idempotents for
$\mathcal W$, and call $\rho\in\Irr(\Gamma)$ good if $\rho(e)\ne0$.
Since $\C[\Gamma]$ is free of rank $[\Gamma:H]$ as a right
$\C[H]$-module,
\[
\dim\C[\Gamma]e=|\mathcal W|[\Gamma:H].
\]
This left ideal is contained in the sum of the matrix blocks of
$\C[\Gamma]$ corresponding to the good representations. Hence
\[
\sum_{\rho\text{ good}}(\dim\rho)^2
\ge|\mathcal W|[\Gamma:H]
\ge\bigl(\elld{\dLat}-\nu_{\dLat}(\Badmax_h)\bigr)[K:K_n].
\]

Every good $\rho$, inflated to $K$, occurs in an irreducible
representation $\pi$ of $G$ by \zcref{rem:Ktype}, and this $\pi$ is
supercuspidal because $\rho(e)\ne0$. Assign each good $\rho$ to one
such $\pi$. The distinct representations assigned to a given $\pi$
occur as a direct sum in $\pi^{K_n}$, so
\[
\sum_{\pi\in\Gsc}(\dim\pi^{K_n})^2
\ge\sum_{\rho\text{ good}}(\dim\rho)^2.
\]
The preceding bound and \zcref{rem:parameters} give \eqref{eq:squares}.
\end{proof}

\section{The elliptic density}\label{sec:positivity}

In this section we prove strict positivity of the elliptic density,
completing the proof of \zcref{thm:main}, and characterize when it equals
$1$ (\zcref{thm:endpoint}).
We also establish the noncriticality statements used in \zcref{sec:incidence} for the improved error exponent.
For \ordinary{} groups (\zcref{def:ordinary}), in
particular whenever $\chr F\ne2$, regular torus-invariant forms
on elliptic maximal tori and on maximal tori of Levi subgroups supply
elliptic forms and noncritical points (\zcref{lem:coroots}). In
characteristic $2$ the adjoint odd orthogonal factors, which are not
\ordinary{}, are treated separately.

The elliptic density $\elld{\dLat}$ is unchanged if $\dLat$ is scaled or conjugated
(\zcref{rem:ball}), but it is otherwise a genuine invariant of the
lattice and not of $\bG$ alone; see \zcref{prop:gl2,prop:iwahori} for $\GL_2$.

Suppose that $\bG=\bG_1\times\bG_2$, and correspondingly $\gd=\g_1^*\oplus\g_2^*$. Then
\begin{equation}\label{eq:prodell}
\Ell{\gd}=\Ell{(\g_1^*)}\times\Ell{(\g_2^*)}.
\end{equation}
For $\bP=\bP_1\times\bP_2$ with Levi subgroup $\bM=\bM_1\times\bM_2$, the map $Q_{\lf,P}$ is the direct sum
of $Q_{\lf|_{\g_1},P_1}$ and $Q_{\lf|_{\g_2},P_2}$, so that
\begin{equation}\label{eq:prodomega}
\Omega_M=\Omega_{M_1}\times\Omega_{M_2}
\end{equation}
(with $\Omega_{G_i}=\g_i^*$).

\subsection{Elementary properties}

\begin{remark}[Scaling, conjugation and positivity]\label{rem:ball}
$\Ell{\gd}$ is stable under multiplication by $F^\times$, so
$\elld{a\dLat}=\elld{\dLat}$ for $a\in F^\times$. Likewise, $\Ell{\gd}$ is
$\Ad^*(G)$-stable, so $\elld{\Ad^*(g)\dLat}=\elld{\dLat}$ for $g\in G$.
Moreover, $\elld{\dLat}>0$ if and only if $\Ell{\gd}\ne\emptyset$.
Indeed, scaling takes any elliptic form into $\dLat$, and openness of
$\Ell{\gd}$ (\zcref{sec:incidencemaps}) then gives a nonempty open
subset of $\dLat\cap\Ell{\gd}$.
\end{remark}

\begin{lemma}[Full-density criterion]\label{lem:one}
$\elld{\dLat}=1$ if and only if $\Omega_M=\emptyset$ for every
$\bP\in\Pc$ (with the notation of \zcref{sec:levi}), i.e., if and only if
every incidence map $\incm_P$ with $\bP\in\Pc$ is everywhere critical. In particular this
condition does not depend on $\Lat$.
\end{lemma}

\begin{proof}
If $\Omega_M=\emptyset$ for all $\bP\in\Pc$ then every point of every
$\Inc_P$ is critical by \eqref{eq:levi}, so $\Par{\gd}=\Crvmax$ has measure zero
by \zcref{thm:sard}, and $\elld{\dLat}=1$. Conversely, if
$\lf\in\Omega_M$ then $d\incm_P$ is invertible at $(\bP,\lf)$, so
$\Par{\gd}\supset \incm_P(\Inc_P)$ contains a nonempty open set, which after
scaling (\zcref{rem:ball}) meets $\dLat$ in a set of positive measure.
\end{proof}

\subsection{Elliptic tori}
A maximal torus $\bT$ of $\bG$ (defined over
$F$) is \emph{elliptic} if its maximal split subtorus is the maximal
split central torus of $\bG$. Equivalently, $\bT$ is contained in no
proper parabolic subgroup of $\bG$.
Elliptic maximal tori exist: this is proved in \cite[Theorem~6.21]{MR1278263}
in characteristic zero and in \cite[\S2.4]{MR2214792} (an argument of
Prasad) in general; cf.~\cite[Theorem~5.3]{MR4617942}.

For $\lf\in\gd$ let $\g_\lf=\{X\in\g:\lf([X,Y])=0\ \forall Y\in\g\}$ be its infinitesimal stabilizer. 

\begin{lemma}[Elliptic forms from elliptic tori]\label{lem:toral}
Let $\bT$ be an elliptic maximal torus with Lie algebra $\liet$ and let
$\lf\in\gd$ be $\bT$-invariant with $\g_\lf=\liet$. Then $\lf\in\Ell{\gd}$.
\end{lemma}

\begin{proof}
Suppose on the contrary that $\lf|_{\uu_P}=0$ for a proper parabolic subgroup $\bP$.
Consider the incidence variety $\bInc_P$ of \zcref{sec:incidencemaps},
and work after base change to $\bar F$
at the point $(\bP,\lf)$. By \eqref{eq:dp}, the
image of the differential of $\incm_P:\bInc_P\to\Lie(\bG)^*$ at this point is
$\uu_P^\perp+\ad^*(\uu_P^-)\lf$. Its annihilator in $\g$ is
$\{X\in\uu_P:\lf([X,\uu_P^-])=0\}=\uu_P\cap\g_\lf$, since
$\lf([\uu_P,\Lie(\bP)])\subset\lf(\uu_P)=0$, and
$\uu_P\cap\g_\lf=\uu_P\cap\liet=0$. Thus the differential is bijective, and hence $(\bP,\lf)$ is an
isolated point of the fiber over $\lf$. The connected group
$\bT_{\bar F}$ fixes $\lf$, hence acts on this fiber and fixes its
isolated points. Therefore $\bT$ normalizes $\bP$, so $\bT\subset\bP$
since parabolic subgroups are self-normalizing. This contradicts the ellipticity of $\bT$.
\end{proof}

For a maximal torus $\bT$ and an absolute root $\alpha$ write
$H_\alpha=d\alpha^\vee(1)\in\liet_{\bar F}$.

\begin{definition}[\Ordinary{} groups]\label{def:ordinary}
We call $\bG$ \emph{\ordinary{}} if
\begin{equation}\label{eq:corootcondition}
H_\alpha\ne0\qquad\text{for every absolute root }\alpha.
\end{equation}
The condition is imposed over $\bar F$ and is independent of the
choice of $\bT$.
\end{definition}

This is the condition called coroot-unramified in \cite[Definition~2.4]{Chi26}.
In characteristic $p>0$, this is equivalent to requiring
$\alpha^\vee\notin pX_*(\bT_{\bar F})$ for every $\alpha$.
Geometrically, it says that the regular semisimple locus in $\gd$ is
nonempty. Indeed, over $\bar F$, identify $\liet^*$ with the forms
vanishing on all root spaces. As observed by Goresky--Kottwitz--MacPherson
\cite[\S6.4]{MR2209851} (see \cite[Lemma~2.2(1) and Corollary~2.3]{Chi26}
for a proof), for $\lf\in\liet^*$ one has
\begin{equation}\label{eq:toralstabilizer}
\g_\lf=\liet\oplus
\bigoplus_{\lf(H_\alpha)=0}\g_\alpha.
\end{equation}
Thus $\g_\lf=\liet$ exactly when every $\lf(H_\alpha)$ is nonzero;
such forms exist exactly under \eqref{eq:corootcondition}.
This is the coadjoint counterpart of the usual root criterion for
regular semisimple elements of $\g$ \cite[\S6.2]{MR2209851}; cf.\ Springer--Steinberg \cite[\S II.3.17]{MR0268192}.
For Jordan decompositions in $\gd$ in arbitrary characteristic, see \cite{ST26}; by
\cite[Theorem~1.1(iv)]{ST26}, the closure of the set of semisimple elements of
$\gd$ is all of $\gd$ unless $\chr F=2$ and $\bG_{\bar F}$ has a direct factor
$\SO_{2d+1}$.

Over $F$ this condition has the following consequence.

\begin{lemma}[Regular torus-invariant forms]\label{lem:coroots}
Assume that $\bG$ is \ordinary{}. Then for every maximal torus $\bT$
of $\bG$ there is a nonempty open subset of the space of $\bT$-invariant
$\lf\in\gd$ with $\g_\lf=\liet$. Consequently
$\Ell{\gd}\ne\emptyset$, and $\Omega_M\ne\emptyset$ for every proper
parabolic subgroup $\bP=\bM\bU$ of $\bG$.
\end{lemma}

\begin{proof}
The projection to the zero weight space identifies $\liet^*$ with the
space of $\bT$-invariant elements of $\gd$. Since each $H_\alpha\ne0$,
$\{\lf\in\liet^*(F):\lf(H_\alpha)=0\}$ is a proper $F$-subspace
and the condition $\lf(H_\alpha)\ne0$ for all $\alpha$ defines a
nonempty open subset of $\liet^*(F)$ (a finite union of proper subspaces
does not cover a vector space over an infinite field).
By \eqref{eq:toralstabilizer}, these forms satisfy $\g_\lf=\liet$.
Taking $\bT$ elliptic (whose existence was recalled at the start of
this subsection), \zcref{lem:toral} gives $\lf\in\Ell{\gd}$. Taking
$\bT\subset\bM$, the pairing $\Bf_{\lf|_\m}$ of
\eqref{eq:levipairing} is, over the splitting field, diagonal with
respect to root vectors with entries $\pm\lf(H_\alpha)\ne0$ for the
roots $\alpha$ of $\uu$, hence non-degenerate. Indeed, $\lf$ is
$\bT$-invariant, so it vanishes on the root spaces; hence
$\lf([X_\alpha,X_{-\beta}])=0$ for $\beta\ne\alpha$, while
$[X_\alpha,X_{-\alpha}]=\pm H_\alpha$ for root vectors from a Chevalley
basis, an identity over $\Z$ valid in every characteristic.
\end{proof}

\subsection{Characteristic two}

Chi \cite[Remark~2.5]{Chi26} notes that the only possible geometric
obstructions to being \ordinary{} are odd orthogonal factors
in characteristic $2$. For our application we need the corresponding
direct-product decomposition over $F$, which we establish below.\footnote{\zcref[S]{lem:split} and its proof are valid for a connected reductive group over an arbitrary field.}

\begin{lemma}[\Oddorth{} and \ordinary{} parts]\label{lem:split}
If $\chr F\ne2$ then $\bG$ is \ordinary{}. If
$\chr F=2$ there is a unique decomposition $\bG=\Godd\times\Gord$ over $F$
such that $\Gord$ is \ordinary{} and $(\Godd)_{\bar F}$ is a
product of adjoint groups of type $\mathrm B_d$, $d\ge1$ (where
$\mathrm B_1=\mathrm A_1$); either factor may be trivial. If $\chr F\ne2$
we put $\Godd=1$ and $\Gord=\bG$. We call $\Godd$ and $\Gord$ the
\emph{\oddorth{} part} and the \emph{\ordinary{} part} of $\bG$. In
particular, $\bG$ is \ordinary{} if and only if $\Godd=1$.
\end{lemma}

\begin{proof}
If $\alpha^\vee\in pX_*$ then $p$ divides $\langle\alpha,\alpha^\vee\rangle=2$,
so $p=2$. Let $p=2$. Since $X_*$ is stable under the Weyl group and every root is
conjugate to a simple root, it suffices to treat simple roots: let
$\alpha$ be a simple root (for some ordering)
of an irreducible component $\Phi_i$ of the absolute root system with
$\alpha^\vee\in2X_*$. Then all the entries $\langle\beta,\alpha^\vee\rangle$
of the column of $\alpha^\vee$ in the Cartan matrix are even, and
inspection of the Cartan matrices shows that this happens only for the
short simple root of $\mathrm B_d$, $d\ge1$. (Dually, simply connected
$\mathrm C_d$ is the only reduced irreducible root datum with a root
divisible in the character lattice; see \cite[proof of Theorem~2.3.10 and
Remark~5.3.11]{MR3362817}. The same reduction to adjoint factors of type
$\mathrm B_d$ appears in \cite[\S4]{ST26}.)

We show that this component splits off as an adjoint factor.
Conceptually, the divisibility of $\alpha^\vee$ forces the cocharacter
lattice of the component to be the full coweight lattice, i.e., the
corresponding normal subgroup is adjoint, and an adjoint factor always
splits off, because the projection of $X_*$ to the span of its coroots
lies in the coweight lattice. In detail:
Let $V_i\subset X_*\otimes\Q$ be the span of the coroots of $\Phi_i$, and let $P_i^\vee$ be the
coweight lattice of $\Phi_i$. In the standard
coordinates, $P_i^\vee=\Z^d$ and $\tfrac12\alpha^\vee=\pm e_j$ for some
$j$; since $X_*$ is stable under the Weyl group, $X_*\cap V_i\supset P_i^\vee$,
hence $X_*\cap V_i=P_i^\vee$. If $V_i'$ is the common kernel of the roots
of $\Phi_i$, then $X_*\otimes\Q=V_i\oplus V_i'$, and the projection of
$X_*$ to $V_i$ pairs integrally with the roots of $\Phi_i$, hence lies
in $P_i^\vee$. Therefore $X_*=P_i^\vee\oplus(X_*\cap V_i')$, and
dually the character lattice and the roots decompose accordingly. Thus
the root datum of $\bG_{\bar F}$ is the direct sum of the adjoint root
datum of type $\mathrm B_d$ and a complementary root datum, and
$\bG_{\bar F}=\bH_i\times\bG_i'$ accordingly.

It remains to descend the product of these factors to $F$.
Let $\Godd$ be the product of the normal subgroups $\bH_i$ over all components $\Phi_i$ containing
a coroot divisible by $2$. Since $\bG$ splits over the separable closure $F^{\mathrm{sep}}$ of $F$, the $\bH_i$ are defined over $F^{\mathrm{sep}}$, and this set of components is stable under $\Gal(F^{\mathrm{sep}}/F)$, so $\Godd$ is defined over $F$, and so is
$\Gord=Z_\bG(\Godd)$, which satisfies $\bG=\Godd\times\Gord$ over $\bar F$
because $\Godd$ has trivial center. Hence $\bG=\Godd\times\Gord$ over $F$,
and no coroot of $\Gord$ is divisible by $2$.

For uniqueness, let $\bG=\bH'\times\bG'$ be another such decomposition,
and work over $\bar F$ with a maximal torus $\bT=\bT_{\bH'}\times\bT_{\bG'}$,
so that $X_*=X_*(\bT_{\bH'})\oplus X_*(\bT_{\bG'})$. Each component
$\Phi_i$ is a component of the root system of $\bH'$ or of $\bG'$. Every
component of $\bH'$ has a coroot divisible by $2$ in $X_*(\bT_{\bH'})$,
hence in $X_*$, since the coroots corresponding to the short roots of an adjoint group of type
$\mathrm B_d$ are divisible by $2$ in its cocharacter lattice. Conversely,
if a coroot $\alpha^\vee$ of $\bG'$ satisfies $\alpha^\vee=2\mu$ with
$\mu\in X_*$, then the projection of $\mu$ to $X_*(\bT_{\bG'})$ shows that
$\alpha^\vee\in2X_*(\bT_{\bG'})$, which is impossible since $\bG'$ is
\ordinary{}. Hence the components of $\bH'$ are exactly those
containing a coroot divisible by $2$, and $\bH'=\Godd$, both being
generated by the root groups of these components. Finally
$\bG'=Z_\bG(\bH')=\Gord$, since $\bH'$ has trivial center.
\end{proof}

\begin{lemma}[Criticality for adjoint odd orthogonal factors]\label{lem:oddorthogonal}
Let $\chr F=2$ and let $\bH$ be an $F$-group such that
$\bH_{\bar F}$ is a product of adjoint groups of type $\mathrm B_d$, $d\ge1$.
Then every incidence map $\incm_P$ attached to a proper parabolic
subgroup of $\bH$ is everywhere critical, so that
$\Par{\lieh^*}$ has measure zero, where $\lieh=\Lie(\bH)(F)$.
\end{lemma}

\begin{proof}
Criticality can be checked over $\bar F$, so we may assume that
$\bH$ is a product of split adjoint groups of type $\mathrm B_d$.
For a product, a proper parabolic is proper in at least one factor,
and $\Bf_\lf$ is the direct sum of the pairings for the factors.
Moreover, a maximal proper parabolic subgroup over $F$ need not remain
maximal over $\bar F$. If $\bP\subset\bQ$ with $\bQ$ maximal proper, then
$\uu_Q\subset\uu_P$ gives a map $\Inc_P\to\Inc_Q$ through which
$\incm_P$ factors, so $\incm_P$ is everywhere critical if $\incm_Q$ is. Thus it suffices to treat maximal proper parabolics
in a single split factor. We will exhibit a degenerate weight block
of $\Bf_\lf$.

We first describe the Lie algebra. Consider a split factor
$\SO(V)$, where $V$ is a vector space of dimension $2d+1$ over $\bar F$
with a non-degenerate quadratic form $Q$, whose polar form $\beta$ has
a one-dimensional radical on which $Q\ne0$. Here $\SO(V)$ denotes the
special orthogonal group of $(V,Q)$ in the sense of
\cite[Definition~C.2.10]{MR3362641}; in characteristic $2$ and odd
dimension it is the smooth connected adjoint group of type $\mathrm B_d$,
and $\mathrm O(Q)=\SO(V)\times\mu_2$ (\cite[Theorem~C.2.11]{MR3362641} and the remark after \cite[Definition~C.2.10]{MR3362641}). The map
\[
\bigwedge^2V\to\mathfrak{so}(V),\qquad
u\wedge v\mapsto\bigl(z\mapsto\beta(v,z)u+\beta(u,z)v\bigr)
\]
is an $\SO(V)$-equivariant isomorphism (cf.~\cite[\S\S3--4]{MR2558787}).
Indeed, every endomorphism $X$ in the image satisfies $\beta(z,Xz)=0$
for all $z$ and $\tr X=0$, so lies in $\mathfrak{so}(V)$.
To see injectivity, decompose $V=\rad\beta\oplus W$ with $\beta|_W$
non-degenerate. Then $\bigwedge^2W$ injects into $\End W$ and
$\rad\beta\wedge W$ into $\Hom(W,\rad\beta)$. Finally, both sides
have dimension $d(2d+1)$.

We now compute the relevant weight block. Let $\bP_i$ be the stabilizer of an isotropic
subspace $E$ of dimension $i$, $1\le i\le d$. Write
$V=E\oplus V_0\oplus E'$ with $E'$ isotropic, $\beta(e_a,f_b)=\delta_{ab}$
for bases $(e_a)$ of $E$ and $(f_b)$ of $E'$, and
$V_0=(E\oplus E')^\perp$ of dimension $2(d-i)+1$. For the cocharacter
acting by $t,1,t^{-1}$ on the three summands, the weight spaces of
$\uu^\pm$ are
\[
\begin{aligned}
\uu_1&=E\wedge V_0, & \uu_2&=\bigwedge^2E,\\
\uu_{-1}&=E'\wedge V_0, & \uu_{-2}&=\bigwedge^2E'.
\end{aligned}
\]
Computing commutators of the endomorphisms above gives, for $v,w\in V_0$,
\[
[e_a\wedge v,f_b\wedge w]=\delta_{ab}\,v\wedge w+\beta_0(v,w)\,e_a\wedge f_b,
\qquad\beta_0=\beta|_{V_0}.
\]
Let $\lf\in\m^*$ and put $a_0(v,w)=\lf(v\wedge w)$,
$C_{ab}=\lf(e_a\wedge f_b)$. Then on $\uu_1\times\uu_{-1}$,
\[
\Bf_\lf(e_a\wedge v,f_b\wedge w)=\delta_{ab}a_0(v,w)+C_{ab}\beta_0(v,w).
\]
Choose a basis of $V_0$, and let $A$ and $B$ be the matrices of
$a_0$ and $\beta_0$. The displayed block has matrix
$I_i\otimes A+C\otimes B$. Over $\bar F$, choose an eigenvector
$x\ne0$ of $C$, with eigenvalue $c$. Since $A+cB$ is alternating
of odd size, there is a nonzero vector $v$ with $(A+cB)v=0$.
Hence
\[
(I_i\otimes A+C\otimes B)(x\otimes v)
=x\otimes(A+cB)v=0.
\]
Thus this block is degenerate for every $\lf$. Since $\Bf_\lf$ is
block diagonal with respect to the weights (as in the proof of
\zcref{lem:levi}), $\Omega_M=\emptyset$.
By \eqref{eq:levi}, $\incm_{P_i}$ is everywhere critical.
The measure-zero assertion follows from \zcref{thm:sard}.
\end{proof}

\begin{proposition}[Peculiar groups]\label{prop:peculiar}
The group $\bG$ is \peculiar{} (\zcref{def:peculiar}) if and only if
$\chr F=2$ and $\bG$ has a direct factor over $F$ which is isotropic over $F$
and whose base change to $\bar F$ is a product of adjoint groups of type
$\mathrm B_d$, $d\ge1$; in other words, if and only if its \oddorth{} part
$\Godd$ (\zcref{lem:split}) is isotropic over $F$.
More precisely, for a proper parabolic
subgroup $\bP=\bP'\times\bP''$ of $\bG=\Godd\times\Gord$, with $\bP'$ and $\bP''$ parabolic subgroups of $\Godd$ and $\Gord$, the map
$\incm_P$ is everywhere critical if and only if $\bP'\ne\Godd$. In
particular, if $\bG$ is not \peculiar{}, then no incidence map of a proper
parabolic subgroup is everywhere critical, while if $\bG$ is \peculiar{} and
$F$-simple, then all of them are.
\end{proposition}

\begin{proof}
Write $\bG=\Godd\times\Gord$ as in \zcref{lem:split}, and let
$\bP=\bP'\times\bP''$ be a proper
parabolic subgroup with Levi subgroup
$\bM=\bM'\times\bM''$.
By \eqref{eq:levi}, $\incm_P$ is everywhere critical if and only if
$\Omega_M=\emptyset$, and by \eqref{eq:prodomega},
$\Omega_M=\Omega_{M'}\times\Omega_{M''}$. The second
factor is nonempty: it is $\gord^*$ if $\bP''=\Gord$, and otherwise
this is \zcref{lem:coroots}, since $\Gord$ is \ordinary{}. If $\bP'=\Godd$, the
first factor is $\godd^*$; otherwise it is empty by
\zcref{lem:oddorthogonal} and \eqref{eq:levi}. This proves the second
assertion. It follows that $\bG$ is \peculiar{} if and only if $\Godd$ has a
proper parabolic subgroup, i.e., is isotropic; and if $\bG$ is \peculiar{} and
$F$-simple, then $\bG=\Godd$. Finally, suppose that $\bG=\bH\times\bH'$ over $F$, where $\bH_{\bar F}$ is a
product of adjoint groups of type $\mathrm B_d$, and let
$\bH'=\bH'_{\mathrm{odd}}\times\bH'_{\mathrm{ord}}$ be the decomposition of
\zcref{lem:split} for $\bH'$. Then $\Godd=\bH\times\bH'_{\mathrm{odd}}$ by the
uniqueness in \zcref{lem:split}. Thus a direct factor of $\bG$ as in the first
formulation is a direct factor of $\Godd$; and $\Godd$ itself
is such a factor when it is isotropic.
\end{proof}

\begin{remark}[Coroots and incidence maps]\label{rem:corootincidence}
Being \ordinary{} implies that $\bG$ is not \peculiar{}, but the converse is not true.
For example, let $\chr F=2$. Then $\bH=\PGL_1(D)$ for a quaternion division
algebra $D$ over $F$ is anisotropic but not \ordinary{}; for $\bG=\PGL_1(D)\times\GL_2$ one
has $\dG{}=1$ and $\eexp=2$.
\end{remark}

\subsection{Positivity and full density}

\begin{theorem}[Positivity and full density]\label{thm:endpoint}
\begin{enumerate}[label=\textup{(\arabic*)},ref=\arabic*]
\item\label{item:end-pos} $\Ell{\gd}\ne\emptyset$, hence $\elld{\dLat}>0$ for every lattice
$\dLat\subset\gd$.
\item\label{item:end-factor} Let $\bG=\Godd\times\Gord$ be as in \zcref{lem:split}. Put $\gord=\Lie(\Gord)(F)$, and let
$\dLat_{\mathrm{ord}}$ be the projection of $\dLat$ to $\gord^*$, with probability Haar measure
$\nu_{\dLat_{\mathrm{ord}}}$. Then
\begin{equation}\label{eq:factor}
\elld{\dLat}=\nu_{\dLat_{\mathrm{ord}}}(\dLat_{\mathrm{ord}}\cap\Ell{(\gord^*)}).
\end{equation}
\item\label{item:end-full} The following are equivalent: $\elld{\dLat}=1$ for some (or every) lattice
$\dLat\subset\gd$; $\Par{\gd}$ is a null set; $\Gord$ is anisotropic modulo its center.
In particular, if $\chr F\ne2$, they hold if and only if $\bG$ is anisotropic modulo its center.
\end{enumerate}
\end{theorem}

\begin{proof}
For \zcref{item:end-factor}, use \eqref{eq:prodell} to write
$\Ell{\gd}=\Ell{(\godd^*)}\times\Ell{(\gord^*)}$, where
$\godd=\Lie(\Godd)(F)$. The first factor is conull by
\zcref{lem:oddorthogonal}, and for trivial $\Godd$ it is the entire
zero-dimensional dual. Thus $\one_{\Ell{\gd}}$ agrees almost everywhere
with the pullback of $\one_{\Ell{(\gord^*)}}$. The projection
$\dLat\to\dLat_{\mathrm{ord}}$ pushes $\nu_{\dLat}$ forward to
$\nu_{\dLat_{\mathrm{ord}}}$, which gives \eqref{eq:factor}.

For \zcref{item:end-pos}, \zcref{lem:coroots} gives an elliptic form on
$\gord$. By \zcref{rem:ball}, the elliptic set has positive measure in
$\dLat_{\mathrm{ord}}$. The formula \eqref{eq:factor} now gives
$\elld{\dLat}>0$, so $\Ell{\gd}\ne\emptyset$.

For \zcref{item:end-full}, note that $\elld{\dLat}=1$ if and only if
$\dLat\cap\Par{\gd}$ is a null set. Since $\Par{\gd}$ is stable under
multiplication by $F^\times$ and $\gd=\bigcup_{k\ge0}\varpi^{-k}\dLat$, this holds for one
lattice if and only if $\Par{\gd}$ is a null set, and then it holds for
every lattice. If $\Gord$ has no proper parabolic subgroup then
$\Ell{(\gord^*)}=\gord^*$ and $\elld{\dLat}=1$ by \eqref{eq:factor}.
Otherwise \zcref{lem:coroots} gives $\Omega_{M''}\ne\emptyset$ for a
maximal proper parabolic subgroup of $\Gord$ with Levi subgroup
$\bM''$, so that
$\nu_{\dLat_{\mathrm{ord}}}(\dLat_{\mathrm{ord}}\cap\Ell{(\gord^*)})<1$ by
\zcref{lem:one} applied to $\Gord$, and $\elld{\dLat}<1$ by
\eqref{eq:factor}. Hence these conditions hold if and only if $\Gord$ has no
proper parabolic subgroup, i.e., $\Gord$ is anisotropic modulo its
center. If $\chr F\ne2$ then $\bG=\Gord$.
\end{proof}

\begin{remark}[Levi densities and the odd orthogonal factor]
Almost every form in $\dLat$ has elliptic projection to
$\godd^*$. For these forms the only parabolics on which the form
annihilates the Lie algebra of the unipotent radical are
$\Godd\times\bP_{\mathrm{ord}}$, where $\bP_{\mathrm{ord}}$ has
the corresponding property for the projected form. Hence the
strata whose Levi class is proper on the odd orthogonal factor
have density zero. The other densities are those of $\Gord$
for $\dLat_{\mathrm{ord}}$, since projection pushes probability
Haar measure to probability Haar measure.
\end{remark}

This completes the proof of \zcref{thm:main}.

\appendix

\section{A quantitative Sard theorem}\label{sec:sard}

We prove a quantitative version of Sard's theorem on the image of critical points as well as a quantitative
inverse function theorem for analytic maps between manifolds of the same dimension, in every characteristic.
The proof provides \bnets{} around the critical values.
Finally, we give an elementary polynomial sublevel estimate.

In positive characteristic, Sard's theorem for analytic maps $X\to Y$
was proved in \cite[Theorem~2.5]{MR3803728} under the condition
$p>\dim X-\dim Y+1$, which is automatic in the equal-dimensional case
considered here. An ultrametric inverse function theorem is given
in \cite[Theorem~2.4]{MR3057121}. The elementary estimates below also
apply when the differential is identically singular. Our quantitative
version is in the spirit of Oesterl\'e's bounds \cite{MR656627}.
In the body of the paper these results are applied only to the incidence
maps, whose coordinate expressions in the atlases of
\zcref{sec:bad,sec:integralatlas} are polynomial; but the analytic setting
is the natural one for the statements and proofs.

\subsection{The local image estimate and inverse theorem}
Endow $F^D$ with the maximum norm and $\OO^D$ with probability Haar
measure. Let $f:\OO^D\to\OO^D$ be given by power series with integral
coefficients converging on $\OO^D$: their coefficients tend to zero
as the total degree tends to infinity. Expanding the monomials at
$x\in\OO^D$ gives, in every characteristic,
\begin{equation}\label{eq:taylor}
f(x+z)=f(x)+df_x(z)+R_x(z),\qquad
R_x(z)=\sum_{\abs\beta\ge2}b_\beta(x)z^\beta,
\end{equation}
with $df_x\in M_D(\OO)$ and $b_\beta(x)\in\OO^D$.
The expansions converge uniformly. Thus, for $z,z'\in\OO^D$,
\begin{equation}\label{eq:remainder}
\norm{R_x(z)}\le\norm z^2,\qquad
\norm{R_x(z)-R_x(z')}\le
\max(\norm z,\norm{z'})\norm{z-z'}.
\end{equation}
Indeed, each difference of monomials of degree at least two contains
a factor $z_i-z_i'$ and at least one further coordinate factor.
Also $J(x)=\det df_x$ is a convergent integral power series, hence
$1$-Lipschitz. Consequently, for $h\ge1$,
\begin{equation}\label{eq:Sh}
S_h(f)=\{x\in\OO^D:\abs{\det df_x}\le q^{-h}\}
\end{equation}
is a union of cosets of $\pp^h\OO^D$.

The second assertion of the next lemma is a quantitative form of the
classical multivariable Newton lemma (cf.\ \cite[\S3, Lemma~2]{MR0207700} and
\cite[Theorem~2.2.1 and Corollary~2.2.1]{MR1743467}).

\begin{lemma}[Image estimate and effective inverse theorem]\label{lem:tube}
For $h\ge1$,
\begin{equation}\label{eq:tube}
\vol\bigl(f(S_h(f))+\pp^{2h}\OO^D\bigr)
\le q^{-h}\vol(S_h(f))\le q^{-h}.
\end{equation}
For every $x\notin S_h(f)$,
\begin{equation}\label{eq:inverseball}
f(x+\pp^h\OO^D)\supset f(x)+\pp^{2h}\OO^D.
\end{equation}
\end{lemma}

\begin{proof}
Fix $x\in S_h(f)$. The remainder estimate \eqref{eq:remainder} gives
\[
f(x+\pp^h\OO^D)+\pp^{2h}\OO^D
\subset f(x)+df_x(\pp^h\OO^D)+\pp^{2h}\OO^D.
\]
If the elementary divisors of $df_x$ have valuations
$s_1,\dots,s_D\in\Z_{\ge0}\cup\{\infty\}$, the set on the right has volume
\[
q^{-Dh-\sum_j\min(s_j,h)}\le q^{-(D+1)h},
\]
since $\sum_js_j\ge h$ implies $\sum_j\min(s_j,h)\ge h$.
There are $q^{Dh}\vol(S_h(f))$ cosets of $\pp^h\OO^D$ in $S_h(f)$;
summing proves \eqref{eq:tube}.

Now let $x\notin S_h(f)$. Cramer's rule gives
$\norm{(df_x)^{-1}}\le\abs{\det df_x}^{-1}\le q^{h-1}$.
For $\delta\in\pp^{2h}\OO^D$, the map
\[
z\longmapsto(df_x)^{-1}(\delta-R_x(z))
\]
maps $\pp^h\OO^D$ to itself and, by \eqref{eq:remainder}, has
Lipschitz constant at most $q^{-1}$. Its fixed point satisfies
$f(x+z)=f(x)+\delta$, proving \eqref{eq:inverseball}.
\end{proof}

\begin{remark}[The Frobenius map]\label{rem:frobenius}
Let $D=1$. The Frobenius map $f(x)=x^p$ in characteristic $p$ has identically
zero differential, so $S_h(f)=\OO$ for every $h$. Its image satisfies
\[
\vol(f(\OO)+\pp^{2h})=q^{\lceil2h/p\rceil-2h}\le q^{-h},
\]
because $x^p\equiv y^p\pmod{\pp^{2h}}$ if and only if
$x\equiv y\pmod{\pp^{\lceil2h/p\rceil}}$. Equality holds for $p=2$.
Thus \zcref{lem:tube} applies even when the differential vanishes
identically, and its bound is attained by a polynomial map in
characteristic $2$. Letting $h\to\infty$ also shows that the image
itself is null, in agreement with \zcref{thm:sard}.
\end{remark}

The lemma holds equally for a convergent integral power series map
between lattices in $D$-dimensional vector spaces, with probability
Haar measures and determinants taken in lattice bases. The conditions
and conclusions are independent of the chosen bases.

\subsection{Finite integral atlases and \bnets{}}\label{sec:bnets}
Let $X$ be an analytic manifold of dimension $D$, let $V$ be an
$F$-vector space of dimension $D$, and let $f:X\to V$ be proper and
analytic. Write $\Vf=f(X)$ and $\Crv=f(\Crit f)$; both are closed.
The inverse function theorem gives $\partial\Vf\subset\Crv$.
Fix a lattice $\dLat\subset V$, let $\nu_{\dLat}$ be the probability Haar measure on $\dLat$, and put
$\XL=f^{-1}(\dLat)$, a compact open subset of $X$.

An \emph{\iatlas{}} of $f$ over $\dLat$ is a finite family of analytic
isomorphisms $\psi_i:\Lat_i\to U_i$ onto open subsets of $\XL$
covering $\XL$, where each $\Lat_i$ is a lattice in a
$D$-dimensional vector space and
$f_i=f\circ\psi_i:\Lat_i\to\dLat$ is a convergent power series with
coefficients in $\dLat$. We use probability Haar measure on each
$\Lat_i$. The empty atlas is allowed when $\XL=\emptyset$.

Such atlases exist: around each $x\in\XL$ choose local coordinates
centered at $x$. The constant term of $f$ lies in $\dLat$; restricting
to a sufficiently small lattice makes all the other coefficients
integral and the series convergent on that lattice
(cf.\ the proof of \cite[Proposition~7.4.1]{MR1743467}). Compactness then gives a finite
cover. For polynomial coordinate maps this is just translation and
rescaling of finitely many coefficients.

\begin{theorem}[Quantitative Sard theorem]\label{thm:sard}
Fix an \iatlas{} $(\psi_i:\Lat_i\to U_i)_{1\le i\le N}$, put
$S_{i,h}=S_h(f_i)$, and define the \emph{\bnets{}}
\begin{equation}\label{eq:badf}
\Bad_h=\bigcup_{i=1}^N\bigl(f_i(S_{i,h})+\pp^{2h}\dLat\bigr)
\subset\dLat,\qquad h\ge1.
\end{equation}
These sets decrease, are unions of cosets of $\pp^{2h}\dLat$, and
every such coset outside $\Bad_h$ is either contained in $\Vf$ or
disjoint from $\Vf$. Moreover,
\begin{equation}\label{eq:sardchart}
\abs{\det df_i(y)}>q^{-h}
\quad\text{if }f_i(y)\in\dLat\setminus\Bad_h,
\end{equation}
\begin{equation}\label{eq:sardmeasure}
\nu_{\dLat}(\Bad_h)\le q^{-h}\sum_{i=1}^N\vol(S_{i,h})\le Nq^{-h},
\end{equation}
and $\bigcap_{h\ge1}\Bad_h=\Crv\cap\dLat$. In particular,
$\Crv$ is a null set.
\end{theorem}

\begin{proof}
The nesting and coset assertions follow from the definition, as does
\eqref{eq:sardchart}. If a coset outside $\Bad_h$ meets $\Vf$, write
a point of intersection as $f_i(y)$. Then $y\notin S_{i,h}$, and
\eqref{eq:inverseball} shows that the entire coset lies in
$f_i(\Lat_i)\subset\Vf$. Summing \eqref{eq:tube} over the charts
gives \eqref{eq:sardmeasure}.

Every critical point in $\XL$ lies in a chart and has zero chart
Jacobian, so $\Crv\cap\dLat\subset\Bad_h$ for all $h$.
Conversely, if $\lf\in\bigcap_h\Bad_h$, finiteness of the atlas
gives an index $i$ and arbitrarily large $h$ with
$\lf=f_i(y_h)+z_h$, where $y_h\in S_{i,h}$ and
$z_h\in\pp^{2h}\dLat$. Along a convergent subsequence in the compact
lattice $\Lat_i$, the limit $y$ satisfies
$\det df_i(y)=0$ and $f_i(y)=\lf$. Hence $\lf\in\Crv$.
The measure bound makes $\Crv\cap\dLat$ null. Applying this to
$\pp^{-k}\dLat$, $k\ge0$, shows that $\Crv$ is a null set.
\end{proof}

Two useful consequences are
\begin{equation}\label{eq:bntube}
(\Crv+\pp^r\dLat)\cap\dLat
\subset\Bad_{\lfloor r/2\rfloor}\qquad(r\ge2),
\end{equation}
since the nets contain $\Crv\cap\dLat$ and are unions of the
corresponding cosets, and, for every compact $\cell\subset\dLat$,
\begin{equation}\label{eq:bncompact}
\cell\cap\Crv=\emptyset\implies
\cell\cap\Bad_h=\emptyset\quad\text{for }h\gg1.
\end{equation}
The latter follows by compactness from the nesting and intersection
assertions. The estimate \eqref{eq:sardmeasure} retains the sizes of
the source sublevel sets; together with the next lemma, this yields
the improved rate in \zcref{sec:incidence}.

\subsection{Polynomial sublevel sets}\label{sec:sublevel}

The following estimate is standard.
For $F=\Q_\ell$ it is stated without proof in \cite[\S3.2, (57)]{MR644559}.
See \cite[Theorem~1.4]{GH26} for a slightly weaker bound, valid in every characteristic,
in a much more general setting. For completeness we provide the details.

\begin{lemma}[Polynomial sublevel estimate]\label{lem:sublevel}
Let $f\in F[x_1,\dots,x_m]$ be a nonzero polynomial of degree $d\ge1$.
Then there is a constant $C_f$ such that
\[
\vol\{x\in\OO^m:\abs{f(x)}\le\epsilon\}\le C_f\epsilon^{1/d}\qquad(\epsilon>0).
\]
\end{lemma}

\begin{proof}
Let $f_d$ be the homogeneous part of degree $d$ of $f$ and choose a
primitive vector $v\in\OO^m$ with $c=f_d(v)\ne0$. Write
$\OO^m=\OO v\oplus W$ for a direct summand $W$, with product Haar
measure. For $w\in W$, the polynomial $t\mapsto f(w+tv)$ has degree
$d$ and leading coefficient $c$, so it factors as
$c\prod_{j=1}^d(t-\alpha_j)$ over an algebraic closure, with the
extended absolute value. If $\abs{f(w+tv)}\le\epsilon$, then
$\abs{t-\alpha_j}\le\rho=(\epsilon/\abs c)^{1/d}$ for some $j$.
Each such set of $t\in F$ is empty or a ball of measure at most
$\rho$. Thus the set of $t\in\OO$ under consideration has measure
at most $d\rho$, uniformly in $w$. Fubini's theorem gives
$C_f=d\abs c^{-1/d}$.
\end{proof}

\section{Technical details for the Jacquet argument}\label{sec:jacquetdetails}

This appendix is part of the proof of \zcref{thm:regular,prop:regularunram}, and hence of
\zcref{thm:main}; it supplies the technical verifications used in
\zcref{sec:regular} to prove that a representation containing the cell
character $\theta_J$ has a nonzero Jacquet module $V_U$.
We verify that the enlarged group $H$ and its characters are well defined,
that an extension of $\theta_J$ can be conjugated to one trivial on $S$,
and that the groups and characters attached to the forms $\Ad^*(t^j)\lf$
form a \cchain{} (\zcref{def:cchain}).
The final subsection compares characters in different charts,
completing the passage from the product chart to the original cell.

We retain the notation of \zcref{sec:localsetting}. Except in
\zcref{lem:adapted} and in \zcref{sec:chartchange}, $(\kappa',\Lat')$ is a
fixed \adapted{} pair;
all depths and norms in the group constructions refer to its lattices.
Thus $b$ is the
depth of $J$ and $r$ is the precision of the cell. Put $n=b+r$;
character values are computed using the additive character
$x\mapsto\psi(\varpi^{-n}x)$. Throughout, $h\ge1$ and $\lf\in\mathcal A_h$:
the form $\lf$ is integral on $\Lat'$, vanishes on $\uu$, and satisfies
$\norm{Q_\lf^{-1}}\le q^h$. We put $a=r-h$, so that $S\subset U_a$,
and assume \eqref{eq:used}.

We first show that \adapted{} pairs exist.
For an \adapted{} pair, the orders of the remainders give the
required depth estimates.
Coordinate errors of depth at least $b$ vanish in additive quotients
such as $\sS/\pp^b\Lat_+$, and only the class modulo $\pp^n$ of the
\emph{phase} $\Phi$ of a character value $\psi(\varpi^{-n}\Phi)$ matters.

\subsection{Integral coordinates and depth estimates}\label{sec:estimates}
\begin{lemma}[Existence of \adapted{} pairs]\label{lem:adapted}
Let $\Lat'$ be an \etacompat{} lattice. Then $(\kappa',\pp^k\Lat')$ is
\adapted{} for all sufficiently large $k$. In particular, for every compact
set $\mathcal K\subset\gd$ there is an \adapted{} pair $(\kappa',\Lat')$ with
$\mathcal K\subset(\Lat')^\vee$.
\end{lemma}

\begin{proof}
The lattices $\pp^k\Lat'$ are \etacompat{}, and $(\kappa',\pp^k\Lat')$ is
standard for $k$ large (\zcref{sec:standardgroups}).
In the coordinates $\xi_\pm$, the linear terms of the conjugations by
$t$ on $U$ and by $t^{-1}$ on $U^-$ are $\Ad(t)|_{\uu}$ and
$\Ad(t^{-1})|_{\uu^-}$; they preserve $\Lat_+$ and $\Lat_-$, since these
lattices are \etacompat{} and the weights of $\eta$ on $\uu$ are positive.
Replacing $\Lat'$ by $\pp^k\Lat'$ multiplies each coefficient of
degree $i\ge2$ by $\varpi^{k(i-1)}$; since these coordinate maps
are polynomial, all their coefficients are integral for $k$ large.
Finally, $\mathcal K\subset(\pp^k\Lat')^\vee$ for $k$ large, since
$\mathcal K$ is bounded.
\end{proof}

For an \adapted{} pair, standardness gives bracket closure and integral
coefficients for the group law, the inverse, and the maps obtained by
multiplying, reordering, and conjugating elements of the three factors.

Let $X,Y$ have depths at least $u,v\ge1$, respectively.
A mixed monomial---one containing at least one coordinate from each
input---has depth at least $u+v$. If its total degree is at least
three, it contains at least one additional coordinate, contributing
depth at least $\min(u,v)$. Its depth is therefore at least
$u+v+\min(u,v)$. These bounds also hold for convergent sums of such
monomials with integral coefficients.

We also need bounds for differences. Suppose $X$ has depth at least
$u$, $Y,Y'$ have depth at least $v$, and $Y-Y'$ has depth at least
$k\ge v$. In the difference of a mixed monomial evaluated at
$(X,Y)$ and $(X,Y')$, each term contains a factor from $Y-Y'$.
Its depth is therefore at least $u+k$, or at least
$u+k+\min(u,v)$ if the monomial has total degree at least three.
We use these elementary observations throughout.

Write $\mu_\pm,\mu_M$ for the group laws in $\xi_\pm,\xi_M$, and
$\mu$ for the group law in $\kappa'$. The other maps needed are
\[
\begin{aligned}
\rho(X,Y)&=\kappa'\bigl(\xi_-(Y)\xi_+(X)\bigr),\\
\gamma_-(X,Y)&=\kappa'\bigl(\xi_-(Y)\xi_+(X)\xi_-(Y)^{-1}\bigr),\\
\gamma_+(X,Y)&=\kappa'\bigl(\xi_+(X)\xi_-(Y)\xi_+(X)^{-1}\bigr).
\end{aligned}
\]
In the product chart $\kappa'$, write the group law as
\[
\mu(X,Y)=X+Y+\beta(X,Y)+O(\deg\ge3),
\]
where $\beta$ is the quadratic term.
Then $\beta(X,Y)-\beta(Y,X)=[X,Y]$, and
\begin{equation}\label{eq:serre}
\kappa'\bigl(\kappa'^{-1}(X)\kappa'^{-1}(Y)
\kappa'^{-1}(X)^{-1}\bigr)
=Y+[X,Y]+O(\deg\ge3)
\end{equation}
(\cite[Part~II, Chapter~IV, \S7, and Chapter~V, \S1]{MR2179691}).
In the product chart $\xi_+(X)\xi_-(Y)$ has coordinate $X+Y$.
Consequently
\begin{equation}\label{eq:quadratic}
\begin{aligned}
\rho(X,Y)&=X+Y+[Y,X]+O(\deg\ge3),\\
\gamma_-(X,Y)&=X+[Y,X]+\mathscr R(X,Y),\\
\gamma_+(X,Y)&=Y+[X,Y]+\mathscr R'(X,Y).
\end{aligned}
\end{equation}
The remainders are mixed and of total degree at least three.
Set $\mathscr R_\lf=\lf\circ\mathscr R$ and $\mathscr R'_\lf=\lf\circ\mathscr R'$. We retain the full pairing
$\Bf_\lf$: the form $\lf$ need not vanish on the
$\uu^-$-component of $[Y,X]$.

\begin{lemma}[Product and phase estimates]\label{lem:estimates}
Let $u,v\ge1$, and let $\lf\in\gd$ be integral on $\Lat'$ and
vanish on $\uu$. In
\zcref{item:est-mixed,item:est-cubic,item:est-diff}, let
$X\in\uu$ and $Y\in\uu^-$ have depths at least $u,v$, respectively.
\begin{enumerate}[label=\textup{(\arabic*)},ref=\arabic*]
\item\label{item:est-group} \emph{Group laws.}
Let $\mu_\bullet$ be any of $\mu_\pm,\mu_M,\mu$.
For inputs $Z,Z'$ in its coordinate space of depths at least $u,v$,
respectively, $\mu_\bullet(Z,Z')-Z-Z'$ has depth at least $u+v$.
If $Z$ has depth at least $u$, the coordinate of its group inverse
plus $Z$ has depth at least $2u$.
\item\label{item:est-mixed} \emph{Reordering and conjugation.}
The errors
$\rho(X,Y)-X-Y$, $\gamma_-(X,Y)-X$, and $\gamma_+(X,Y)-Y$
have depth at least $u+v$. More generally, conjugating an element
of depth at least $u$ in one factor by an element of depth at least
$v$ in another changes its coordinates by an element of depth at
least $u+v$.
\item\label{item:est-cubic} \emph{The commutator phase.}
\begin{equation}\label{eq:mixedphase}
\begin{aligned}
\lf(\gamma_-(X,Y))&=-\Bf_\lf(X,Y)+\mathscr R_\lf(X,Y),\\
\val(\mathscr R_\lf(X,Y))&\ge u+v+\min(u,v).
\end{aligned}
\end{equation}
The phase $\lf(\rho(X,Y))-\lf(Y)$ has the same leading term and
error bound; for $\lf(\gamma_+(X,Y))-\lf(Y)$ the leading term is
$\Bf_\lf(X,Y)$, with the same bound.
\item\label{item:est-diff} \emph{Differences of remainders.}
If $Y'\in\uu^-$ satisfies $\val(Y')\ge v$ and
$\val(Y-Y')\ge k\ge v$, then
\begin{equation}\label{eq:difference}
\val\bigl(\mathscr R_\lf(X,Y)-\mathscr R_\lf(X,Y')\bigr)
\ge u+k+\min(u,v).
\end{equation}
\end{enumerate}
\end{lemma}

\begin{proof}
The group-law and conjugation errors are mixed, and the inversion
errors have degree at least two. This gives
\zcref{item:est-group,item:est-mixed}.
Applying $\lf$ to \eqref{eq:quadratic} gives the leading phases in
\zcref{item:est-cubic}, since $\lf|_\uu=0$.
Every remaining monomial has multidegree at least $(2,1)$ or $(1,2)$,
giving $\min(2u+v,u+2v)=u+v+\min(u,v)$.
For differences, expand $Y^\beta-Y'^\beta$ with one factor from
$Y-Y'$; the same argument gives \zcref{item:est-diff}.
\end{proof}

\begin{lemma}[Contraction estimates]\label{lem:contraction}
For $d\ge1$,
\[
tU_dt^{-1}\subset U_{d+1},\qquad
t^{-1}U^-_dt\subset U^-_{d+1},\qquad
\bigcup_{j\ge0}t^{-j}U_dt^j=U.
\]
If $u\ge2w_++1$, then conjugation by $t^{-1}$ maps every element of $U$ of
depth at least $u$ to an element of depth at least $u-w_+$.
\end{lemma}

\begin{proof}
In $\xi_+$ coordinates, conjugation by $t$ is
$\phi(X)=A_+X+\varepsilon(X)$, where
$A_+=\Ad(t)|_\uu$ has weights $\varpi^s$, $1\le s\le w_+$, and
$\varepsilon$ is integral of degree at least two. Thus an input of
depth $d\ge1$ has image of depth at least $d+1$. The argument for
$t^{-1}$ on $U^-$ is identical. The exhaustion follows from the
contracting action of $\eta$ on $U$.

On $\pp^{w_++1}\Lat_+$ the map
$X\mapsto X+A_+^{-1}\varepsilon(X)$ is a bijection: its nonlinear
part has Lipschitz constant at most $q^{-1}$, so the equation for
each prescribed image is solved by contraction. Hence
\[
\phi(\pp^{w_++1}\Lat_+)
=A_+(\pp^{w_++1}\Lat_+)\supset\pp^{2w_++1}\Lat_+.
\]
For such $X$,
$\val(\varepsilon(X))\ge2\val(X)>\val(A_+X)$, so
$\val(\phi(X))=\val(A_+X)\le\val(X)+w_+$.
Injectivity of conjugation now gives the inverse-depth bound for
every element of depth at least $2w_++1$.
\end{proof}

The remaining proofs use the depth inequalities \eqref{eq:used}.
They imply $a=r-h\ge h+2w_+\ge2w_++1$ since $h\ge1$, $a\le b$, and
$r\ge2h+1$ since $w_+\ge1$.
All inputs below therefore have positive depth. In particular,
quadratic errors of depth $2a$ vanish modulo depth $b$, and errors
of depth $2b\ge n$ do not change the character values.

\subsection{The groups and the extension}\label{sec:groupdetails}
Recall $\sS$ from \eqref{eq:Slattice}. We will also need its
opposite counterpart:
\begin{equation}\label{eq:ST}
\begin{gathered}
\sS=\pp^r(Q^*)^{-1}\Lat_-^\vee,\qquad
\sT=\{Y\in\uu^-:\Bf_\lf(\pp^b\Lat_+,Y)
\subset\pp^n\}=\pp^rQ^{-1}\Lat_+^\vee,\\
S=\xi_+(\sS),\qquad S^-=\xi_-(\sT),\qquad R=U^-_b.
\end{gathered}
\end{equation}
Since $Q$ is integral and $\norm{Q^{-1}}\le q^h$, both lattices
lie between the indicated depth lattices:
\[
\begin{aligned}
\pp^b\Lat_+&\subset\pp^r\Lat_+\subset\sS\subset\pp^a\Lat_+,\\
\pp^b\Lat_-&\subset\pp^r\Lat_-\subset\sT\subset\pp^a\Lat_-.
\end{aligned}
\]
Set $J=U_bM_bR$ and $H=SM_bR$. The product factorizations are unique
by the big cell. On $H$ define $\theta_\lf$ by \eqref{eq:thetaJ}.

\begin{lemma}[Groups and character extension]\label{lem:groups}
The sets $S,S^-,R,J,H$ are compact subgroups; $J,H$ are open, and
$J$ is normal in $H$. The coordinate maps give
\[
H/J\simeq\sS/\pp^b\Lat_+,\qquad
S^-/R\simeq\sT/\pp^b\Lat_-.
\]
Moreover, $S^-$ normalizes $H,J,R$, and $\theta_\lf$ is a character
of $H$, trivial on $S$.
\end{lemma}

\begin{proof}
\emph{Subgroups.}
In $S$ and $S^-$, multiplication and inversion agree with addition and
negation modulo depth $2a\ge b$. Since the defining lattices
contain $\pp^b\Lat_\pm$, they are closed under these operations.
The same argument applies to $U_b^\pm$ and $M_b$. In the product
chart, multiplication and inversion preserve $\pp^b\Lat'$,
so $J=\kappa'^{-1}(\pp^b\Lat')$ is a compact open subgroup.
The error in its group law has depth $2b\ge n$; therefore
\[
\theta_J(x)=\psi(\varpi^{-n}\lf(\kappa'(x)))
\]
is a character of $J$, as in \zcref{sec:cells}.

\smallskip\noindent\emph{Character invariance.}
For $s=\xi_+(X)\in S$ and $j\in J$, the mixed estimate shows that
conjugation changes the coordinates of $j$ only at depth
$a+b\ge b$. Applying this also to $s^{-1}$ shows that $S$
normalizes $J$. To prove invariance of $\theta_J$, it suffices to
consider the three factors of $J$.
On $U_b$ both character values are $1$. For $j\in M_b$,
$sjs^{-1}=s'j$ with $s'\in U_b$: indeed $s'\in U$ by the Levi
action, and $s'\in J$ by normalization. Thus the values again agree.
For $j=\xi_-(Y)\in R$, the change of phase is
\[
\Bf_\lf(X,Y)+\mathscr R'_\lf(X,Y).
\]
The first term has depth at least $n$ by the definition of $\sS$.
Since $a\le b$, the remainder has depth at least
$2a+b\ge n$; here $2a\ge b\ge r$. Thus the change of phase lies in $\pp^n$
and does not change the character value.

\smallskip\noindent\emph{Extension.}
It follows that $H=SJ$ is a compact open subgroup and $J$ is normal
in $H$. Since $S\cap J=U_b$, the coordinate map identifies
$H/J$ with $\sS/\pp^b\Lat_+$. The formula
$\theta_\lf(sj)=\theta_J(j)$ is well defined because $\theta_J$
is trivial on $U_b$. It is a character because $S$ fixes $\theta_J$:
\[
\begin{aligned}
\theta_\lf(s_1j_1s_2j_2)
&=\theta_J(s_2^{-1}j_1s_2)\theta_J(j_2)\\
&=\theta_J(j_1)\theta_J(j_2).
\end{aligned}
\]
This agrees with \eqref{eq:thetaJ} and is trivial on $S$.

\smallskip\noindent\emph{Normalization by $S^-$.}
For $y\in S^-$ and $s=\xi_+(X)\in S$, conjugation changes $X$ only
at depth $2a\ge b$ and introduces $\m$- and $\uu^-$-coordinates
of depth at least $b$. Hence $ysy^{-1}\in H$. Similarly,
$yJy^{-1}\subset J$, since the coordinate change has depth at least
$a+b\ge b$. Applying the same argument to $y^{-1}$ proves that
$S^-$ normalizes $H$ and $J$. Also $R=J\cap U^-$, so $S^-\subset U^-$
normalizes $R$. Finally, the group law of $S^-$ is addition modulo
$\pp^b\Lat_-$, giving the asserted description of $S^-/R$.
\end{proof}

\subsection{Normalizing an extension}\label{sec:normalizationdetails}
An extension of $\theta_J$ differs from $\theta_\lf$ by a character
of $S/U_b$. Write its parameter as $\delta$. Conjugation by
$y=\xi_-(Y)\in S^-$ changes the leading parameter from $\delta$ to
$\delta-QY$. Thus non-degeneracy of $Q$ suggests a unique class of
$Y$ for which the character becomes trivial on $S$. The proof below
shows that the nonlinear terms preserve this conclusion, then checks
that the new linear form satisfies the same bounds.

\begin{lemma}[Quantitative normalization]\label{lem:normalizationbounds}
Let $\theta$ extend $\theta_J$ to $H$. There is a unique class
$yR\in S^-/R$ such that $\theta^y|_S=1$, where
$\theta^y(z)=\theta(yzy^{-1})$. For this $y$ there is (a not necessarily unique)
$\lf'\in\mathcal A_h$ such that $\theta^y=\theta_{\lf'}$ and
\begin{equation}\label{eq:Aprime}
\begin{gathered}
\val(\lf'-\lf)\ge a,\qquad
\norm{Q_{\lf'}^{-1}}\le q^h,\\
(Q_{\lf'}^*)^{-1}\Lat_-^\vee=(Q^*)^{-1}\Lat_-^\vee.
\end{gathered}
\end{equation}
In particular $\lf'$ defines the same groups $S,H$.
\end{lemma}

\begin{proof}
\emph{Parameterizing the extensions.}
With respect to the pairing
$(\delta,X)\mapsto\psi(\varpi^{-n}\delta(X))$, the annihilators of
$\pp^b\Lat_+$ and $\sS$ in $\uu^*$ are, respectively,
$\pp^r\Lat_+^\vee$ and $\pp^bQ\Lat_-$. Consequently
\[
\widehat{S/U_b}\simeq
\mathcal Z:=\pp^r\Lat_+^\vee/\pp^bQ\Lat_-.
\]
Every extension is therefore of the form
\[
\theta_\delta\bigl(\xi_+(X)\xi_M(X_0)\xi_-(Y)\bigr)
=\psi\bigl(\varpi^{-n}(\lf(X_0+Y)+\delta(X))\bigr),
\]
where we choose a representative
$\delta\in\pp^r\Lat_+^\vee$ of its class in $\mathcal Z$.

For $y=\xi_-(Y)\in S^-$ and $X\in\sS$, the phase of
$\theta_\delta^y(\xi_+(X))$ is, modulo $\pp^n$,
\begin{equation}\label{eq:normalphase}
\delta(X)-\Bf_\lf(X,Y)+\mathscr R_\lf(X,Y).
\end{equation}
Indeed, let $\pi_+:\g\to\uu$ be the projection along $\m\oplus\uu^-$.
The phase differs from \eqref{eq:normalphase} by
$\delta(\pi_+\gamma_-(X,Y)-X)$. Since $X$ and $Y$ have depth at least $a$,
$\pi_+\gamma_-(X,Y)-X$ has depth at least $2a$ by
\zcref{item:est-mixed} of \zcref{lem:estimates}, and since $\delta$ has depth at
least $r$, this term has depth at least $r+2a\ge r+b=n$ by \eqref{eq:used}.

\smallskip\noindent\emph{Finding the normalizing element.}
First take $X\in\pp^b\Lat_+$, $Y\in\sT$. The two leading
terms in \eqref{eq:normalphase} have depth at least $n$: this follows
from $r+b=n$ for $\delta(X)$ and from the definition of $\sT$ for
$\Bf_\lf(X,Y)$. The remainder has depth at least
$b+2a\ge n$ by \eqref{eq:mixedphase} and \eqref{eq:used}.
Thus $\theta_\delta^y|_S$ is trivial on $U_b$.
It depends only on $yR$, because conjugation by $R\subset H$ fixes
every character of $H$. We obtain a map of finite sets
\[
\Upsilon:\sT/\pp^b\Lat_-\longrightarrow\widehat{S/U_b},
\qquad [Y]\longmapsto\theta_\delta^{\xi_-(Y)}|_S.
\]
The linear map $Q$ identifies its source with $\mathcal Z$, so source
and target have the same cardinality. We prove injectivity by finding
an element of $S$ on which two distinct images differ.

Suppose $Y\not\equiv Y'\pmod{\pp^b\Lat_-}$, and put
$k=\val(Y-Y')<b$. Choose $\zeta\in\Lat_-^\vee$ with
$\val(\zeta(Y-Y'))=k$. For $z\in\OO$ set
\[
X_z=z\varpi^{n-k-1}(Q^*)^{-1}\zeta.
\]
For this choice, the leading difference in \eqref{eq:normalphase} is
\[
-\Bf_\lf(X_z,Y-Y')=-z\varpi^{n-k-1}\zeta(Y-Y').
\]
As $z$ varies, these values run through $\pp^{n-1}$.
Also $X_z\in\sS$, since $n-k-1\ge r$, and
\[
\val(X_z)\ge u:=n-k-h-1\ge a,\qquad k\ge a.
\]
Apply \eqref{eq:difference} with $v=a$. Its lower bound is
\[
u+k+a=n+r-2h-1,
\]
which is at least $n$ by \eqref{eq:used}. Hence the nonlinear
difference does not affect the character. Since $\psi$ is nontrivial
on $\pp^{-1}$, some $z$ distinguishes the two images of
$\Upsilon$. The map is therefore injective, hence bijective. The
trivial character has a unique preimage, proving existence and
uniqueness of $yR$.

\smallskip\noindent\emph{Controlling the new parameter.}
Choose such a $y=\xi_-(Y)$ and put $\theta'=\theta_\delta^y$.
Set $\tilde\lf=\Ad^*(y^{-1})\lf$.
Conjugation by $y$ preserves $M_bR=J\cap P^-$, on which
$\theta_\delta=\theta_\lf$. In coordinates it has linear term
$\Ad(y)$ and an error of depth at least $2b\ge n$. Thus
\[
\theta'(z)=\psi(\varpi^{-n}\tilde\lf(\kappa'(z)))
\qquad(z\in M_bR).
\]
Define $\lf'$ to agree with $\tilde\lf$ on $\m\oplus\uu^-$ and
to vanish on $\uu$. The integral conjugation series gives
$\val(\Ad(y^{\pm1})-1)\ge\val(Y)\ge a$. Hence
$\val(\lf'-\lf)\ge a$. In particular, $\lf'$ is integral
on $\Lat'$ and vanishes on $\uu$.
Write $Q_{\lf'}=Q+\Delta$. Then
\[
\val(\Delta)\ge a,\qquad
\val(Q^{-1}\Delta)\ge r-2h\ge1.
\]
Thus $D=1+Q^{-1}\Delta$ is a lattice automorphism of $\Lat_-$
and $Q_{\lf'}=QD$. It follows that
\[
\begin{aligned}
(Q_{\lf'}^*)^{-1}\Lat_-^\vee
&=(Q^*)^{-1}(D^*)^{-1}\Lat_-^\vee
=(Q^*)^{-1}\Lat_-^\vee,\\
\norm{Q_{\lf'}^{-1}}&=\norm{D^{-1}Q^{-1}}\le q^h.
\end{aligned}
\]
This proves \eqref{eq:Aprime}. The first equality preserves $S$,
and hence $H$. By \zcref{lem:groups}, $\theta_{\lf'}$ is a character
of this same $H$. It agrees with $\theta'$ on $M_bR$, and both are
trivial on $S$, so $\theta'=\theta_{\lf'}$.
\end{proof}

\subsection{Verification of the moving chain}\label{sec:chaindetails}
We verify that the groups and characters constructed at each stage
form a \cchain{} (\zcref{def:cchain}).
Set $t=\eta(\varpi)$, $\lf_j=\Ad^*(t^j)\lf$ and
$Q_j=Q_{\lf_j}$. Define $\sS_j,S_j,H_j=S_jM_bR$ and
$\theta_j=\theta_{\lf_j}$ as before, with $\lf_j$ in place of
$\lf$. Thus $S_0=S$, $H_0=H$ and $\theta_0=\theta_\lf$.
We first show that the inverse bound for $Q_j$ is uniform in $j$.
The rest of the proof verifies the four conditions of \zcref{def:cchain}.

\begin{lemma}[The moving chain]\label{lem:chain}
For every $j\ge0$, the form $\lf_j$ lies in $\mathcal A_h$ and agrees
with $\lf$ on $\m$.
Thus \zcref{lem:groups} applies to $\lf_j$, and
$\pp^r\Lat_+\subset\sS_j\subset\pp^a\Lat_+$.
Moreover, $(H_j,\theta_j)_{j\ge0}$ is a \cchain{} for $(\bP,t)$, with
$L=M_b$ and $R=U^-_b$.
\end{lemma}

\begin{proof}
\emph{Uniform control of the pairing.}
Since $t$ centralizes $M$ and $\Ad(t^{-j})$ contracts $\uu^-$,
the form $\lf_j$ agrees with $\lf$ on $\m$, still vanishes on
$\uu$, and remains integral on $\Lat'$. Put
$A_+=\Ad(t)|_{\uu}$ and $A_-=\Ad(t^{-1})|_{\uu^-}$; on the
weight spaces $\g_s$ and $\g_{-s}$ these are multiplication by
$\varpi^s$, $s>0$. Equivariance of the pairing gives
\begin{equation}\label{eq:Qj}
Q_j=(A_+^*)^{-j}Q_0A_-^j,\qquad
Q_j^{-1}=A_-^{-j}Q_0^{-1}(A_+^*)^j.
\end{equation}
The block of $Q_0$ from $\g_{-s'}$ to $\g_s^*$ vanishes when
$s>s'$, because $[\g_s,\g_{-s'}]\subset\uu$ is killed by $\lf$.
Thus $Q_0$ is block triangular (as in the proof of \zcref{lem:levi}), and
so is its inverse. In \eqref{eq:Qj}, every possibly nonzero block is
multiplied by a nonnegative power of $\varpi$: by $\varpi^{j(s'-s)}$ for
$Q_0$, and by $\varpi^{j(s-s')}$ for the inverse block from
$\g_s^*$ to $\g_{-s'}$. Since the lattices split into weight spaces,
\[
Q_j\Lat_-\subset\Lat_+^\vee,\qquad
\norm{Q_j^{-1}}\le q^h.
\]
This proves the asserted uniformity and the bounds on $\sS_j$.
By \zcref{lem:groups}, applied to $\lf_j$, $H_j=S_jM_bR$ is a compact open
subgroup and $\theta_j$ is a character of $H_j$, trivial on $S_j$.
Also $M_b$ normalizes $R=J\cap U^-$. It remains to verify the four
conditions of \zcref{def:cchain}.

\smallskip\noindent\emph{Condition \zcref{item:cchain-contract}.}
From \eqref{eq:Qj},
\[
(Q_{j+1}^*)^{-1}
=A_+(Q_j^*)^{-1}(A_-^*)^{-1}.
\]
Since $(A_-^*)^{-1}\Lat_-^\vee\supset\pp^{-1}\Lat_-^\vee$,
we get $\Ad(t)\sS_j\subset\pp\sS_{j+1}$.
The error of conjugation has depth $2a\ge r$, so it too lies in
$\sS_{j+1}$. Hence $tS_jt^{-1}\subset S_{j+1}$.
Moreover, $t^{-1}Rt\subset R$ by \zcref{lem:contraction}.

\smallskip\noindent\emph{Condition \zcref{item:cchain-character}.}
For $u=\xi_-(Y)\in R$, the error of conjugation by $t^{-1}$ has
depth $2b\ge n$. It does not affect the character, and therefore
\[
\theta_j(t^{-1}ut)
=\psi(\varpi^{-n}\lf_j(\Ad(t^{-1})Y))
=\theta_{j+1}(u).
\]
The equality on $M_b$ follows because $t$ centralizes $M$ and the
forms agree there. This proves compatibility on $M_bR$.

\smallskip\noindent\emph{Condition \zcref{item:cchain-separate}.}
Take $x=\xi_+(X)\in t^{-1}S_{j+1}t\setminus S_j$.
By \zcref{lem:contraction},
\[
\val(X)\ge a':=a-w_+.
\]
Since $X\notin\sS_j$, there is $Y_0\in\pp^b\Lat_-$ such that
$\Bf_{\lf_j}(X,Y_0)\notin\pp^n$. The ideal
$\varpi^{-n}\Bf_{\lf_j}(X,Y_0)\OO$ therefore contains $\pp^{-1}$.
Since $\psi$ is nontrivial on $\pp^{-1}$, some $z\in\OO$ gives
$\psi(\varpi^{-n}\Bf_{\lf_j}(X,zY_0))\ne1$.
Put $Y=zY_0\in\pp^b\Lat_-$ and $u=\xi_-(Y)$.
The mixed coordinate error in $xux^{-1}$ has depth at least
$a'+b\ge b$. Thus its $\m$- and $\uu^-$-coordinates have depth
at least $b$, and its $\uu$-coordinate lies in
$\pp^r\Lat_+\subset\sS_j$. Hence $xux^{-1}\in H_j$.

It remains to check that the nonlinear remainder lies in $\pp^n$
and hence does not change the character value. Since
$a'\le b$, its depth is at least $2a'+b$. The excess over the
required depth $n$ is
\[
2a'+b-n=r-2h-2w_+\ge0.
\]
This follows from \eqref{eq:used}.
Consequently
\[
\theta_j(xux^{-1})\theta_j(u)^{-1}
=\psi(\varpi^{-n}\Bf_{\lf_j}(X,Y))\ne1,
\]
which proves separation.

\smallskip\noindent\emph{Condition \zcref{item:cchain-exhaust}.}
The inclusion $S_j\supset U_r$ holds for every $j$, while
$\bigcup_jt^{-j}U_rt^j=U$ by \zcref{lem:contraction}.
\end{proof}

\subsection{Changing the chart}\label{sec:chartchange}
The following lemma compares cell characters in two \nchart{}s.
It is used in \zcref{sec:returncell} to pass from the statement in the
product chart to the original cell. Its constants are made explicit
in the unramified case in the proof of \zcref{prop:regularunram}.
\begin{lemma}[Change of \nchart{}]\label{lem:compare}
Let $\kappa_1$ be another \nchart{}, $\Lat_1$ another
lattice in $\g$, and $\mathcal B\subset\gd$ a bounded set. Fix an integer $s\ge0$
with $\pp^s\Lat_1\subset\Lat$. There are constants $b_0,C_1$ such
that for $\lf\in\mathcal B$ and integers $b,n$ with $b\ge b_0$ and $2b-n\ge C_1$, the formulas
\[
\chi_{n,\lf}(g)=\psi(\varpi^{-n}\lf(\kappa(g))),\qquad
\chi'(g)=\psi(\varpi^{-n}\lf(\kappa_1(g)))
\]
define characters of $K_b=\kappa^{-1}(\pp^b\Lat)$ and
$J_b=\kappa_1^{-1}(\pp^{b+s}\Lat_1)$, respectively.
Moreover, $J_b\subset K_b$ and $\chi'=\chi_{n,\lf}|_{J_b}$.
Consequently
\begin{equation}\label{eq:compare}
e_{J_b}^{\chi'}*e_{K_b}^{\chi_{n,\lf}}=e_{K_b}^{\chi_{n,\lf}}=e_{K_b}^{\chi_{n,\lf}}*e_{J_b}^{\chi'}.
\end{equation}
The constants can be taken uniform when $\kappa_1,\Lat_1$ are
transported by conjugation by elements of a fixed compact subset of $G$,
with $s$ chosen so that the lattice inclusion holds throughout the family.
\end{lemma}

\begin{proof}
Choose $b_0$ so that $\pp^{b_0}\Lat$ and $\pp^{b_0+s}\Lat_1$ lie in the
images of $\kappa$ and $\kappa_1$, respectively. For $b$ sufficiently large,
$J_b$ is a subgroup by \zcref{sec:standardgroups}. Since both charts are normalized, their
transition map has the form
\[
\kappa\circ\kappa_1^{-1}(X)=X+R(X),\qquad R(X)=O(X^2).
\]
Thus, for $X\in\pp^{b+s}\Lat_1\subset\pp^b\Lat$,
\[
R(X)\in\pp^{2b-C_0}\Lat
\]
with a fixed constant $C_0$. For large $b$ this lies in $\pp^b\Lat$,
so $J_b\subset K_b$.

Boundedness of $\mathcal B$ gives an integer $c\ge0$ such that
$\lf(\Lat)\subset\pp^{-c}$ for every $\lf\in\mathcal B$. Hence
\[
\val(\lf(R(X)))\ge2b-C_0-c.
\]
If $2b-n\ge C_0+c$, the transition error does not change the
character value, and $\chi'=\chi_{n,\lf}|_{J_b}$.
The standard group-law estimate \eqref{eq:law} also shows that
$\chi_{n,\lf}$ is a character when $2b-n\ge c$.
Thus $\chi'$ is a character by restriction, and the convolution
identities \eqref{eq:compare} follow from the subgroup inclusion
and agreement of the characters.

For uniformity, replace $\kappa_1$ and $\Lat_1$ by
$x\mapsto\Ad(g)\kappa_1(g^{-1}xg)$ and $\Ad(g)\Lat_1$, with $g$
in a fixed compact set. Compactness permits a common choice of $s$.
The coordinate maps depend analytically on $g$; a finite cover of
the compact set gives a common neighborhood and uniform bounds
for the group laws and transition remainders. The same proof
therefore gives uniform $b_0,C_1$.
\end{proof}

\section{Examples}\label{sec:examples}

We explicate the elliptic set in a few cases and compute the elliptic
densities for three lattices in $\Mat_2(F)$.
We state a conjecture on the Levi densities for $\GL_d$ (\zcref{sec:gldconj}).
We also calculate the Plancherel masses for $\GL_2$ and for $\PGL_2$
in characteristic $2$, obtaining exact formulas that also measure the
sharpness in the general error estimate.

\subsection{The elliptic set}

The following elementary lemma computes the elliptic densities of
$\GL_2$ for three lattices, in every characteristic.
\begin{lemma}[Random monic quadratics]\label{lem:quadratics}
Let $\rho_{\OO}$, resp.\ $\rho_{\pp}$, be the probability that
$T^2+bT+c$ is irreducible over $F$ when $(b,c)$ is uniformly distributed
in $\OO^2$, resp.\ in $\pp\times\pp$. Then
\[
\rho_{\OO}=\frac{q+2}{2(q+1)},\qquad\rho_{\pp}=\frac{2q+1}{2(q+1)}.
\]
Moreover, let $x\in\pp$, $y\in\OO^\times$, and
$X=\begin{psmallmatrix}x&y\\z&w\end{psmallmatrix}$. If $(z,w)$ is uniformly
distributed in $\pp\times\pp$, then the characteristic polynomial of $X$ is
irreducible with probability $\rho_{\pp}$.
\end{lemma}

\begin{proof}
A monic quadratic polynomial over $F$ is reducible if and only if it has
a root in $F$, and a root of a monic polynomial with coefficients in $\OO$
lies in $\OO$.

Let $(b,c)$ be uniform in $\pp\times\pp$. With probability $(q-1)/q$ we
have $\val(c)=1$, and $T^2+bT+c$ is an Eisenstein polynomial, hence
irreducible. With probability $1/q$ we have $c\in\pp^2$; then
$T^2+bT+c=\varpi^2\bigl(S^2+b'S+c'\bigr)$ with $S=T/\varpi$, $b'=b/\varpi$,
$c'=c/\varpi^2$, and conditionally $(b',c')$ is uniform in $\OO^2$. Hence
$\rho_{\pp}=\frac{q-1}q+\frac1q\rho_{\OO}$.

Let $(b,c)$ be uniform in $\OO^2$, and reduce $f=T^2+bT+c$ modulo $\pp$:
$\bar f$ is uniformly distributed among the $q^2$ monic quadratic polynomials
over $k$. If $\bar f$ is irreducible, which happens for $(q^2-q)/2$ of them,
then $f$ is irreducible. If $\bar f$ has two distinct roots in $k$, which
happens for $(q^2-q)/2$ of them, then $f$ has a root in $\OO$ by Hensel's
lemma. In the remaining $q$ cases $\bar f=(T-\bar a)^2$ for some
$\bar a\in k$. Fix a lift $a\in\OO$ of $\bar a$. Then
$f(T+a)=T^2+(b+2a)T+(c+ab+a^2)$, and $(b,c)\mapsto(b+2a,c+ab+a^2)$ is a
measure preserving bijection from the residue class of $(b,c)$ onto
$\pp\times\pp$: it is affine-triangular with unit diagonal coefficients
in $b$ and $c$, and its image consists of the coefficients of a lift of
$\bar f(T+\bar a)=T^2$. Hence
$\rho_{\OO}=\frac{q-1}{2q}+\frac1q\rho_{\pp}$. The two linear relations
give the stated values.

For the last assertion, the characteristic polynomial of $X$
is $T^2-(x+w)T+(xw-yz)$, and
$(z,w)\mapsto(-(x+w),xw-yz)$ is a measure preserving bijection of
$\pp\times\pp$ onto itself, since $w\mapsto-(x+w)$ is one of $\pp$ and, for
fixed $w\in\pp$, $z\mapsto xw-yz$ is one of $\pp$, as $y$ is a unit.
\end{proof}

\begin{proposition}[The hyperspecial $\GL_2$ lattice]\label{prop:gl2}
For every non-archimedean local field $F$, $\elld[\GL_2]{\Mat_2(\OO)}=\frac12$.
\end{proposition}

\begin{proof}
By \zcref{ex:gld}, $\elld[\GL_2]{\Mat_2(\OO)}$ is the probability
$P$ that the characteristic polynomial of a uniformly distributed
$X\in\Mat_2(\OO)$ is irreducible over $F$. Conjugation by $\GL_2(\OO)$
preserves Haar measure on $\Mat_2(\OO)$ and this property, and
$\GL_2(\OO)\to\GL_2(k)$ is surjective. We stratify according to the
conjugacy class of $\bar X\in\Mat_2(k)$ \cite[5.3~Proposition,
6.1]{MR2234120}; conditionally on $\bar X$, the matrix $X$ is uniformly
distributed in a coset of $\pp\Mat_2(\OO)$.
\begin{enumerate}[label=\textup{(\alph*)},ref=\alph*]
\item $\bar X$ has irreducible characteristic polynomial: there are
$(q^2-q)/2$ such polynomials, each the characteristic polynomial of
$q(q-1)$ matrices, so this has probability $\frac{(q-1)^2}{2q^2}$, and
then the characteristic polynomial of $X$ is irreducible.
\item $\bar X$ has two distinct eigenvalues in $k$: probability
$\frac{q(q-1)}2\cdot\frac{q(q+1)}{q^4}=\frac{q^2-1}{2q^2}$, and then the
characteristic polynomial of $X$ has a root in $F$ by Hensel's lemma.
\item $\bar X=\bar a+\bar N$ with $\bar a\in k$ and $\bar N\ne0$
nilpotent: probability $\frac{q(q^2-1)}{q^4}=\frac{q^2-1}{q^3}$. After
conjugation by $\GL_2(\OO)$ we may assume
$\bar X=\begin{psmallmatrix}\bar a&1\\0&\bar a\end{psmallmatrix}$. Fix a
lift $a\in\OO$ of $\bar a$; then
$X-a=\begin{psmallmatrix}\varpi s&1+\varpi t\\\varpi u&\varpi v\end{psmallmatrix}$
with $s,t,u,v$ uniform in $\OO$, and the characteristic polynomial of $X$
is irreducible if and only if that of $X-a$ is. By the last assertion
of \zcref{lem:quadratics} (with $x=\varpi s$, $y=1+\varpi t$), this has
conditional probability $\rho_{\pp}$.
\item $\bar X=\bar a$ is scalar: probability $q^{-3}$. Then
$X=a+\varpi Y$ with $Y$ uniform in $\Mat_2(\OO)$, and the characteristic
polynomial of $X$ is irreducible if and only if that of $Y$ is; the
conditional probability is $P$.
\end{enumerate}
Since $\frac{q^2-1}{q^3}\rho_{\pp}=\frac{(q-1)(2q+1)}{2q^3}$, we obtain
\[
P=\frac{(q-1)^2}{2q^2}+\frac{(q-1)(2q+1)}{2q^3}+q^{-3}P
=\frac{1-q^{-3}}2+q^{-3}P,
\]
whence $P=\frac12$.
\end{proof}

For $F=\Q_p$, including $p=2$, the same value follows from the joint
eigenvalue distribution of Haar random matrices in $\Mat_2(\Z_p)$
computed in \cite{SVP26}: by \cite[Theorem~1.7]{SVP26} with $n=m=2$
and $K_1=K_2=\Q_p$, the density of ordered pairs of distinct eigenvalues
in $\Z_p^2$ is $(1+p^{-1})\abs{x-y}$, so the probability that both eigenvalues lie in
$\Z_p$ is $\frac{1+p^{-1}}2\int_{\Z_p^2}\abs{x-y}\,dx\,dy=\frac12$. The
exact computation of $\scd{K_n}$ in \zcref{prop:gl2exact} below gives
$\lim_n\scd{K_n}=\frac12$ for every $F$, in accordance with
\zcref{thm:main} and \zcref{prop:gl2}.

To see the dependence on the lattice, consider
\begin{equation}\label{eq:iwahorilattices}
\Lat_1=\begin{pmatrix}\OO&\OO\\\pp&\OO\end{pmatrix},\qquad
\Lat_2=\begin{pmatrix}\OO&\pp^{-1}\\\OO&\OO\end{pmatrix}.
\end{equation}
$\Lat_1$ is the Lie algebra of the Iwahori group scheme of $\GL_2$,
whose $\OO$-points form the Iwahori subgroup $\Iw$ of \zcref{sec:gl2exact},
and the tower of $\Lat_1$ consists of its congruence kernels
$\{g\in\Iw:g\equiv1\ (\pp^n),\ g_{21}\equiv0\ (\pp^{n+1})\}$
(\zcref{lem:model}).
By contrast, $\Lat_2$ is not closed under the bracket; its tower is
defined for $n$ large, as the tower of a rescaled standard pair
$(\kappa,\pp^{n_0}\Lat_2)$ (\zcref{sec:standardgroups}).

\begin{proposition}[Iwahori-type lattices]\label{prop:iwahori}
For every non-archimedean local field $F$,
\[
\elld{\Lat_1^\vee}=\frac q{q+1},\qquad \elld{\Lat_2^\vee}=\frac1{q+1}.
\]
\end{proposition}

\begin{proof}
For $i=1,2$ let $\phi_i$ be the probability that the characteristic
polynomial of a uniformly distributed element of $\Lat_i$ is irreducible
over $F$. Under the trace pairing $\Lat_1^\vee=\Lat_2$ and
$\Lat_2^\vee=\Lat_1$, so $\elld{\Lat_1^\vee}=\phi_2$ and
$\elld{\Lat_2^\vee}=\phi_1$ by \zcref{ex:gld}.

We first compute $\phi_1$. Let
$Y=\begin{psmallmatrix}a&b\\c&d\end{psmallmatrix}$ be uniformly distributed in
$\Lat_1$, so that $a,b,d$ are uniform in $\OO$ and $c$ is uniform in $\pp$. The
reduction of $Y$ is upper triangular with diagonal entries $\bar a,\bar d$.
If $\bar a\ne\bar d$, which has probability $(q-1)/q$, the characteristic
polynomial of $Y$ has a root in $F$ by Hensel's lemma. If $\bar a=\bar d$, fix
a lift $a_0\in\OO$; then
$Y-a_0=\begin{psmallmatrix}\varpi s&b\\\varpi c'&\varpi v\end{psmallmatrix}$
with $s,b,c',v$ uniform in $\OO$. If $b\in\OO^\times$, which has
probability $(q-1)/q$, the conditional probability of irreducibility is
$\rho_{\pp}$ by \zcref{lem:quadratics} with $x=\varpi s$ and $y=b$. If
$b\in\pp$, then $Y-a_0=\varpi Y_1$ with $Y_1$ uniform in $\Mat_2(\OO)$, and the
conditional probability is $\frac12$ by \zcref{prop:gl2}. Hence
\[
\phi_1=\frac1q\Bigl(\frac{q-1}q\,\rho_{\pp}+\frac1{2q}\Bigr)
=\frac1q\cdot\frac{(q-1)(2q+1)+(q+1)}{2q(q+1)}=\frac1{q+1}.
\]

For $\phi_2$ let $\mathfrak A=\Mat_2(\OO)$ and
$\mathfrak A'=g\mathfrak Ag^{-1}
=\begin{psmallmatrix}\OO&\pp^{-1}\\\pp&\OO\end{psmallmatrix}$, where
$g=\diag(1,\varpi)$. Then $\Lat_1=\mathfrak A\cap\mathfrak A'$ and
$\Lat_2=\mathfrak A\cup\mathfrak A'\cup E$, where $E$ is the set of
$Y=\begin{psmallmatrix}a&b\\c&d\end{psmallmatrix}\in\Lat_2$ with
$\val(b)=-1$ and $c\in\OO^\times$. With respect to the probability Haar
measure on $\Lat_2$, the sets $\mathfrak A$ and $\mathfrak A'$ have measure
$q^{-1}$, their intersection $\Lat_1$ has measure $q^{-2}$, and $E$ has
measure $(1-q^{-1})^2$.

Every $Y\in E$ is elliptic since $\det Y=ad-bc$ has valuation $-1$ while $\tr Y\in\OO$.
On $\mathfrak A$ the conditional probability of irreducibility is $\frac12$
by \zcref{prop:gl2}. The same holds on $\mathfrak A'$, since conjugation by
$g$ preserves the Haar measure of $\Mat_2(F)$ and the characteristic
polynomial. On $\Lat_1$ the conditional probability is $\phi_1$. By
inclusion--exclusion,
\[
\phi_2=\Bigl(1-\frac1q\Bigr)^2+\frac1{2q}+\frac1{2q}-\frac{\phi_1}{q^2}
=\frac{(q-1)^2(q+1)+q(q+1)-1}{q^2(q+1)}=\frac q{q+1}.\qedhere
\]
\end{proof}

\subsection{Levi densities for \texorpdfstring{$\GL_d$}{GLd}}\label{sec:gldconj}
Let $\bG=\GL_d$ and $\Lat=\Mat_d(\OO)$, so that $\dLat=\Mat_d(\OO)$ under the
trace pairing. For a partition $\mathbf d=(d_1,\dots,d_k)$ of $d$ let
$[M_{\mathbf d}]$ be the class of $\GL_{d_1}\times\dots\times\GL_{d_k}$, and
put $z_{\mathbf d}=\prod_{i\ge1}i^{m_i}m_i!$, where $m_i$ is the number of
parts of $\mathbf d$ which are equal to $i$. Thus $1/z_{\mathbf d}$ is the
proportion of the permutations of cycle type $\mathbf d$ in the symmetric
group on $d$ letters. By \zcref{ex:gldlevi}, the density
$\levid{[M_{\mathbf d}]}{\Mat_d(\OO)}$ is the probability that the
characteristic polynomial of a random matrix in $\Mat_d(\OO)$ is a product of
irreducible polynomials of degrees $d_1,\dots,d_k$, and by \zcref{thm:levi}
it is the limit of the proportion of the Plancherel mass of the fixed vectors
of the principal congruence subgroups of $\GL_d(\OO)$ which is carried by the
tempered representations whose supercuspidal support, or discrete series
support, lies in $[M_{\mathbf d}]$.

\begin{conjecture}[Universality]\label{conj:gld}
For every partition $\mathbf d$ of $d$ there exist monic palindromic
polynomials $N_{\mathbf d},D_{\mathbf d}\in\Z[x]$ of the same degree, where
$D_{\mathbf d}$ is a product of cyclotomic polynomials other than $x-1$, such
that for every non-archimedean local field $F$, with residue field of size $q$,
\[
\levid{[M_{\mathbf d}]}{\Mat_d(\OO)}=R_{\mathbf d}(q),\qquad\text{where }
R_{\mathbf d}=\frac{N_{\mathbf d}}{z_{\mathbf d}D_{\mathbf d}}.
\]
In particular, $R_{\mathbf d}(1/x)=R_{\mathbf d}(x)$ and
$R_{\mathbf d}(q)\to1/z_{\mathbf d}$ as $q\to\infty$.
\end{conjecture}

Here a polynomial $P$ of degree $e$ is called palindromic if $x^eP(1/x)=P(x)$.

\begin{remark}[Known cases and related work]\label{rem:gldconj}
\begin{enumerate}[label=\textup{(\arabic*)},ref=\arabic*]
\item For $d=2$ the conjecture holds, with
$R_{(2)}=R_{(1,1)}=\frac12$, by \zcref{prop:gl2}.
\item\label{item:gldconj-split}
The conjecture holds for the completely split partition $(1^d)$.
Write $R_d(q)$ for its density, and put
$Z_d(q)=\int_{\OO^d}\prod_{i<j}\abs{x_i-x_j}\,dx$, with probability Haar
measure. For every $F$,
\[
d!\,R_d(q)=\prod_{i=1}^d\frac{1-q^{-i}}{1-q^{-1}}\,Z_d(q).
\]
The argument in the proof of \cite[Theorem~12.1]{SVP26} works over
every non-archimedean local field $F$. The integral $Z_d(q)$ also has an
explicit rooted-tree formula \cite[Corollary~5.7]{MR4842997}; see
\cite{MR2229382} for a tree formula for the monic splitting probability.
The asserted form follows by induction from the residue-class recurrence
for $Z_d(q)$, together with the reciprocity and limit at infinity for
$Z_d(q)=d!\,r_d(q)$ from \cite[\S4]{MR2229382}, where $r_d(q)$ is the
splitting probability for random monic polynomials.
\item Let $F=\Q_p$ and $d=3$. By \cite[Corollary~9.2]{SVP26} the expected
number of eigenvalues in $\Q_p$ of a random matrix in $\Mat_d(\Z_p)$ is $1$.
Together with \zcref{item:gldconj-split} this gives
\[
R_{(1,1,1)}=\frac{x^4+x^2+1}{6\,\Phi_5},\qquad
R_{(2,1)}=\frac{x^4+2x^3+x^2+2x+1}{2\,\Phi_5},\qquad
R_{(3)}=\frac{x^4+x^2+1}{3\,\Phi_5},
\]
where $\Phi_5=x^4+x^3+x^2+x+1$.
\item \zcref{conj:gld} is the analogue for random matrices of the conjecture
of Bhargava, Cremona, Fisher and Gajovi\'c \cite[Conjecture~1.2]{MR4422614}
on the densities of the polynomials over $\Z_p$ with a given factorization
type, which was proved, for the finite extensions of $\Q_p$, by G, Wei and Yin
\cite[Theorem~1.4]{GWY25}. The two distributions on characteristic
polynomials are different, already for $d=2$. The invariance under
$p\mapsto1/p$ of several statistics of the eigenvalues of random matrices over
$\Z_p$ was observed in \cite[\S1.6]{SVP26}.
\item We do not know in general whether the density
$\levid{[M_{\mathbf d}]}{\Mat_d(\OO)}$ depends only on $q$.
\end{enumerate}
\end{remark}

\subsection{The Plancherel mass of \texorpdfstring{$\GL_2$}{GL2}}\label{sec:gl2exact}

We compute the nonsupercuspidal mass directly from the Plancherel
formula and subtract it from the total mass. The calculation applies
to every non-archimedean local field, including those of residual
characteristic $2$.
We freely use the standard results on the representations of $\GL_2(F)$, for
which we refer to \cite{MR2234120} and \cite{MR1431508}; for the explicit
Plancherel measure see \cite[Theorem~VIII.1.1]{MR1989693} and
\cite[Theorems~4.3 and~4.4]{MR2131140}.

Let $G=\GL_2(F)$, $K=\GL_2(\OO)$, let $\Iw\subset K$ be the Iwahori
subgroup of matrices which are upper triangular modulo $\pp$, and
let $K_n=1+\pp^n\Mat_2(\OO)$, $n\ge1$; this is the tower
\eqref{eq:tower} for $\kappa(g)=g-1$ and $\Lat=\Mat_2(\OO)$. We
normalize $dg$ by $\vol(K)=1$, so that
\begin{equation}\label{eq:volumes}
\vol(\Iw)=\frac1{q+1},\qquad
\vol(K_n)=\frac1{\abs{\GL_2(k)}\,q^{4(n-1)}}
=\frac1{q(q-1)^2(q+1)q^{4(n-1)}}.
\end{equation}
Let $B=TU$ be the Borel subgroup of upper triangular matrices and
$\delta_B$ its modulus character. For a unitary character $\chi$ of
$F^\times$ let $\cond(\chi)$ be the smallest $m\ge0$ such that $\chi$ is
trivial on $1+\pp^m$ (with $1+\pp^0=\OO^\times$). For
unitary characters $\chi_1,\chi_2$ let
$\pi(\chi_1,\chi_2)=i_B^G(\chi_1\otimes\chi_2)$ be
the unitary principal series.
The tempered, irreducible representations of $G$ are the representations $\pi(\chi_1,\chi_2)$,
with $\pi(\chi_1,\chi_2)\simeq\pi(\chi_2,\chi_1)$, the unitary twists $\St\otimes\chi=\St\otimes(\chi\circ\det)$ of the Steinberg representation $\St$,
and the unitary supercuspidal representations.

Let $\Xc$ be the group of unitary characters of $F^\times$ and
$\Xc_m=\{\chi\in\Xc:\cond(\chi)\le m\}$, $m\ge0$. These are open
subgroups; $\Xc_0$ is the circle of unramified characters, and
$\Xc_m/\Xc_0$ is the dual of $\OO^\times/(1+\pp^m)$. We fix the
Haar measure $d\chi$ on $\Xc$ with $\int_{\Xc_0}d\chi=1$, so that
\begin{equation}\label{eq:massXm}
\int_{\Xc_m}d\chi=[\OO^\times:1+\pp^m]=(q-1)q^{m-1}\quad(m\ge1).
\end{equation}

For unitary characters $\chi,\chi_1,\chi_2$ of $F^\times$ and $n\ge1$ we have
\begin{equation}\label{eq:fixed}
\begin{aligned}
\dim\pi(\chi_1,\chi_2)^{K_n}&=\begin{cases}
q^{n-1}(q+1)&\text{if $\cond(\chi_1),\cond(\chi_2)\le n$,}\\
0&\text{otherwise,}\end{cases}\\
\dim(\St\otimes\chi)^{K_n}&=\begin{cases}
q^{n-1}(q+1)-1&\text{if $\cond(\chi)\le n$,}\\
0&\text{otherwise,}\end{cases}\\
\dim\pi(\chi_1,\chi_2)^{\Iw}&=\begin{cases}
2&\text{if $\chi_1$ and $\chi_2$ are unramified,}\\
0&\text{otherwise,}\end{cases}\\
\dim(\St\otimes\chi)^{\Iw}&=\begin{cases}
1&\text{if $\chi$ is unramified,}\\
0&\text{otherwise,}\end{cases}\\
\dim\pi(\chi_1,\chi_2)^{K}&=\begin{cases}
1&\text{if $\chi_1$ and $\chi_2$ are unramified,}\\
0&\text{otherwise,}\end{cases}
\end{aligned}
\end{equation}
while $(\St\otimes\chi)^K=0$, and $\pi^K=\pi^{\Iw}=0$ for every supercuspidal
representation $\pi$.

The Plancherel formula reads
\begin{equation}\label{eq:plancherelgl2}
\begin{aligned}
(\varphi/dg)(1)={}&\int_{\Gsc}\tr\pi(\varphi)\,d\mu_G(\pi)
+\Cst\int_{\Xc}\tr(\St\otimes\chi)(\varphi)\,d\chi\\
&+\frac{\Cps}{2}\int_{\Xc\times\Xc}\mu(\chi_1\chi_2^{-1})
\tr\pi(\chi_1,\chi_2)(\varphi)\,d\chi_1\,d\chi_2,\ \ \ \varphi\in\Hk(G),
\end{aligned}
\end{equation}
for suitable constants $\Cps,\Cst>0$,
where $\mu$ is Harish-Chandra's $\mu$-function of the unitary principal series, normalized so that
\begin{equation}\label{eq:mu}
\mu(\chi)=q^{\cond(\chi)}\,
\frac{L(1,\chi)L(1,\chi^{-1})}{L(0,\chi)L(0,\chi^{-1})}
=\begin{cases}
\dfrac{\abs{1-z}^2}{\abs{1-q^{-1}z}^2}
&\text{if $\cond(\chi)=0$ and $z=\chi(\varpi)$,}\\[1.5ex]
q^{\cond(\chi)}&\text{if $\cond(\chi)\ge1$.}
\end{cases}
\end{equation}

To determine $\Cps$ we substitute $\varphi=e_K$ and use
\begin{equation}\label{eq:I0}
\int_{\Xc_0}\mu(\chi)\,d\chi
=\frac{1}{2\pi}\int_0^{2\pi}\frac{\abs{1-e^{i\theta}}^2}{\abs{1-q^{-1}e^{i\theta}}^2}\,d\theta
=\frac{2q}{q+1}.
\end{equation}
By \eqref{eq:fixed} we get
\[
1=\frac{\Cps}{2}\int_{\Xc_0\times\Xc_0}\mu(\chi_1\chi_2^{-1})\,d\chi_1\,d\chi_2
=\frac{\Cps}{2}\int_{\Xc_0}\mu(\chi)\,d\chi=\frac{q}{q+1}\,\Cps,
\]
i.e., $\Cps=(q+1)/q$.
Similarly, taking $\varphi=e_{\Iw}$ gives
\[
q+1=\frac{\Cps}{2}\cdot2\int_{\Xc_0}\mu(\chi)\,d\chi+\Cst\int_{\Xc_0}d\chi
=2+\Cst,
\]
i.e., $\Cst=q-1$.

Taking $\varphi=e_{K_n}$ in \eqref{eq:plancherelgl2}, and using
$\tr\pi(e_{K_n})=\dim\pi^{K_n}$ and $(e_{K_n}/dg)(1)=\vol(K_n)^{-1}$, we get the decomposition
\begin{equation}\label{eq:threeparts}
1=\scd{K_n}+\Pst+\Pps
\end{equation}
of the normalized Plancherel mass of the $K_n$-fixed vectors into its
supercuspidal, Steinberg and principal series parts, where
\begin{align*}
\Pst&=\vol(K_n)\Cst\int_{\Xc}\dim(\St\otimes\chi)^{K_n}\,d\chi,\\
\Pps&=\vol(K_n)\frac{\Cps}{2}\int_{\Xc\times\Xc}\mu(\chi_1\chi_2^{-1})
\dim\pi(\chi_1,\chi_2)^{K_n}\,d\chi_1\,d\chi_2.
\end{align*}

In order to compute these masses consider
\[
\int_{\Xc_m}\mu(\chi)\,d\chi,\ \ m\ge1.
\]
By \eqref{eq:mu} and \eqref{eq:massXm}, the set
$\Xc_1\setminus\Xc_0$ has measure $q-2$ and contributes $q(q-2)$, and
for $m\ge2$ the set $\Xc_m\setminus\Xc_{m-1}$ has measure
$(q-1)^2q^{m-2}$ and contributes $(q-1)^2q^{2m-2}$, which is
$\frac{(q-1)q^{2m}}{q+1}-\frac{(q-1)q^{2m-2}}{q+1}$.
Using \eqref{eq:I0} and induction on $m$ we infer that
\begin{equation}\label{eq:Sm}
\int_{\Xc_m}\mu(\chi)\,d\chi=\frac{(q-1)q^{2m}}{q+1},\ \ m\ge1.
\end{equation}

By \eqref{eq:fixed}, \eqref{eq:massXm}, \eqref{eq:volumes}, we infer
\begin{align*}
\Pps&=\vol(K_n)\frac{\Cps}{2}q^{n-1}(q+1)\int_{\Xc_n\times\Xc_n}\mu(\chi_1\chi_2^{-1})\,d\chi_1\,d\chi_2\\
&=\vol(K_n)\frac{\Cps}2q^{n-1}(q+1)\cdot(q-1)q^{n-1}\cdot\frac{(q-1)q^{2n}}{q+1}=\frac12,
\end{align*}
and
\[
\Pst=\vol(K_n)(q-1)\bigl(q^{n-1}(q+1)-1\bigr)(q-1)q^{n-1}=\frac{q^{n-1}(q+1)-1}{q(q+1)q^{3n-3}}.
\]
Taking into account \eqref{eq:threeparts} we conclude

\begin{proposition}[Exact Plancherel masses for $\GL_2$]\label{prop:gl2exact}
For every non-archimedean local field $F$ and every $n\ge1$,
\[
\Pps=\frac12,\qquad
\Pst=\frac{q^{n-1}(q+1)-1}{q(q+1)q^{3n-3}},\qquad
\scd{K_n}=\frac12-q^{1-2n}+\frac1{(q+1)q^{3n-2}}.
\]
In particular $\lim_n\scd{K_n}=\frac12$, and
$q^{-2n}<\frac12-\scd{K_n}<q^{1-2n}$ for every $n\ge1$. Moreover,
$\dsd{K_n}=\frac12$ for every $n\ge1$.
\end{proposition}

For $n=1$ this gives $\scd{K_1}=\frac{q-1}{2(q+1)}$, in agreement with the
$\frac12q(q-1)$ cuspidal representations of $\GL_2(k)$, each of dimension $q-1$.

\subsection{\texorpdfstring{$\PGL_2$}{PGL2} in characteristic two} \label{ex:pgl2}
Let $\chr F=2$ and $\bG=\PGL_2$. In $\mathfrak{pgl}_2=\mathfrak{gl}_2/F$
we have $[E_{12},E_{21}]=E_{11}-E_{22}\equiv0$, since
$E_{11}+E_{22}$ is the identity. Thus the pairing
\eqref{eq:levipairing} vanishes identically, the incidence map of the
Borel subgroup is everywhere critical, and $\elld[\PGL_2]{\dLat}=1$ for
every $\dLat$ by \zcref{lem:one}. Thus \zcref{rem:lower} gives
$1-\scd{K_n}\ll q^{-n/4}$.

One can also see the full elliptic density directly. The trace pairing identifies
$\mathfrak{pgl}_2^*$ with $\mathfrak{sl}_2=\{\begin{psmallmatrix}a&b\\c&a\end{psmallmatrix}\}$,
and such a matrix is non-elliptic if and only if it has an $F$-rational
eigenline, i.e., if and only if $bc$ is a square in $F$. The subgroup
$F^2=\{x^2:x\in F\}$ is a closed additive subgroup of $F$ of measure zero.
The same phenomenon occurs for the \oddorth{} parts $\Godd$ of all groups (\zcref{lem:split}).

In terms of \zcref{thm:levi}, $\Lcl$ consists of $[G]$ and of the class $[T]$
of the split torus. Since the incidence map of the Borel subgroup is
everywhere critical, $\Estr_{[T]}=\Pim{[T]}=\Par{\gd}$ is contained in the null set $\Crv$, so that
$\levid{[T]}{\dLat}=0$, while $\Estr_{[G]}=\Ell{\gd}$ and $\levid{[G]}{\dLat}=1$,
for every lattice $\dLat$. Thus the tempered representations with
supercuspidal or discrete series support in the torus, i.e., the
constituents of the unitary principal series and, for the supercuspidal
support, also the twists of the Steinberg representation, carry asymptotically no
Plancherel mass, in accordance with \zcref{prop:pgl2}. This shows that the limits in \zcref{thm:levi} need
not be positive, and that the toral description of \zcref{rem:toral} does not
extend to all characteristics.

On the spectral side, squaring changes the distribution of conductors:
many characters $\chi$ have $\chi^2$ of smaller conductor. The
following proposition computes this distribution and the resulting
principal-series mass explicitly.

\begin{proposition}[Exact masses for $\PGL_2$ in characteristic two]\label{prop:pgl2}
Let $\chr F=2$, $G=\PGL_2(F)$, and let
$K_n=\ker(\PGL_2(\OO)\to\PGL_2(\OO/\pp^n))$, $n\ge1$, be the tower
\eqref{eq:tower} of the smooth model $\PGL_2$ over $\OO$
(\zcref{lem:model}). Then, with $\ell=\lfloor n/2\rfloor$,
\[
\Pps=\frac{q^{-\ell}}2,\qquad
\Pst=\frac{q^{\ell}\bigl(q^{n-1}(q+1)-1\bigr)}{q(q+1)q^{3n-3}},\qquad
1-\scd{K_n}=\frac{q^{-\ell}}2\bigl(1+O(q^{-n})\bigr).
\]
Moreover, $\dsd{K_n}=1-\frac{q^{-\ell}}2$.
\end{proposition}

\begin{proof}
The representations of $G$ are the representations of $\GL_2(F)$ with
trivial central character, and $K_n$ is the image of $1+\pp^n\Mat_2(\OO)$,
so by the classification recalled above and \eqref{eq:fixed}, $\Gtemp$
consists of $\Gsc$, of the twists $\St\otimes\chi$ with $\chi^2=1$, and
of the principal series $\pi(\chi,\chi^{-1})$, $\chi\in\Xc$ modulo
$\chi\mapsto\chi^{-1}$, with $K_n$-fixed vectors of dimensions
$q^{n-1}(q+1)-1$, resp.\ $q^{n-1}(q+1)$, if $\cond(\chi)\le n$ and $0$
otherwise. Put $K=\PGL_2(\OO)$ and normalize $dg$ by $\vol(K)=1$.
Then
$\vol(\Iw)=(q+1)^{-1}$ and
$\vol(K_n)^{-1}=\abs{\PGL_2(\OO/\pp^n)}=q(q^2-1)q^{3(n-1)}$.
The Plancherel formula becomes
\[
\begin{aligned}
(\varphi/dg)(1)={}&\sum_{\pi\in\Gsc}\mu_G(\{\pi\})\tr\pi(\varphi)
+\mu_G(\{\St\})\sum_{\chi^2=1}\tr(\St\otimes\chi)(\varphi)\\
&+\frac{\Cps'}{2}\int_{\Xc}\mu(\chi^2)\tr\pi(\chi,\chi^{-1})(\varphi)\,d\chi,
\end{aligned}
\]
where $\mu$ is as in \eqref{eq:mu}. Here, for a discrete series
representation $\pi$, the mass $\mu_G(\{\pi\})$ is the formal degree of $\pi$
with respect to $dg$, i.e., $\fdeg{\pi}=\mu_G(\{\pi\})\,dg$; it is the same
for $\St$ and its twists. The constant $\Cps'$ refers
to the present normalization; a priori it could have a different value from $\Cps$ of \zcref{sec:gl2exact}.
The constants are determined as before: $\varphi=e_K$ gives
\[
1=\frac{\Cps'}{2}\int_{\Xc_0}\mu(\chi^2)\,d\chi=\frac{\Cps'}{2}\int_{\Xc_0}\mu(\chi)\,d\chi=\frac{q}{q+1}\Cps',
\]
since $z\mapsto z^2$ preserves the Haar measure of the circle, and
$\varphi=e_{\Iw}$ gives
\[
q+1=2+2\mu_G(\{\St\}),
\]
because exactly two characters with $\chi^2=1$ are unramified ($\chi(\varpi)=\pm1$). Hence
$\Cps'=\frac{q+1}q$ and $\mu_G(\{\St\})=\frac{q-1}2$.

It remains to compute the distribution of $\cond(\chi^2)$ on $\Xc_n$.
Squaring is injective on $\OO^\times$, and $(1+y)^2=1+y^2$, so for
$0\le j\le n$, with the convention $1+\pp^0=\OO^\times$,
\[
(1+\pp^j)^2\cap(1+\pp^n)=(1+\pp^{\max(j,\lceil n/2\rceil)})^2,
\]
and the image of $(1+\pp^j)^2$ in $\OO^\times/(1+\pp^n)$ has
order $[1+\pp^j:1+\pp^{\lceil n/2\rceil}]=q^{\lceil n/2\rceil-j}$
for $1\le j\le\lceil n/2\rceil$, and order
$(q-1)q^{\lceil n/2\rceil-1}$ for $j=0$.

A character $\chi\in\Xc_n$ satisfies $\cond(\chi^2)\le j$ if and
only if $\chi$ is trivial on $(1+\pp^j)^2$, and the characters of
$\OO^\times/(1+\pp^n)$ trivial on a subgroup of order $N$ form a
set of cardinality $(q-1)q^{n-1}/N$. Writing $\Xc_n=\Xc_n^{\circ}\cdot\Xc_0$
with $\Xc_n^\circ$ the characters trivial on $\varpi$, so that
$\cond(\chi^2)$ depends only on the $\Xc_n^\circ$-component, we get for
the number $N_j$ of $\xi\in\Xc_n^\circ$ with $\cond(\xi^2)\le j$:
\[
N_0=q^\ell,\qquad N_j=(q-1)q^{\ell+j-1}\quad(1\le j\le\lceil n/2\rceil),
\]
and $N_{\lceil n/2\rceil}=(q-1)q^{n-1}=\#\Xc_n^\circ$. Hence, by
\eqref{eq:mu}, \eqref{eq:I0} and the invariance of $d\chi$ on $\Xc_0$
under $\chi\mapsto\chi^2$,
\[
\int_{\Xc_n}\mu(\chi^2)\,d\chi
=N_0\frac{2q}{q+1}+\sum_{j=1}^{\lceil n/2\rceil}q^j(N_j-N_{j-1}).
\]
The sum is evaluated as in the proof of \eqref{eq:Sm}: the term $j=1$
is $q^{\ell+1}(q-2)$ and the terms $j\ge2$ are $(q-1)^2q^{\ell+2j-2}$, so
that the right-hand side equals $q^\ell$ times
the integral \eqref{eq:Sm} with $m=\lceil n/2\rceil$, namely $(q-1)q^{\ell+2\lceil n/2\rceil}/(q+1)=(q-1)q^{n+\lceil n/2\rceil}/(q+1)$.
Therefore
\[
\Pps=\vol(K_n)\frac{\Cps'}2q^{n-1}(q+1)\cdot\frac{(q-1)q^{n+\lceil n/2\rceil}}{q+1}
=\frac{q^{2n+\lceil n/2\rceil-1}}{2q^{3n-1}}=\frac{q^{-\ell}}2,
\]
and, since the characters with $\chi^2=1$ and $\cond(\chi)\le n$ number
$2N_0=2q^\ell$,
\[
\Pst=\vol(K_n)\frac{q-1}2\bigl(q^{n-1}(q+1)-1\bigr)2q^\ell
=\frac{q^\ell\bigl(q^{n-1}(q+1)-1\bigr)}{q(q+1)q^{3n-3}}
\le q^{\ell-2n+1}\le 2q^{1-n}\cdot\frac{q^{-\ell}}2,
\]
since $2\ell\le n$. This gives the formula for $1-\scd{K_n}$. Finally,
the discrete series of $G$ are the supercuspidal representations and the
twists $\St\otimes\chi$ with $\chi^2=1$, so $\dsd{K_n}=1-\Pps$.
\end{proof}

The $2q^\ell$ twists $\St\otimes\chi$ with $K_n$-fixed vectors come from
quadratic characters $\chi$ of unbounded conductor. Thus, although $\PGL_2$
has trivial center, its non-supercuspidal discrete series do not have
bounded depth. In the general comparison, their mass is supported on the exceptional cells.

%\bibliographystyle{../../Bibfiles/amsalphapure}
%\bibliography{../../Bibfiles/all}

\begin{thebibliography}{DdSMS99}

\bibitem[AP05]{MR2131140}
Anne-Marie Aubert and Roger Plymen, \emph{Plancherel measure for {${\rm
  GL}(n,F)$} and {${\rm GL}(m,D)$}: explicit formulas and {B}ernstein
  decomposition}, J. Number Theory \textbf{112} (2005), no.~1, 26--66.
  \MR{2131140}

\bibitem[BCFG22]{MR4422614}
Manjul Bhargava, John Cremona, Tom Fisher, and Stevan Gajovi\'{c}, \emph{The
  density of polynomials of degree {$n$} over {$\Bbb Z_p$} having exactly {$r$}
  roots in {$\Bbb Q_p$}}, Proc. Lond. Math. Soc. (3) \textbf{124} (2022),
  no.~5, 713--736. \MR{4422614}

\bibitem[BGMR06]{MR2229382}
Joe Buhler, Daniel Goldstein, David Moews, and Joel Rosenberg, \emph{The
  probability that a random monic {$p$}-adic polynomial splits}, Experiment.
  Math. \textbf{15} (2006), no.~1, 21--32. \MR{2229382}

\bibitem[BH06]{MR2234120}
Colin~J. Bushnell and Guy Henniart, \emph{The local {L}anglands conjecture for
  {$\rm GL(2)$}}, Grundlehren der Mathematischen Wissenschaften [Fundamental
  Principles of Mathematical Sciences], vol. 335, Springer-Verlag, Berlin,
  2006. \MR{2234120}

\bibitem[BHK11]{MR2767510}
Colin~J. Bushnell, Guy Henniart, and Philip~C. Kutzko, \emph{Types and explicit
  {P}lancherel formul\ae \ for reductive {$p$}-adic groups}, On certain
  {$L$}-functions, Clay Math. Proc., vol.~13, Amer. Math. Soc., Providence, RI,
  2011, pp.~55--80. \MR{2767510}

\bibitem[BK98]{MR1643417}
Colin~J. Bushnell and Philip~C. Kutzko, \emph{Smooth representations of
  reductive {$p$}-adic groups: structure theory via types}, Proc. London Math.
  Soc. (3) \textbf{77} (1998), no.~3, 582--634. \MR{1643417}

\bibitem[BM83]{MR737927}
Walter Borho and Robert MacPherson, \emph{Partial resolutions of nilpotent
  varieties}, Analysis and topology on singular spaces, {II}, {III} ({L}uminy,
  1981), Ast\'{e}risque, vol. 101, Soc. Math. France, Paris, 1983, pp.~23--74.
  \MR{737927}

\bibitem[Bor76]{MR0444849}
Armand Borel, \emph{Admissible representations of a semi-simple group over a
  local field with vectors fixed under an {I}wahori subgroup}, Invent. Math.
  \textbf{35} (1976), 233--259. \MR{MR0444849 (56 \#3196)}

\bibitem[Bor91]{MR1102012}
\bysame, \emph{Linear algebraic groups}, second ed., Graduate Texts in
  Mathematics, vol. 126, Springer-Verlag, New York, 1991. \MR{MR1102012
  (92d:20001)}

\bibitem[Bou98]{MR1728312}
Nicolas Bourbaki, \emph{Lie groups and {L}ie algebras. {C}hapters 1--3},
  Elements of Mathematics (Berlin), Springer-Verlag, Berlin, 1998, Translated
  from the French, Reprint of the 1989 English translation. \MR{MR1728312
  (2001g:17006)}

\bibitem[BP16]{MR3463423}
Rapha\"el Beuzart-Plessis, \emph{A short proof of the existence of
  supercuspidal representations for all reductive {$p$}-adic groups}, Pacific
  J. Math. \textbf{282} (2016), no.~1, 27--34. \MR{3463423}

\bibitem[BT65]{MR0207712}
Armand Borel and Jacques Tits, \emph{Groupes r\'eductifs}, Inst. Hautes
  \'Etudes Sci. Publ. Math. (1965), no.~27, 55--150. \MR{MR0207712 (34 \#7527)}

\bibitem[Bum97]{MR1431508}
Daniel Bump, \emph{Automorphic forms and representations}, Cambridge Studies in
  Advanced Mathematics, vol.~55, Cambridge University Press, Cambridge, 1997.
  \MR{1431508}

\bibitem[BZ76]{MR0425030}
I.~N. Bern\v{s}te\u{\i}n and A.~V. Zelevinski\u{\i}, \emph{Representations of
  the group {$GL(n,F),$} where {$F$} is a local non-{A}rchimedean field},
  Uspehi Mat. Nauk \textbf{31} (1976), no.~3(189), 5--70. \MR{0425030}

\bibitem[BZ77]{MR0579172}
I.~N. Bernstein and A.~V. Zelevinsky, \emph{Induced representations of
  reductive {${\germ p}$}-adic groups. {I}}, Ann. Sci. \'Ecole Norm. Sup. (4)
  \textbf{10} (1977), no.~4, 441--472. \MR{0579172 (58 \#28310)}

\bibitem[Car79]{MR522009}
Henri Carayol, \emph{Repr\'esentations supercuspidales de {${\rm GL}\sb{n}$}},
  C. R. Acad. Sci. Paris S\'er. A-B \textbf{288} (1979), no.~1, A17--A20.
  \MR{MR522009 (80b:22023)}

\bibitem[Car84]{MR760676}
H.~Carayol, \emph{Repr\'esentations cuspidales du groupe lin\'eaire}, Ann. Sci.
  \'Ecole Norm. Sup. (4) \textbf{17} (1984), no.~2, 191--225. \MR{MR760676
  (86f:22019)}

\bibitem[CG97]{MR1433132}
Neil Chriss and Victor Ginzburg, \emph{Representation theory and complex
  geometry}, Birkh\"auser Boston, Inc., Boston, MA, 1997. \MR{1433132
  (98i:22021)}

\bibitem[CGP15]{MR3362817}
Brian Conrad, Ofer Gabber, and Gopal Prasad, \emph{Pseudo-reductive groups},
  second ed., New Mathematical Monographs, vol.~26, Cambridge University Press,
  Cambridge, 2015. \MR{3362817}

\bibitem[Chi26]{Chi26}
Jingren Chi, \emph{Equivalued affine {S}pringer fibers in mixed
  characteristic}, 2026, arXiv:2606.31083, version 2.

\bibitem[Con14]{MR3362641}
Brian Conrad, \emph{Reductive group schemes}, Autour des sch\'emas en groupes.
  {V}ol. {I}, Panor. Synth\`eses, vol. 42/43, Soc. Math. France, Paris, 2014,
  pp.~93--444. \MR{3362641}

\bibitem[DdSMS99]{MR1720368}
J.~D. Dixon, M.~P.~F. du~Sautoy, A.~Mann, and D.~Segal, \emph{Analytic
  pro-{$p$} groups}, second ed., Cambridge Studies in Advanced Mathematics,
  vol.~61, Cambridge University Press, Cambridge, 1999. \MR{1720368
  (2000m:20039)}

\bibitem[DeB06]{MR2214792}
Stephen DeBacker, \emph{Parameterizing conjugacy classes of maximal unramified
  tori via {B}ruhat-{T}its theory}, Michigan Math. J. \textbf{54} (2006),
  no.~1, 157--178. \MR{2214792}

\bibitem[Dil23]{MR4617942}
Peter Dillery, \emph{Rigid inner forms over local function fields}, Adv. Math.
  \textbf{430} (2023), Paper No. 109204, 100. \MR{4617942}

\bibitem[FS21]{MR4325284}
Jessica Fintzen and Sug~Woo Shin, \emph{Congruences of algebraic automorphic
  forms and supercuspidal representations}, Camb. J. Math. \textbf{9} (2021),
  no.~2, 351--429. \MR{4325284}

\bibitem[Ger75]{MR396859}
Paul G\'erardin, \emph{Construction de s\'eries discr\`etes {$p$}-adiques},
  Lecture Notes in Mathematics, vol. Vol. 462, Springer-Verlag, Berlin-New
  York, 1975, Sur les s\'eries discr\`etes non ramifi\'ees des groupes
  r\'eductifs d\'eploy\'es $p$-adiques. \MR{396859}

\bibitem[GH26]{GH26}
Itay Glazer and Yotam~I. Hendel, \emph{A lower bound on the analytic
  log-canonical threshold over local fields of positive characteristic}, Proc.
  Amer. Math. Soc. (2026), Electronically published September 30, 2026;
  arXiv:2511.01270.

\bibitem[GKM06]{MR2209851}
Mark Goresky, Robert Kottwitz, and Robert MacPherson, \emph{Purity of
  equivalued affine {S}pringer fibers}, Represent. Theory \textbf{10} (2006),
  130--146 (electronic). \MR{2209851 (2007i:22025)}

\bibitem[Glo13]{MR3057121}
Helge Gl\"{o}ckner, \emph{Invariant manifolds for analytic dynamical systems
  over ultrametric fields}, Expo. Math. \textbf{31} (2013), no.~2, 116--150.
  \MR{3057121}

\bibitem[Gre66]{MR0207700}
Marvin~J. Greenberg, \emph{Rational points in {H}enselian discrete valuation
  rings}, Inst. Hautes \'Etudes Sci. Publ. Math. (1966), no.~31, 59--64.
  \MR{0207700 (34 \#7515)}

\bibitem[GWY25]{GWY25}
Asvin G., Yifan Wei, and John Yin, \emph{A {C}hebotarev density theorem over
  local fields}, 2025, arXiv:2212.00294v5.

\bibitem[How77a]{MR0579176}
Roger~E. Howe, \emph{Kirillov theory for compact {$p$}-adic groups}, Pacific J.
  Math. \textbf{73} (1977), no.~2, 365--381. \MR{0579176 (58 \#28314)}

\bibitem[How77b]{MR0492088}
\bysame, \emph{Some qualitative results on the representation theory of {${\rm
  Gl}_{n}$} over a {$p$}-adic field}, Pacific J. Math. \textbf{73} (1977),
  no.~2, 479--538. \MR{0492088 (58 \#11242)}

\bibitem[Igu00]{MR1743467}
Jun-ichi Igusa, \emph{An introduction to the theory of local zeta functions},
  AMS/IP Studies in Advanced Mathematics, vol.~14, American Mathematical
  Society, Providence, RI, 2000. \MR{1743467 (2001j:11112)}

\bibitem[Kim04]{MR2075916}
Ju-Lee Kim, \emph{Dual blobs and {P}lancherel formulas}, Bull. Soc. Math.
  France \textbf{132} (2004), no.~1, 55--80. \MR{2075916 (2005f:22025)}

\bibitem[Kim07]{MR2276772}
\bysame, \emph{Supercuspidal representations: an exhaustion theorem}, J. Amer.
  Math. Soc. \textbf{20} (2007), no.~2, 273--320 (electronic). \MR{2276772
  (2008c:22014)}

\bibitem[KP23]{MR4520154}
Tasho Kaletha and Gopal Prasad, \emph{Bruhat-{T}its theory---a new approach},
  New Mathematical Monographs, vol.~44, Cambridge University Press, Cambridge,
  2023. \MR{4520154}

\bibitem[Mil17]{MR3729270}
J.~S. Milne, \emph{Algebraic groups}, Cambridge Studies in Advanced
  Mathematics, vol. 170, Cambridge University Press, Cambridge, 2017, The
  theory of group schemes of finite type over a field. \MR{3729270}

\bibitem[Mor88]{MR950594}
Lawrence Morris, \emph{{$P$}-cuspidal representations}, Proc. London Math. Soc.
  (3) \textbf{57} (1988), no.~2, 329--356. \MR{950594 (89j:22038)}

\bibitem[Moy17]{MR3753921}
Allen Moy, \emph{Computations with {B}ernstein projectors of {${\rm SL}(2)$}},
  Representation theory, number theory, and invariant theory, Progr. Math.,
  vol. 323, Birkh\"auser/Springer, Cham, 2017, pp.~503--533. \MR{3753921}

\bibitem[MP94]{MR1253198}
Allen Moy and Gopal Prasad, \emph{Unrefined minimal {$K$}-types for {$p$}-adic
  groups}, Invent. Math. \textbf{116} (1994), no.~1-3, 393--408. \MR{MR1253198
  (95f:22023)}

\bibitem[MP96]{MR1371680}
\bysame, \emph{Jacquet functors and unrefined minimal {$K$}-types}, Comment.
  Math. Helv. \textbf{71} (1996), no.~1, 98--121. \MR{MR1371680 (97c:22021)}

\bibitem[MV24]{MR4842997}
Joshua Maglione and Christopher Voll, \emph{Flag {H}ilbert-{P}oincar\'e{}
  series and {I}gusa zeta functions of hyperplane arrangements}, Israel J.
  Math. \textbf{264} (2024), no.~1, 177--233. \MR{4842997}

\bibitem[Oes82]{MR656627}
Joseph Oesterl{\'e}, \emph{R\'eduction modulo {$p^{n}$} des sous-ensembles
  analytiques ferm\'es de {${\bf Z}^{N}_{p}$}}, Invent. Math. \textbf{66}
  (1982), no.~2, 325--341. \MR{656627 (83j:12014)}

\bibitem[PR94]{MR1278263}
Vladimir Platonov and Andrei Rapinchuk, \emph{Algebraic groups and number
  theory}, Pure and Applied Mathematics, vol. 139, Academic Press Inc., Boston,
  MA, 1994, Translated from the 1991 Russian original by Rachel Rowen.
  \MR{1278263 (95b:11039)}

\bibitem[Ren10]{MR2567785}
David Renard, \emph{Repr\'esentations des groupes r\'eductifs {$p$}-adiques},
  Cours Sp\'ecialis\'es [Specialized Courses], vol.~17, Soci\'et\'e
  Math\'ematique de France, Paris, 2010. \MR{2567785 (2011d:22019)}

\bibitem[Ser81]{MR644559}
Jean-Pierre Serre, \emph{Quelques applications du th\'eor\`eme de densit\'e de
  {C}hebotarev}, Inst. Hautes \'Etudes Sci. Publ. Math. (1981), no.~54,
  323--401. \MR{644559 (83k:12011)}

\bibitem[Ser06]{MR2179691}
\bysame, \emph{Lie algebras and {L}ie groups}, Lecture Notes in Mathematics,
  vol. 1500, Springer-Verlag, Berlin, 2006, 1964 lectures given at Harvard
  University, Corrected fifth printing of the second (1992) edition.
  \MR{2179691 (2006e:17001)}

\bibitem[SS70]{MR0268192}
T.~A. Springer and R.~Steinberg, \emph{Conjugacy classes}, Seminar on
  {A}lgebraic {G}roups and {R}elated {F}inite {G}roups ({T}he {I}nstitute for
  {A}dvanced {S}tudy, {P}rinceton, {N}.{J}., 1968/69), Lecture Notes in
  Mathematics, Vol. 131, Springer, Berlin, 1970, pp.~167--266. \MR{0268192 (42
  \#3091)}

\bibitem[ST26]{ST26}
Loren Spice and Cheng-Chiang Tsai, \emph{Jordan decompositions in {L}ie
  algebras and their duals}, 2026, arXiv:2601.07168.

\bibitem[SVP26]{SVP26}
Jiahe Shen and Roger Van~Peski, \emph{Eigenvalues of {$p$}-adic random
  matrices}, 2026, arXiv:2601.06283.

\bibitem[TVVW18]{MR3803728}
David Taylor, V.~S. Varadarajan, Jukka Virtanen, and David Weisbart,
  \emph{Temperedness of measures defined by polynomial equations over local
  fields}, Pacific J. Math. \textbf{296} (2018), no.~1, 227--256. \MR{3803728}

\bibitem[Wal03]{MR1989693}
J.-L. Waldspurger, \emph{La formule de {P}lancherel pour les groupes
  {$p$}-adiques (d'apr\`es {H}arish-{C}handra)}, J. Inst. Math. Jussieu
  \textbf{2} (2003), no.~2, 235--333. \MR{1989693}

\bibitem[Xue09]{MR2558787}
Ting Xue, \emph{Nilpotent orbits in the dual of classical {L}ie algebras in
  characteristic 2 and the {S}pringer correspondence}, Represent. Theory
  \textbf{13} (2009), 609--635. \MR{2558787}

\end{thebibliography}

\def\cprime{$'$} 
\providecommand{\bysame}{\leavevmode\hbox to3em{\hrulefill}\thinspace}
\providecommand{\MR}{\relax\ifhmode\unskip\space\fi MR }
% \MRhref is called by the amsart/book/proc definition of \MR.
\providecommand{\MRhref}[2]{%
  \href{http://www.ams.org/mathscinet-getitem?mr=#1}{#2}
}
\providecommand{\href}[2]{#2}

\end{document}